\documentclass[11pt]{article}

\RequirePackage[T1]{fontenc}
\RequirePackage[french,english]{babel}           %
\RequirePackage{amsmath,amssymb,dsfont}
\RequirePackage[numbers]{natbib}
\RequirePackage[colorlinks,citecolor=blue,urlcolor=blue]{hyperref}

\RequirePackage[utf8]{inputenc}        % pour les accents ;
\RequirePackage{textcomp}                % pour les accents sur les majuscules ;
\RequirePackage{listings,multicol}       % pour la mise en page ;
\RequirePackage{graphicx,color}          % pour inserer des images ; et afficher directement des graphes ;
\RequirePackage{wrapfig}                 % pour habiller une image ;
\RequirePackage[rightcaption]{sidecap}   % pour mettre les legendes des figures a droite des images (avec SCfigure)
\RequirePackage{amsthm}
\RequirePackage{ulem}                    % pour personnaliser les soulignements
\RequirePackage{fancyhdr}                % pour personnaliser les entetes et pieds de pages ;
\RequirePackage{subeqnarray}
\RequirePackage{graphicx}
\RequirePackage{caption}
\RequirePackage{subcaption}

\newcommand{\pa}[1]{\ensuremath{\left( #1 \right)}}
\newcommand{\cro}[1]{\ensuremath{\left[ #1 \right]}}
\newcommand{\ac}[1]{\ensuremath{\left\{ #1 \right\}}}
\newcommand{\abs}[1]{\ensuremath{\left| #1 \right|}}

\newcommand{\bl}{\ $\bullet$ \ }

\newcommand{\Nbc}{\varphi^{coinc}}

\newcommand{\fct}[4]{\ensuremath{\left( \begin{array}{ccc}
                                 #1 & \longrightarrow & #2 \\
                                 #3 &   \longmapsto   & #4 \\
                               \end{array} \right) }}

\newcommand{\tq}{\ /\ }

\newcommand{\N}{\ensuremath{\mathds{N}}}

\newcommand{\R}{\ensuremath{\mathds{R}}}

\newcommand{\X}{\ensuremath{\mathds{X}}}

\renewcommand{\L}{\ensuremath{\mathds{L}}}
\newcommand{\calX}{\ensuremath{\mathcal{X}}}
\newcommand{\calD}{\ensuremath{\mathcal{D}}}
\newcommand{\calL}{\ensuremath{\mathcal{L}}}

\newcommand\independent{\protect\mathpalette{\protect\independenT}{\perp}}
\def\independenT#1#2{\mathrel{\rlap{$#1#2$}\mkern2mu{#1#2}}}

\newcommand{\ds}[1]{\ensuremath{\mathds{#1}}} 			% majuscules rayées en mode mathématique.
\newcommand{\mc}[1]{\ensuremath{\mathcal{#1}}}          % majuscules rondes en mode mathématique.
\newcommand{\proba}[1]{\ensuremath{\mathds{P}\left(#1\right)}}			% probabilite ;
\newcommand{\esp}[1]{\ensuremath{\mathds{E}\left[  #1 \right]}} 		% esperance ;
\newcommand{\Estar}[1]{\ensuremath{\mathds{E}^*\!\!\left[  #1 \right]}} 		% esperance conditionnellement à l'échantillon de départ;
\newcommand{\esps}[2]{\ensuremath{\mathds{E}_{#1}\left[ #2 \right]}}	% esperance sous qqch;

\newcommand{\var}[1]{\ensuremath{\mathds{V}\text{ar}\left(  #1 \right)}} 		% esperance ;

\newcommand{\1}[1]{\ensuremath{\mathds{1}_{ #1 }}}		% fonction caracteristique d'un ensemble ;

\newcommand{\PP}{\ensuremath{\mathds{P}}}

\newcommand{\Nm}{\ensuremath{\mc{N}(0,1)}}              % pour la loi normale ;
\newcommand{\loi}[1]{\ensuremath{\mc{L}\pa{#1}}}

\newcommand{\norm}[1]{\| #1 \|}
\newcommand{\rperm}{\ensuremath{ {\Pi_n} }}

\newcommand{\cv}[1]{\underset{#1}{\longrightarrow}}
\newcommand{\cvf}[1]{\underset{#1}{\Longrightarrow}}
\newcommand{\cvps}[1]{\overset{a.s.}{\underset{#1}{\longrightarrow}}}

\newcommand{\cvproba}[1]{\overset{\PP}{\underset{#1}{\longrightarrow}}}

\newcommand{\cvloi}[1]{\overset{\mc{L}}{\underset{#1}{\longrightarrow}}}

\newcommand{\Sn}[1]{\ensuremath{\mathfrak{S}_{#1}}}

\definecolor{mag}{cmyk}{0,1,0.5,0}

\definecolor{orange}{cmyk}{0,0.5,1,0.3}

\definecolor{vert}{cmyk}{1,0.3,1,0.3}

\numberwithin{equation}{section}
\theoremstyle{plain}
\newtheorem{thm}{Theorem}[section]
\newtheorem{prop}{Proposition}[section]
\newtheorem{lm}{Lemma}[section]
\newtheorem{coro}{Corollary}[section]

\begin{document}

\title{\textsc{Bootstrap and permutation tests of independence for point processes.}}
\author{M\'elisande Albert\thanks{Univ. Nice Sophia Antipolis, CNRS, LJAD, UMR 7351, 06100 Nice, France.},~ Yann Bouret\thanks{Univ. Nice Sophia Antipolis, CNRS, LPMC, UMR 7336, 06100 Nice, France.},~ Magalie Fromont\thanks{Univ. Européenne de Bretagne, CNRS, IRMAR, UMR 6625, 35043 Rennes Cedex, France.} \\
and~ Patricia Reynaud-Bouret\thanks{Univ. Nice Sophia Antipolis, CNRS, LJAD, UMR 7351, 06100 Nice, France.}}
\date{}

\maketitle

\vspace{-18pt}
\small{{\bf Abstract:}~Motivated by a neuroscience question about synchrony detection in spike train analysis, we deal with the independence testing problem for point processes. We introduce non-parametric test statistics, which are rescaled general $U$-statistics, whose corresponding critical values are constructed from bootstrap and randomization/permutation approaches, making as few assumptions as possible on the underlying distribution of the point processes. We derive general consistency results for the bootstrap and for the permutation w.r.t. to Wasserstein's metric, which induce weak convergence as well as convergence of second order moments. The obtained bootstrap or permutation independence tests are thus proved to be asymptotically of the prescribed size, and to be consistent against any reasonable alternative.  A simulation study is performed  to illustrate the derived theoretical results, and to compare the performance of our new tests with existing ones in the neuroscientific literature.}

\medskip

\noindent{\bf Mathematics Subject Classification:} Primary: 62M07, 62F40, 62E20, 60G55, 60F05; secondary: 62P10

\smallskip

\noindent{\bf Keywords:} Independence test, $U$-statistics, point processes, bootstrap, randomization, permutation, neuroscience, spike train analysis.

% ---------------------------------------------------------------------------------------------------- %
\section{Introduction}
\label{Intro}
% ---------------------------------------------------------------------------------------------------- %

Inspired by neuroscience problems, the present work is devoted to independence tests for point processes. The question of testing whether two random variables are independent is of course largely encountered in the statistical literature, as it is one of the central goals of data analysis. From the historical Pearson's chi-square test of independence (see \cite{Pearson1, Pearson2}) to the modern test of \cite{Gretton} using kernel methods in the spirit of statistical learning, many non-parametric independence tests have been developed for real valued random variables or even random vectors. Among them, of particular interest are the tests based on the randomization/permutation principle introduced by Fisher \cite{Fisher}, and covered thereafter in the series of papers by Pitman \cite{Pitman1, Pitman3}, Scheffe \cite{Scheffe}, Hoeffding \cite{Hoeffding52} for instance, or bootstrap approaches derived from Efron's \cite{Efron} "naive" one. 
Note that permutation and bootstrap-based tests have a long history of applications, of which independence tests are just a very small part (see for instance \cite{EfronTib, PesarinSalmaso, Romano88, Romano89} for some reviews, or \cite{AntochHuskova,  JanssenPauls, KlebanovEtal, Kirch, FLLR2012} for more recent works).
Focusing on independence tests, two families of permutation or bootstrap-based tests may be distinguished at least: the whole family of rank tests including the tests of Hotelling and Pabst \cite{Hot}, Kendall \cite{Kendall},  Wolfowitz \cite{Wolfowitz} or Hoeffding \cite{Hoeffding48} on the one hand, the family of Kolmogorov-Smirnov type tests, like Blum, Kiefer, and Rosenblatt's \cite{BKR61}, Romano's \cite{Romano89} or Van der Vaart and Wellner's \cite{vdvwellner96} ones on the other hand. 

To describe the properties of these tests, let us recall and fix a few definitions, which are furthermore  used throughout this article. Tests are said to be {\it non-parametric} if they are free from the underlying distribution of the observed variables. For any prescribed $\alpha$ in $(0,1)$, tests are said to be {\it exactly of level} $\alpha$ if their first kind error rate is less than or equal to $\alpha$ whatever the number of observations.  This is a non-asymptotic property. Tests are also said to be {\it asymptotically of size} $\alpha$ if their first kind error rate tends to $\alpha$ when the number of observations tends to infinity. 
They are said to be {\it consistent} against some alternative if, under this alternative, their second kind error rate tends to $0$ or equivalently their power tends to $1$, when the number of observations tends to infinity. Finally, {\it bootstrap} refers here to bootstrap with replacement. It is thus different from {\it permutation}, which appears sometimes in the literature as bootstrap without replacement.
In this respect, the above mentioned tests of independence are all non-parametric and asymptotically of the prescribed size. Moreover, the tests based on permutation are exactly of the desired level. Some of these tests are proved to be consistent against many alternatives, such as Hoeffding's \cite{Hoeffding48} one and the family of Kolmogorov-Smirnov type tests.

Detecting dependence is also a fundamental old point in the neuroscientific literature (see e.g. \cite{GerPerkel}). The neuroscience problem we were initially interested in consists in detecting interactions between occurrences of action potentials on two different neurons simultaneously recorded on $n$ independent trials, as described in \cite{GrunB}.  Each recorded set of time occurrences of action potentials for each neuron is usually referred to as a spike train, the spikes being the time occurrences themselves. It is commonly accepted that these spikes are one of the main components of the brain activity (see \cite{Singer93}). So, when observing two spike trains coming from two different neurons, one of the main elementary problem is to assess whether these two spike trains are independent or not. Unfortunately, even if the real recordings of spike trains are discretized in time and thus belong to finite dimensional spaces, due to the record resolution, the dimension of these spaces is so huge (from ten thousand up to one million) that it is neither realistic nor reasonable to model them by finite dimensional variables, and to apply usual independence tests. Several methods, such as the classical Unitary Events method (see \cite{GrunB} and the references therein), consist in binning the spike trains at first in order to deal with vectorial data with reduced dimension. However, it has been shown that these dimension reduction methods involve an information loss of more than 60\% in some cases, making this kind of preprocessing quite proscribed despite its simplicity of use. It is therefore more realistic and reasonable to model recordings of spike trains by  finite point processes, and to use independence tests specifically dedicated to such point processes. Asymptotic tests of independence between point processes have already been introduced in \cite{MTGAUE}, but in the particular case of homogeneous Poisson processes. Such a parametric framework is necessarily restrictive and even possibly inappropriate here, as the very existence of any precise underlying distribution for the point processes modelling spike train data is subject to broad debate  (see \cite{Pouzat09,RRGT}). We thus focus on non-parametric tests of independence for point processes. 
In this spirit, particular bootstrap methods under the name of {\it trial-shuffling} have been proposed in \cite{Pipa2003, Pipaet2003} for binned data with relatively small dimension, without proper mathematical justification. Besides the loss of information the binning data pre-processing involves, it appears that the test statistics chosen in these papers do not lead to  tests of asymptotic prescribed size as shown by our simulation study.

We here propose to construct new non-parametric tests of independence between two point processes, from the observation of $n$ independent copies of these point processes, with as few assumptions as possible on their underlying distributions. Our test statistics are based on $U$-statistics (see \cite[Chapter 5]{Serfling} for a key reference on $U$-statistics). The corresponding critical values are obtained from bootstrap or permutation approaches. It has been acknowledged that when both bootstrap and permutation approaches are available, permutation should be preferred, since the corresponding tests are exactly of the desired level \cite[p. 218]{EfronTib}. Nevertheless, we keep investigating them  together, as bootstrap methods -~through {\it trial-shuffling}~- are the usual references in neuroscience. Moreover, for specific $U$-statistics, the corresponding tests share the same properties: both are proved to be asymptotically of the prescribed size and consistent against any reasonable alternative, despite the fact that different tools are used to obtain these results. Indeed, the distance between the bootstrapped distribution and the initial distribution under independence is here directly studied for the bootstrap approach, unlike the permutation approach. Finally both procedures have good performance in practice when the sample size is moderate to small, as is often the case in neuroscience due to biological or economical reasons.

As $U$-statistics are usual tools for non-parametric statistical inference, many works deal with the application of bootstrap or permutation % type methods 
to $U$-statistics. From the original work of Arvesen \cite{Arvesen} about the Jackknife of $U$-statistics, to the recent one of Leucht and Neumann \cite{LeuchtNeumann09}, several papers \cite{Bickel-Freedman,Bretagnolle,agbst,Dehling94} have been devoted to the general problem of bootstrapping a $U$-statistic.  The use of bootstrap or permutation of $U$-statistics is specially considered in testing problems %(see \cite{HorvathHuskova,ChungRomano} for instance)
\cite{HorvathHuskova,ChungRomano}, in particular in dependence detection problems with the Kolmogorov-Smirnov type tests cited above %(see \cite{Romano89, vdvwellner96}).
\cite{Romano89, vdvwellner96}.

But  all those works exclusively focus on $U$-statistics of i.i.d. real valued random variables or vectors. Up to our knowledge, there is no previous work on the bootstrap or permutation of general $U$-statistics for i.i.d. pairs of point processes, as considered in the present paper.
 The main difficulty thus lies in the nature of the mathematical objects we handle here, that is point processes and their associated point measures which are random measures. The proofs of our results, although inspired by Romano's \cite{RomanoTR, Romano89} work and Hoeffding's \cite{Hoeffding52} precursor results on the permutation, are therefore more technical and complex on many aspects detailed in the sequel. In addition, we aim at obtaining the asymptotic distribution of the bootstrapped or permuted test statistics under independence, but also under dependence (see Theorem \ref{Convergence} and Theorem \ref{consistencyperm}). As concerns the permutation approach, such a result is, as far as we know, new even for more classical settings than point processes. It thus partially solves a problem stated as open question in \cite{vdvwellner96}.

\smallskip

This paper is organized as follows. 

We first present in Section \ref{ft} the testing problem, and introduce the main notations. Starting from existing works in neuroscience, we introduce our test statistics, based on general kernel-based $U$-statistics.

Section \ref{bootsec} is devoted to our bootstrap approach. Are given new general results about the consistency of the bootstrap for the considered $U$-statistics, expressed in terms of Wasserstein's metric as in \cite{Bickel-Freedman}. The convergence is studied under independence as well as under dependence. The corresponding bootstrap independence tests are therefore shown to be asymptotically of the desired size, and consistent against any reasonable alternative. 
 The impact of using Monte Carlo methods to approximate the bootstrap quantiles is also investigated in this section.

Section \ref{permsec} is devoted to the permutation approach which leads, by nature, to non-parametric independence tests exactly of the desired level, and this, even when a Monte Carlo method is used to approximate the permutation quantiles. Are then given new general results about the consistency of the permutation approach when the kernel of the $U$-statistic has a specific form. These results are still expressed in terms of Wasserstein's metric. As a consequence the corresponding permutation independence tests are proved to satisfy the same asymptotic properties as the bootstrap ones under the null hypothesis as well as under the same alternatives. 

As a comparison of the performance of our tests with existing ones in neuroscience, especially when the sample sizes are moderate or even small, a simulation study  is presented in Section \ref{simusec}.

A conclusion is given  in the last section. 

Finally notice that all proofs and some additional technical results can be found in a supplementary material.

% ---------------------------------------------------------------------------------------------------- %
\section{From neuroscience interpretations to general test statistics}
\label{ft}
% ---------------------------------------------------------------------------------------------------- %
 
% ------------------------------------------------------------ %
\subsection{The testing problem} 
% ------------------------------------------------------------ %

Throughout this article we consider finite point processes defined on a probability space $(\Omega,\mc{A},\PP)$ and observed on $[0,1]$, i.e. random point processes on $[0,1]$ whose total number of points is almost surely finite (see \cite{DV} for instance). Typically, in a neuroscience framework, such finite point processes may represent spike trains recorded on a given finite interval of time, and rescaled so that their values may be assumed to belong to $[0,1]$. 
The set $\calX$ of all their possible values consists of the countable subsets of $[0,1]$. It is  equipped with a metric $d_\calX$ that we introduce in \eqref{defmetric}. This metric, issued from the Skorohod topology, makes $\calX$ separable and allows to define accordingly borelian sets on $\calX$ and by extension on $\calX^2$ through the product metric.

The point measure $dN_x$ associated with an element $x$ of $\calX$ is defined for all measurable real-valued function $f$ by $\int_{[0,1]} f(u) dN_x(u)=\sum_{u\in x} f(u).$
In particular, the total number of points of $x$, denoted by $\# x$, is equal to $\int_{[0,1]} dN_x(u)$. Moreover, for a finite  point process $X$ defined on $(\Omega,\mc{A},\PP)$ and observed on $[0,1]$, 
$\int f(u) dN_X(u)$ becomes a real random variable, defined on the same probability space $(\Omega,\mc{A},\PP)$.

A pair $X=(X^1,X^2)$ of finite point processes defined on $(\Omega,\mc{A},\PP)$ and observed on $[0,1]$, has joint distribution $P$, with marginals $P^1$ and $P^2$ if $P(\mc{B})=\PP(X\in \mc{B})$, $P^1(\mc{B}^1)=\PP(X^1\in \mc{B}^1)$, and $P^2(\mc{B}^2)=\PP(X^2\in \mc{B}^2)$, for every borelian set $\mc{B}$ of $\calX^2$, and all borelian sets $\mc{B}^1$, $\mc{B}^2$ of $\calX$.

Given the observation of an i.i.d. sample  $\X_n=(X_1,\ldots,X_n)$ from the same distribution $P$ as $X$, with $X_i=(X^1_i,X^2_i)$ for every $i=1\ldots n$, we aim at testing $(H_0)$ $X^1$ and $X^2$ are independent against $(H_1)$ $X^1$ and $X^2$ are not independent,
which can also be written as
$$(H_0)\ P=P^1\otimes P^2 \quad \textrm{against}\quad (H_1)\ P\neq P^1\otimes P^2.$$

% ------------------------------------------------------------ %
\subsection{\label{coincsec} Independence test based on coincidences in neuroscience}
% ------------------------------------------------------------ %

In the neuroscience issue which initially motivated this work, the i.i.d. sample $\X_n=(X_1,\ldots,X_n)$ models pairs of rescaled spike trains issued from two distinct and simultaneously recorded neurons during $n$ trials. Those data are usually recorded on living animals that are repeatedly subject to the same stimulus  or that are repeatedly executing the same task. Because there are periods of rest between the records, it is commonly admitted that the $n$ trials are i.i.d. and that the considered i.i.d. sample model is actually realistic. Then, the main dependence feature that needs to be detected between both neurons corresponds to synchronization in time, referred to as coincidences~\cite{GrunB}. 
More precisely, neuroscientists expect to detect if such coincidences occur significantly, that is more than what may be due to chance. They speak in this case of a detected synchrony. 

In  \cite{MTGAUE}, the notion of coincidence count between two point processes $X^1$ and $X^2$ with delay $\delta$ ($\delta>0$) is defined by
\begin{equation}
\label{defnbc}
\Nbc_\delta(X^1,X^2)\!=\!\int_{[0,1]^2}  \!\!\!\!\! \1{|u-v|\leq \delta}dN_{X^1}(u) dN_{X^2}(v)=\!\!\!\!\!\!\!\!\sum_{u\in X^1,v\in X^2} \!\!\!\!\!\!\!\!\1{|u-v|\leq \delta}.
\end{equation}

Notice that other coincidence count functions have been used in the neuroscience literature such as the binned coincidence count function (i.e. based on binned data) introduced in \cite{Grunt} or its  shifted version \cite{GDGRA} (see also \cite{MTGAUE} for explicit formulae).  A further example of possible function used to detect dependence in neuroscience \cite{LaureCh} is of the form
\begin{equation}
\label{defw}
\varphi^{w}(X^1,X^2)=\int_{[0,1]^2} w(u,v) dN_{X^1}(u) dN_{X^2}(v).
\end{equation}
Under the assumption that both $X^1$ and $X^2$ are homogeneous Poisson processes, the independence test of \cite{MTGAUE} rejects $(H_0)$ when a test statistic based on $\sum_{i=1}^{n}\! \Nbc_\delta\!\!\left(X^1_i,X^2_i\right)$ is larger than a given critical value.
This critical value is deduced from the asymptotic Gaussian distribution of the test statistic under $(H_0)$. 
The test is proved to be asymptotically of the desired size, but only under the homogeneous Poisson processes assumption. However, it is now well-known that this assumption, as well as many other model assumptions, fails to be satisfied in practice for spike trains \cite{Pouzat09,RRGT}.

\subsection{General $U$-statistics as independence test statistics}

In the parametric homogeneous Poisson framework of \cite{MTGAUE}, the expectation of $\Nbc_\delta\left(X^1_i,X^2_i\right)$ has a simple expression as a function of $\delta$ and the intensities $\lambda_1$ and $\lambda_2$ of $X^1$ and $X^2$. Since $\lambda_1$ and $\lambda_2$ can be easily estimated, an estimator of this expectation can thus be obtained using the plug-in principle, and subtracted from $\Nbc_\delta\left(X^1_i,X^2_i\right)$ to lead to a test statistic with a centered asymptotic distribution under $(H_0)$.

In the present non-parametric framework where we want to make as few assumptions as possible on the point processes $X^1$ and $X^2$, such a centering plug-in tool is not available. We propose to use instead a self-centering trick, which amounts, combined with a rescaling step, to considering the statistic
\begin{equation}\label{statistic}
\frac{1}{n(n-1)}\sum_{i\neq i'\in\{1,\ldots,n\}} \left(\Nbc_\delta\left(X^1_i,X^2_i\right)-\Nbc_\delta\left(X^1_i,X^2_{i'}\right)\right).
\end{equation}
It is clear  that the function $\Nbc_\delta$ used in \cite{MTGAUE} suits the dependence feature the neuroscientists expect to detect in a spike train analysis. However, it is not necessarily the best choice for other kinds of dependence features to be detected in a general point processes analysis. Note furthermore that the statistic \eqref{statistic} can be written as a $U$-statistic of the i.i.d. sample $\X_n=(X_1,\ldots,X_n)$ with a symmetric kernel, as defined by Hoeffding \cite{Hoeffding48bis}.

Let us therefore consider the general independence test statistics which are $U$-statistics of the form
\begin{equation}\label{defstat}
U_{n,h}(\X_n)=\frac{1}{n(n-1)}\sum_{i\neq i'\in\{1,\ldots,n\}} h\left(X_i,X_{i'}\right),
 \end{equation}
where  $h: (\calX^2)^2 \to \R$ is a symmetric kernel such that:

\smallskip

$\pa{\mc{A}_{Cent}}\quad \textrm{\begin{tabular}{|l}  For all $n\geq 2$, $U_{n,h}(\X_n)$ is zero mean under $(H_0)$,\\
i.e. for $X_1$ and $X_2$, i.i.d. with distribution $P^1\otimes P^2$ on $\calX^2$,\\ $\esp{h\left(X_1,X_2\right)}=0$.
\end{tabular}}$

\smallskip

In the sequel, we call \emph{Coincidence case} the case where $h=h_{\Nbc_\delta}$, with
\begin{multline}\label{hphicoinc}
h_{\Nbc_\delta}(x,y)= \frac{1}{2}\ \big( {\Nbc_\delta}\left(x^1,x^2\right)+{\Nbc_\delta}\left(y^1,y^2\right)\\
-{\Nbc_\delta}\left(x^1,y^2\right)-{\Nbc_\delta}\left(y^1,x^2\right)\big),
\end{multline}
so that $U_{n,h_{\Nbc_\delta}}(\X_n)$ is equal to the statistic  \eqref{statistic}.

\smallskip

A more general choice, which of course includes the above \emph{Coincidence case}, is obtained by replacing $\Nbc_\delta$ by any generic integrable function $\varphi$. This is the \emph{Linear case}.
For any integrable function $\varphi$, the kernel $h$ is then taken equal to $h_\varphi$, with
\begin{equation}
\label{hphi}
h_\varphi(x,y)=\frac{1}{2}\left(\varphi\left(x^1,x^2\right)+\varphi\left(y^1,y^2\right)-\varphi\left(x^1,y^2\right)-\varphi\left(y^1,x^2\right)\right).
\end{equation}
This example is of utmost importance in the present work since it provides a first proved case of consistency for the permutation approach under the null hypothesis as well as under the alternative (see Theorem \ref{consistencyperm}). In this case, note that  $\pa{\mc{A}_{Cent}}$ is straightforwardly satisfied, i.e. $U_{n,h_\varphi}(\X_n)$ is zero mean under $(H_0)$. Note furthermore that $U_{n,h_\varphi}(\X_n)$ is an unbiased estimator of $$\int\int \varphi\left(x^1,x^2\right)\left(dP(x^1,x^2)-dP^1(x^1)dP^2(x^2)\right),$$ without any assumption on the underlying point processes. This is therefore a reasonable independence test statistic.
If $X^1$ and $X^2$ were finite dimensional variables with continuous distributions w.r.t. the Lebesgue measure, this test statistic would be closely related to generalized Kolmogorov-Smirnov tests of independence. For instance, the test statistics of Blum, Kiefer, and Rosenblatt \cite{BKR61}, Romano \cite{Romano89}, Van der Vaart and Wellner in \cite{vdvwellner96}  are equivalent  to
$\sqrt{n}\sup_{v^1\in \mc{V}^1,v^2\in \mc{V}^2}\left|U_{n,h_{\varphi_{(v^1,v^2)}}}(\X_n)\right|, $
where, respectively: 
\begin{itemize}
\item $\mc{V}^1=\mc{V}^2=\R$ and $\varphi_{(v^1,v^2)}(x^1,x^2)=\1{]-\infty,v^1]}(x^1)\1{]-\infty,v^2]}(x^2)$, 
\item $\mc{V}^1$ and $\mc{V}^2$ are countable V.-C. classes of subsets of $\R^d$,\\ and $\varphi_{(v^1,v^2)}(x^1,x^2)=\1{v^1}(x^1)\1{v^2}(x^2)$, 
\item $\mc{V}^1$ and $\mc{V}^2$ are well-chosen classes of real-valued functions,\\ and $\varphi_{(v^1,v^2)}(x^1,x^2)=v^1(x^1)v^2(x^2)$.
\end{itemize}
Note also the work of \cite{mein} based on integrals instead of supremum of similar quantities with $\varphi_{(v^1,v^2)}(x^1,x^2)= e^{iv^1 x^1}e^{iv^2x^2}$.
Thus, up to our knowledge, the existing  test statistics are  based on functions $\varphi$ of product type. However, as seen in Section \ref{coincsec}, when dealing with  point processes, natural functions $\varphi$, as for instance $\varphi_\delta^{coinc}$, are not of this type.

\subsection{Non-degeneracy of the $U$-statistics under $(H_0)$\label{nondegsec}}
Following the works of Romano \cite{Romano89} or Van der Vaart and Wellner \cite{vdvwellner96}, the tests we propose here are based on bootstrap and permutation approaches for the above general $U$-statistics. Most of the assumptions on $h$ depend on the chosen method (permutation or bootstrap) and are postponed to the corresponding section. However another assumption is common, besides $\pa{\mc{A}_{Cent}}$:

\smallskip

\noindent$\pa{\mc{A}_{non-deg}}\ \textrm{\begin{tabular}{|l}  For all $n\geq 2$, $U_{n,h}(\X_n)$ is non-degenerate under $(H_0)$,\\  i.e. for all $X_1$ and $X_2$, i.i.d. with distribution $P^1\otimes P^2$ on $\calX^2$,\\
$\var{\esp{h(X_1,X_2)|X_1}}\not =0$. 
\end{tabular}}$

\smallskip

This assumption is needed in all results with weak convergence to a Gaussian limit, as  its variance has to be strictly positive (see e.g. Proposition~\ref{TCLUstat} or Theorem~\ref{consistencyperm}). Since under $(H_0)$, $U_{n,h}(\X_n)$ is assumed to have zero mean, it is degenerate under $(H_0)$ if and only if for $X$ with distribution $P^1\otimes P^2$ and for $P^1\otimes P^2$-almost every $x$ in $\calX^2$, $\esp{h(x,X)}=0.$

In the \emph{Linear case}, this condition implies a very particular link between $\varphi$ and the distribution of the bivariate point process $X$, which is unknown. The following result gives some basic condition to fulfill $\pa{\mc{A}_{non-deg}}$ when $\varphi$ is the coincidence count function.
\begin{prop}\label{Prop:non_deg_coinc}
If the empty set is charged by the marginals, i.e. if $P^1(\{\emptyset\})>0$ and $P^2(\{\emptyset\})>0$ and if 
$\Nbc_\delta(X_1,X_2)$ (see \eqref{defnbc}) is not almost surely null under $(H_0)$, then  when $h$ is given by \eqref{hphicoinc}, $\pa{\mc{A}_{non-deg}}$  is satisfied.
\end{prop}
The proof can be found in the supplementary material together with a more informal discussion on the \emph{Linear case} with $\varphi=\varphi^w$ as given by  \eqref{defw}.

With respect to neuronal data, assuming that the processes may be empty is an obvious assumption  as there often exist trials (usually short) where, just by chance, no spikes have been detected. 
Moreover, practitioners usually choose $\delta$ large enough such that coincidences are observed in practice and therefore $\Nbc_\delta(X_1,X_2)$ is not almost surely null. 
Hence in practice, the non-degeneracy assumption is always satisfied in the \emph{Coincidence case}.

\medskip 

Throughout this article, $\pa{X_i}_{i}$ denotes a sequence of i.i.d. pairs of point processes, with $X_i=(X^1_i,X^2_i)$ of distribution $P$ on $\calX^2$, whose marginals  are $P^1$ and $P^2$ on $\calX$. For $n\geq 2$, let $\X_n=(X_1,\ldots,X_n)$ and $U_{n,h}(\X_n)$ as in \eqref{defstat}, with a fixed measurable symmetric kernel $h$ satisfying $\pa{\mc{A}_{Cent}}$. To shorten mathematical expression,  $U_n(\X_n)$ refers from now on to $U_{n,h}(\X_n)$.

% ---------------------------------------------------------------------------------------------------- %
\section{Bootstrap tests of independence}
\label{bootsec}
% ---------------------------------------------------------------------------------------------------- %

Since the distribution of the test statistic  $U_n(\X_n)$ is not free from the unknown underlying marginal distributions $P^1$ and $P^2$ under the null hypothesis $(H_0)$, we turn to a  classical bootstrap approach, which aims at mimicking it, for large, but also moderate or small sample sizes. 

To describe this bootstrap approach, and to properly state our results, 
we give below additional notations, and 
discuss the main assumptions. 

% ------------------------------------------------------------ %
\subsection{Additional notations: bootstrap and convergence formalism}
% ------------------------------------------------------------ %

For $j$ in $\{1,2\}$, let $P_n^j$ be the empirical marginal distribution defined by \
\begin{equation}\label{pnj} 
P_n^j=\frac{1}{n}\sum_{i=1}^n \delta_{X_i^j}.
 \end{equation}
A bootstrap sample from $\X_n$ is denoted by $\X_n^*=\left(X_{n,1}^*,\ldots,X_{n,n}^*\right)$, with $X_{n,i}^*=(X_{n,i}^{*1},X_{n,i}^{*2}),$ and is defined as an $n$ i.i.d. sample from the distribution $P_n^1\otimes P_n^2$. Then, the bootstrap distribution of interest is 
the conditional distribution of  $\sqrt{n}U_n(\X_n^*)$ given $\X_n$ to be compared 
with the initial distribution of $\sqrt{n}U_n(\X_n)\!$ under $(H_0)$. To state our convergence results as concisely as possible, we use the following classical formalism.

\begin{itemize}
\item For any functional $Z\!:\! (\calX^2)^n\! \to\! \R$, $\calL\left(Z,Q\right)$ denotes the distribution of $Z(\mathbb{Y}_n)$, where $\mathbb{Y}_n$ is an i.i.d.  sample from the distribution $Q$ on $\calX^2$. In particular, the distribution of $\sqrt{n}U_n(\X_n)$ under $(H_0)$ is denoted by $\calL\!\left(\!\sqrt{n}U_n,P^1\!\otimes\! P^2\right)$.
\item If the distribution $Q=Q(W)$ depends on a random variable $W$, $\calL\left(\left.Z,Q\right|W\right)$ is the conditional distribution of $Z(\mathbb{Y}_n)$, $\mathbb{Y}_n$ being an i.i.d. sample from the distribution $Q=Q(W)$, given $W$.

In particular, the conditional distribution of  $\sqrt{n}U_n(\X_n^*)$ given $\X_n$ is denoted by $\calL\left(\left.\sqrt{n}U_n,P^1_n\otimes P^2_n\right|\X_n\right)$.
\item " $Q$-a.s. in $(X_i)_i$"  at the end of a statement means that the statement only depends on the sequence $(X_i)_{i}$, where the $X_i$'s are i.i.d with distribution $Q$, and that 
there exists an event $\mathcal{C}$ only depending on $(X_i)_i$ such that $\proba{\mathcal{C}}=1$, on which the statement is true. Here $Q$ is usually equal to $P$.
\item "$Q_n \cvf{n\to+\infty} Q$" means that the sequence of distributions $(Q_n)_n$ converges towards $Q$ in the weak sense, that is for any real valued, continuous and bounded function $g$, $\int g(z) dQ_n(z) \to_{n\to+\infty} \int g(z) dQ(z).$
\item As usual, $\Estar{\cdot}$ stands for the conditional expectation given $\X_n$.
\end{itemize}

One of the aims of this work is to prove that the conditional distribution $\calL\left(\left.\sqrt{n}U_n,P^1_n\otimes P^2_n\right|\X_n\right)$ is asymptotically close to $\calL\left(\sqrt{n}U_n,P^1\otimes P^2\right)$. Following the historical paper by Bickel et Freedman \cite{Bickel-Freedman}, the closeness between these two distributions, which are both distributions on $\R$, is here measured via  the $\L^2$-Wasserstein's metric (also called Mallows' metric): % defined by
\begin{equation}\label{Wasserstein}
d_2^2(Q,Q')= \inf{} \big\{\esp{(Z-Z')^2},\ (Z,Z') \textrm{ with marginals }Q\textrm{ and }Q'\big\},
\end{equation}
for all the  distributions $Q$, $Q'$ with finite second order moments. 
Recall that convergence w.r.t. $d_2$  is equivalent to both weak convergence and convergence of second order moments.

% ------------------------------------------------------------ %
\subsection{Main assumptions}\label{ma}
% ------------------------------------------------------------ %

The random variables we deal with are not real-valued variables but point processes, so the assumptions needed in our results may be difficult to interpret in this setting. We therefore devote this whole section to their description and discussion.

% ------------------------------ %
In addition to Assumption $\pa{\mc{A}_{Cent}}$,  we need its following empirical version:

\smallskip

$\pa{\mc{A}_{Cent}^*}\quad \textrm{\begin{tabular}{|l} For $x_1=(x_1^1,x_1^2),\dots, x_n=(x_n^1,x_n^2)$ in $\calX^2$,\\
$\sum_{i_1,i_2,i'_1,i_2'=1}^n h\pa{\pa{x_{i_1}^1,x_{i_2}^2},\pa{x_{i'_1}^1,x_{i'_2}^2}}=0.$\end{tabular}}$

\smallskip

\noindent Notice that this assumption, as well as \pa{\mc{A}_{Cent}}, is fulfilled in the \emph{Linear case} where $h$ is of the form $h_\varphi$ given by \eqref{hphi}, but $\pa{\mc{A}_{Cent}^*}$ does not imply that $h$ is of the form $h_\varphi$ (see the supplementary material for a counterexample).

% ------------------------------ %
\paragraph{Moment assumptions} Due to the $\L^2$-Wasserstein's metric used here to study the consistency of the bootstrap approach, moment assumptions are required. In particular, the variance of $U_n(\X_n)$ should exist, i.e.

\smallskip

$\pa{\mc{A}_{Mmt}}\quad  \textrm{\begin{tabular}{|l}  For $X_1$ and $X_2$, i.i.d. with distribution $P$ on $\calX^2$,\\
 $\esp{h^2\pa{X_1,X_2}}<+\infty,$\end{tabular}}$

\smallskip

\noindent and more generally we  need:

\smallskip

$\pa{\mc{A}_{Mmt}^{*}}\quad\textrm{\begin{tabular}{|l} 
For $X_1,X_2,X_3,X_4$ i.i.d. with distribution $P$ on $\calX^2$,\\
and for $i_1,i_2,i_1',i_2'$ in $\ac{1,2,3,4}$, \\
$\esp{h^2\pa{\big(X_{i_1}^1,X_{i_2}^2\big),\big(X_{i_1'}^1,X_{i_2'}^2\big)}}<+\infty.$
\end{tabular}}$

\smallskip

\noindent Notice that when $\pa{\mc{A}_{Mmt}^{*}}$ is satisfied, this implies that
\begin{itemize}
\item $\pa{\mc{A}_{Mmt}}$ is satisfied (taking $i_1=i_2$, $i_1'=i_2'$, and $i_1'\neq i_1$),
\item for  $X\sim P$, $\esp{h^2\pa{X,X}}<+\infty$ (taking  $i_1=i_2=i_1'=i_2'$),
\item for $X_1$, $X_2$ i.i.d with distribution $P^1\otimes P^2$, $\esp{h^2\pa{X_1,X_2}} < +\infty$ (taking $i_1,i_2,i_1',i_2'$ all different).
\end{itemize}
A sufficient condition for $\pa{\mc{A}_{Mmt}^{*}}$  and  $\pa{\mc{A}_{Mmt}}$  to be satisfied is that there exist positive constants $\alpha_1$, $\alpha_2,$ $C$ such that for every $x=(x^1,x^2),y=(y^1,y^2)$ in $\calX^2$,
$ |h(x,y)|\leq C\left((\# x^1)^{\alpha_1} + (\#y^1)^{\alpha_1}\right)\left((\# x^2)^{\alpha_2} + (\#y^2)^{\alpha_2}\right),$
with\\  $\ds{E}[ (\#X^1)^{4\alpha_1}]<+\infty$ and $\ds{E}[(\#X^2)^{4\alpha_2}]<+\infty$.

In the \emph{Linear case} where $h$ is of the form $h_\varphi$ given by \eqref{hphi}, a possible sufficient condition is 
that there exist some positive constants $\alpha_1$, $\alpha_2$, and $C$ such that for every $x^1,x^2$ in $\calX$, 
$|\varphi(x^1,x^2)|\leq C (\# x^1)^{\alpha_1}(\# x^2)^{\alpha_2},$
with $\ds{E}[ (\#X^1)^{4\alpha_1}]<+\infty$ and $\ds{E}[(\#X^2)^{4\alpha_2}]<+\infty$. 
In particular, in the \emph{Coincidence case}, the coincidence count function $\Nbc_\delta$ satisfies: for every $x^1,x^2$ in $ \calX$,
$|\Nbc_\delta(x^1,x^2)|\leq (\# x^1)(\# x^2).$ So, $\pa{\mc{A}_{Mmt}^{*}}$ and  $\pa{\mc{A}_{Mmt}}$ are satisfied as soon as $\ds{E}[ (\#X^1)^{4}]<+\infty$ and $\ds{E}[(\#X^2)^{4}]<+\infty$.

Such moment bounds for the total number of points of the processes are in fact satisfied by many kinds of point processes: discretized point processes at resolution $0<r<1$ (see \cite{MTGAUE} for a definition), which have at most $1/r$ points, Poisson processes, whose total number of points obeys a Poisson distribution having exponential moments of any order, and point processes with bounded conditional intensities, which can be constructed by thinning homogeneous Poisson processes (see \cite{ogatathin}).
Similar moment bounds can also be obtained (see \cite{HRBR}) 
 for linear stationary Hawkes processes with positive interaction functions that are classical models in spike train analysis (see e.g. \cite{PSCR,MTGAUE}). 
 This finally may be extended to point processes whose conditional intensities are upper bounded by intensities of linear stationary Hawkes processes with positive interaction functions, by thinning arguments. This includes more general Hawkes processes (see \cite{BM}) 
 and in particular Hawkes processes used to model inhibition in spike train analysis (see \cite{HRBR, MTGAUE, ieee, RRGT}).

% ------------------------------ %
\paragraph{Continuity of the kernel} 
The set $\calX$ can be embedded in the space $\calD$ of c\`adl\`ag functions on $[0,1]$ through the identification $$N:x\in\calX\mapsto \pa{N_x:t\mapsto \int_0^1\1{u\leq t} dN_x(u)}\in\calD.$$ 
Notice that the quantity $N_x$ is actually the counting process associated with $x$ (see \cite{Bremaud} for instance): at time $t$, $N_x(t)$ is the number of points of $x$ less than $t$. 
Now consider the uniform Skorohod topology on $\calD$ (see \cite{Billingsley2}), associated with the metric $d_\calD$ defined by
$$d_\calD(f,g)=\inf \ac{\varepsilon>0\ ;\ \exists\lambda\in\Lambda, \left\{\begin{array}{l}\sup_{t\in[0,1]} |\lambda(t)-t|\leq \varepsilon,\\ \sup_{t\in [0,1]}|f(\lambda(t))-g(t)|\leq \varepsilon \end{array}\right.},$$
where $\Lambda$ is the set of strictly increasing, continuous mappings of $[0,1]$ onto itself. 
Notice that here, $\lambda$ represents a uniformly small deformation of the time scale. 
Thanks to the identification $N$ above, $\calX$ can then be endowed with the topology induced by  $d_{\calX}$ defined on $\calX$ by
\begin{equation}\label{defmetric}
d_{\calX}(x,x')=d_{\calD}(N(x),N(x'))\quad \textrm{for every $x$, $x'$ in $\calX$}.\end{equation}
As an illustration, if $x$ and $x'$ are in $\calX$, for $\varepsilon$ in $(0,1)$, $d_{\calX}(x,x')\leq \varepsilon$ implies that $x$ and $x'$ have the same cardinality, and for $k$ in $\ac{1,\dots,\#x}$, the $k$\up{th} point of $x$ is at distance less than $\varepsilon$ from the $k$\up{th} point of $x'$.
Since $(\calD,d_\calD)$ is a separable metric space, so are  $(\calX,d_{\calX})$, $\left(\calX^2,d_{\calX^2}\right),$ where $d_{\calX^2}$ is the  product metric defined from $d_{\calX}$ (see \cite[p 32]{Dudley}),  and $\left(\calX^2\times\calX^2,d\right),$ where $d$, the product metric defined from $d_{\calX^2}$,  is given by
\begin{equation}\label{distance}
d\big((x,y),(x',y')\big)=\sup \ac{\sup_{j=1,2}\ac{d_{\calX}(x^j,x^{'j})},\sup_{j=1,2}\ac{d_{\calX}(y^j,y^{'j}) }},
\end{equation}
 for every $x=(x^1,x^2)$, $y=(y^1,y^2)$, $x'=(x^{'1},x^{'2})$, $y'= (y^{'1},y^{'2})$ in $\calX^2$.\\
The kernel $h$ chosen to define the $U$-statistic $U_n(\X_n)$ in \eqref{defstat}  should satisfy:

\smallskip

$\pa{\mc{A}_{Cont}}\quad \textrm{\begin{tabular}{|l}
There exists a subset $\mathcal{C}$ of $\calX^2\times \calX^2$, such that \\
 (i) $h$ is continuous on $\mathcal{C}$ for the topology
induced by $d$,  \\
(ii) $(P^1\otimes P^2)^{\otimes 2}(\mathcal{C})=1$.
\end{tabular}}$

\smallskip

Here are some examples in the \emph{Linear case} for which $\pa{\mc{A}_{Cont}}$ holds. 

\begin{prop}\label{hinHw}
Let $w : [0,1]^2 \to \R$ be a continuous integrable function. Then the kernel $h_{\varphi^w}$ defined on $\calX^2\times \calX^2$ by \eqref{defw} and  \eqref{hphi} is continuous w.r.t. the topology induced by $d$, defined by \eqref{distance}.
\end{prop}
\noindent The above result  does not apply to $h_{\Nbc_\delta}$ but the following one  holds. 
\begin{prop}\label{hinH}
The coincidence count kernel $h_{\Nbc_\delta}$ defined on $\calX^2\times \calX^2$ by \eqref{defnbc} and  \eqref{hphi} is continuous w.r.t. the topology induced by $d$,  on 
\begin{multline}
\label{1indep2} \mc{C}_\delta=\big\{\pa{(x^1,x^2),(y^1,y^2)}\in \calX^2\times \calX^2\ ;\\ 
\pa{\ac{x^1}\cup\ac{y^1}} \cap \pa{\ac{x^2 \pm \delta}\cup\ac{y^2 \pm\delta}} =\emptyset\big\}.
\end{multline}
\end{prop}

As suggested in \cite{MTGAUE}, when dealing with discretized point processes at resolution $r$, the right choice for $\delta$ is $kr+r/2$ for an integer $k$, so $(P^1\otimes P^2)^{\otimes 2}(\mathcal{C}_\delta)=1$, and $h_{\Nbc_\delta}$ satisfies $\pa{\mc{A}_{Cont}}$. Furthermore, when dealing with independent point processes with conditional intensities,  those processes may be constructed by thinning  two independent Poisson processes $X$ and $X'$. Hence, in this case, the probability $(P^1\otimes P^2)^{\otimes 2}$ of $\mc{C}_\delta$ in \eqref{1indep2} is larger than $\proba{X\cap(X'\pm\delta)=\emptyset}, $ whose value is 1. So when dealing with point processes with conditional intensities,   $h_{\Nbc_\delta}$ also satisfies $\pa{\mc{A}_{Cont}}$.

% ------------------------------------------------------------ %
\subsection{Consistency of the bootstrap approach}
% ------------------------------------------------------------ %

The validity of the bootstrap approach for our independence tests is due to the following consistency result.

\begin{thm}\label{Convergence}
For every $n\geq 2$, let $P_n^j$ for $j=1,2$ be the empirical marginal distributions defined by \eqref{pnj}. Then, under $\pa{\mc{A}_{Cent}}$, $\pa{\mc{A}_{Cent}^*}$,  $\pa{\mc{A}_{Mmt}^{*}}$ and $\pa{\mc{A}_{Cont}}$,   
$$d_2\left(\calL\left( \left. \sqrt{n} U_{n}, P_n^1\otimes P_n^2 \right|\X_n\right) \!, \calL\left( \sqrt{n}U_{n},P^1\!\otimes \! P^2\right)  \right) \!\!\cv{n\to+\infty}\! \!0,\ \textrm{$P$-a.s. in $(X_i)_i$}.$$\end{thm}

The proof follows similar arguments to the ones of  \cite{Bickel-Freedman} for the bootstrap of the mean, or to \cite{Dehling94} and \cite{LeuchtNeumann09} for the bootstrap of $U$-statistics. The main novel point here consists in using the identification \eqref{distance} and the properties of  the separable Skorohod metric space $(\calD,d_\calD)$, where weak convergence of sample probability distributions is available (see \cite{Varadarajan58}). This theorem derives in fact from the following two propositions which may be useful in various frameworks. The first one states a non-asymptotic result, while the second one gives rather natural results of convergence.

\begin{prop}\label{inegalite_Wasserstein}
Under $\pa{\mc{A}_{Cent}}$, $(\mc{A}_{Cent}^*)$,  $\pa{\mc{A}^*_{Mmt}}$, with the notation of Theorem \ref{Convergence}, there exists an absolute constant $C>0$ such that for $n\geq2$,
\begin{multline*}
d_2^2\left(\calL\left( \left. \sqrt{n} U_{n}, P_n^1\otimes P_n^2 \right|\X_n\right) , 
\calL\left( \sqrt{n}U_{n},P^1\otimes P^2\right)  \right)\\
\leq C \inf \Big\{\Estar{\pa{h\pa{ Y_{n,a}^*, Y_{n,b}^*}-h\pa{Y_a, Y_b}}^2},
Y_{n,a}^*\sim P_n^1\otimes P_n^2,\  Y_a \sim P^1\otimes P^2,\\ \textrm{and } (Y_{n,b}^*,Y_b) \textrm{ is an independent copy of }(Y_{n,a}^*,Y_a)\Big\}.
\end{multline*}
\end{prop}

\emph{Comment.} In the above proposition, the infimum is taken over all the possible distributions of $(Y_{n,a}^*,Y_a)$ having the correct marginals, $(Y_{n,b}^*,Y_b)$ being just an independent copy of $(Y_{n,a}^*,Y_a)$. In particular, $Y_{n,a}^*$ is not necessarily independent of $Y_a$.

 \begin{prop}
\label{LGNhcarre}
If $\esp{|h(X_1,X_2)|}<+\infty$, then
\begin{equation}\label{cvh}
U_{n}(\X_n)\! \!\! \cv{n\to+\infty} \! \! \!\esp{h(X_1,X_2)}\!\!=\!\!\!\int\!\! h(x,x')dP(x) dP(x'),\textrm{$P$-a.s. in $(X_i)_i$.}
\end{equation}
Under $\pa{\mc{A}_{Mmt}^{*}}$, one  moreover obtains that $P$-a.s. in $(X_i)_i$,
$$\frac{1}{n^4} \sum_{i,j,k,l=1}^n h^2\pa{\pa{X_i^1,X_j^2},\pa{X_k^1,X_l^2}} \\ \cv{n\to+\infty} \esp{h^2\pa{\pa{X_1^1,X_2^2},\pa{X_3^1,X_4^2}}}.$$
\end{prop}

% ------------------------------------------------------------ %
\subsection{Convergence of cumulative distribution functions (c.d.f.) and quantiles}
% ------------------------------------------------------------ %

As usual, $\mathcal{N}(m,v)$ stands for the Gaussian distribution with mean $m$ and variance $v$, $\Phi_{m,v}$ for its c.d.f. and $\Phi_{m,v}^{-1}$ for its quantile function.
From the results of Rubin and Vitale \cite{RubinVitale80} generalizing Hoeffding's \cite{Hoeffding48bis}  decomposition of non-degenerate $U$-statistics to the case where the $X_i$'s are non necessarily real valued random vectors, a Central Limit Theorem for $U_n(\X_n)$ can be easily derived. It is expressed here using the $\L^2$-Wasserstein's metric, and is thus
slightly stronger than the one stated in Equation (1.1) of \cite{HuskovaJanssen}.
 
 \begin{prop}\label{TCLUstat}
Assume that $h$ satisfies $\pa{\mc{A}_{non-deg}}$, $\pa{\mc{A}_{Cent}}$, and $\pa{\mc{A}_{Mmt}}$.
Let $\sigma^2_{P^1 \otimes P^2}$ be defined by 
\begin{equation}\label{sigma} \sigma^2_{P^1 \otimes P^2}= 4\var{\esp{h\left(X_1,X_2\right)|X_1}},
\end{equation} 
when  $X_1$ and $X_2$ are $P^1\otimes P^2$-distributed. Then 
$$d_2\pa{\calL\left(\sqrt{n}U_{n},P^1\otimes P^2 \right), \mathcal{N}(0, \sigma^2_{P^1 \otimes P^2})} \cv{n\to+\infty} 0.$$
 
\end{prop}
\emph{Comments.} 

{\it (i)} Notice that $\pa{\mc{A}_{non-deg}}$ is equivalent to $\sigma^2_{P^1 \otimes P^2}>0$.  In the case where  $\pa{\mc{A}_{non-deg}}$ does not hold, i.e. if $\sigma^2_{P^1 \otimes P^2}=0$, the quantity $\sqrt{n}U_{n}(\X_n)$ tends in probability towards 0. In this case, Theorem \ref{Convergence} implies that the two distributions  $ \calL\left( \left. \sqrt{n} U_{n}, P_n^1\otimes P_n^2 \right|\X_n\right)$ and $\calL\left( \sqrt{n}U_{n},P^1\!\otimes \! P^2\right)$ are not only close, but that they 
are actually both tending to the Dirac mass in $0$. Indeed, degenerate U-statistics of order $2$ have a faster rate of convergence than  $\sqrt{n}$ (see \cite{ArconesGineLimit} for instance for explicit limit theorems). So in this degenerate case, one could not use $\sqrt{n}U_{n}(\X_n)$ as a test statistic anymore (without changing the normalization). But as mentioned above, $\pa{\mc{A}_{non-deg}}$ is usually satisfied in practice  (see Section \ref{nondegsec} for the {\it Coincidence case}). 

{\it (ii)} Let us introduce, as in \cite{HuskovaJanssen}, an estimator of $\sigma^2_{P^1 \otimes P^2}$, but which is here corrected to be unbiased under $(H_0)$, namely:
$$\hat{\sigma}^2=\frac{4}{n(n-1)(n-2)} \sum_{i,j,k\in \{1,\ldots,n\}, \#\{i,j,k\}=3} h(X_i,X_j)h(X_i,X_k),$$
and the statistic:
\begin{equation}
\label{stat_S_n}
S_{n}=\sqrt{n}U_{n}(\X_n)/\hat\sigma.
\end{equation}
From Proposition \ref{TCLUstat} combined with Slutsky's lemma and the law of large numbers for $U$-statistics of order $3$, one easily derives  that under $(H_0)$, 
$S_n$ converges in distribution to $\mc{N}(0,1)$. 
This leads to a rather simple but asymptotically satisfactory test: the test which rejects $(H_0)$ when  $|S_n|\geq \Phi_{0,1}^{-1}(1-\alpha/2)$  is indeed asymptotically of size $\alpha$. It is  also consistent against any reasonable alternative $P$, satisfying $\pa{\mc{A}_{Mmt}}$ and  such that $\esp{h(X,X')}\neq 0$, for $X$, $X'$ i.i.d. with distribution $P$. Such a purely asymptotic test may of course suffer from a lack of power when the sample size $n$ is small or even moderate, which is typically the case for the application in neuroscience described in Section \ref{ft} for biological reasons (from few tens up to few hundreds at best). 
Though the bootstrap approach is mainly justified by asymptotic arguments, the simulation study  presented in Section \ref{simu} shows its efficiency in a non-asymptotic context, compared to this simpler test.

\smallskip

As Proposition \ref{TCLUstat} implies that the limit distribution of $\sqrt{n}U_{n}(\X_n)$ has a continuous c.d.f., the convergence of the conditional c.d.f. or quantiles of the considered bootstrap distributions holds. Note that these conditional bootstrap distributions are discrete, so the corresponding quantile functions are to be understood as the generalized inverses of the cumulative distribution functions.

\begin{coro}\label{coroBoot}
For $n\geq 2$, with the notation of Theorem \ref{Convergence}, let $\X_n^*$ be a bootstrap sample,  i.e. an i.i.d $n$-sample from the distribution $P^1_n\otimes P^2_n$. 
Let $\X_n^{\independent}$ be another i.i.d. $n$-sample from the distribution $P^1\otimes P^2$ on $\calX^2$. Under $\pa{\mc{A}_{non-deg}}$ and the assumptions of Theorem \ref{Convergence},  $$\sup_{z\in\R}\left| \PP\!\pa{ \left.\!\sqrt{n} U_{n}\left(\X_n^*\right)\leq z\right|\X_n}\!-\!\PP\!\pa{\!\sqrt{n}U_{n}(\X_n^{\independent})\leq z}\right|
\cv{n\to+\infty} 0,\ P\textrm{-a.s. in }(X_i)_i.$$
If moreover, for $\eta$ in $(0,1)$, $q_{\eta,n}^*(\X_n)$ denotes the conditional $\eta$-quantile of $\sqrt{n}U_{n}(\X_n^*)$ given $\X_n$ and  $q_{\eta,n}^{\independent}$ denotes the $\eta$-quantile of $\sqrt{n}U_{n}(\X_n^{\independent})$,  
\begin{equation}
\label{cvquantBoot}
\ |q_{\eta,n}^*(\X_n)-q_{\eta,n}^{\independent}|\cv{n\to+\infty} 0, \textrm{ $P$-a.s. in $(X_i)_i$.}
\end{equation}
\end{coro}

% ------------------------------------------------------------ %
\subsection{\label{asymp}Asymptotic properties of the bootstrap tests}
% ------------------------------------------------------------ %

We are interested in the asymptotic behavior of sequences of tests all based on test statistics of the form $\sqrt{n}U_{n}(\X_n)$. The bootstrap approach, whose consistency is studied above, allows to define bootstrap-based critical values for these tests. Note that the permutation approach studied in Section \ref{permsec} is based on the same test statistics, but with critical values obtained by permutation. Hence we introduce here a condensed and common formalism for the upper-, lower- and two-tailed tests considered in this work, taking into account that the only change in our two considered approaches concerns the critical values. This will help to state our results in the shortest manner.

Let $\alpha$ be fixed in $(0,1)$, and $q$ be a sequence of upper and lower critical values:
$$q=\left(q^+_{\alpha,n}(\X_n), q^-_{\alpha,n}(\X_n)\right)_{n\geq 2}.$$
From this sequence $q$, let us now define the family $\Gamma(q)$ of three sequences of tests $\Delta^{+}=(\Delta_{\alpha,n}^{+})_{n\geq 2}$, $\Delta^{-}=(\Delta_{\alpha,n}^{-})_{n\geq 2}$, and $\Delta^{+/-}=(\Delta_{\alpha,n}^{+/-})_{n\geq 2}$, where
 \begin{equation}  \label{deftest}
\left\{
\begin{array}{l l l l}
\Delta_{\alpha,n}^{+}(\X_n)&=&\1{\sqrt{n}U_{n}(\X_n)>q^+_{\alpha,n}(\X_n)} &\mbox{(upper-tailed test)},\\
\Delta_{\alpha,n}^{-}(\X_n)&=&\1{\sqrt{n}U_ {n}(\X_n) < q^-_{\alpha,n}(\X_n)} &\mbox{(lower-tailed test)},\\
\Delta_{\alpha,n}^{+/-}(\X_n)&=&\max\pa{\Delta_{\alpha/2,n}^{+}(\X_n),\Delta_{\alpha/2,n}^{-}(\X_n)}&\mbox{(two-tailed test)},\end{array}\right.
\end{equation}
the last test being implicitly defined by the corresponding choices in $\alpha/2$.

Of course, $q$, $\Gamma(q)$, as well as $\Delta^+$, $\Delta^-$ and $\Delta^{+/-}$, depend on the choice of $\alpha$, but since $\alpha$ is fixed at the beginning, to keep the notation as simple as possible, this dependence is, like the one in $h$, omitted in the notations. 

Depending on the choice of $q$, the classical asymptotic properties that can be expected to be satisfied by $\Gamma(q)$ are $\pa{\mc{P}_{size}}$ and $\pa{\mc{P}_{consist.}}$ defined by:
\begin{align*}
\pa{\mc{P}_{size}} &\ \textrm{\begin{tabular}{|l}
Each $\Delta=\pa{\Delta_{\alpha,n}}_{n\geq 2}$ in $\Gamma(q)$ is asymptotically of size $\alpha$, \\
 i.e.  $\proba{\Delta_{\alpha,n}(\X_n)=1}\to_{n\to+\infty}\! \alpha$  if $P=P^1\otimes P^2$;
\end{tabular}}\\
\\
\pa{\mc{P}_{consist.}}&\ \textrm{\begin{tabular}{|l}
Each $\Delta=\pa{\Delta_{\alpha,n}}_{n\geq 2}$ in $\Gamma(q)$ is consistent,\\
 i.e. $\proba{\Delta_{\alpha,n}(\X_n)=1}\to_{n\to+\infty} 1$, 
 for  every $P$ such that \\
 ~~\bl $\int h(x,x') dP(x) dP(x') >0$ if $\Delta= \Delta^{+}$,\\
 ~~\bl $\int h(x,x') dP(x) dP(x') <0$ if $\Delta= \Delta^{-}$,\\
 ~~\bl $\int h(x,x') dP(x) dP(x') \neq 0$ if $\Delta= \Delta^{+/-}$.\\
\end{tabular}}
\end{align*}

Following Corollary~\ref{coroBoot}, our bootstrap tests family is defined from \eqref{deftest} by $\Gamma(q^*)$, with
\begin{equation}\label{qetoile}
q^*=\pa{q_{1-\alpha,n}^*(\X_n), q_{\alpha,n}^*(\X_n)}_{n\geq 2}.
\end{equation}

%Note that $q^*$ is random, depending on $\X_n$, and that it may be exactly computed by considering the $n^{2n}$ possible bootstrap samples. The algorithmic complexity of such an exact computation is usually so large that a Monte Carlo  of the bootstrap quantiles, based on resampling from the original  $\X_n$, is preferred  in practice. This Monte Carlo step is considered in Section \ref{MCB}.

\begin{thm}\label{thniveauconsistance}
Let $\Gamma(q^*)$ be the family of  tests defined by \eqref{deftest} and \eqref{qetoile}.  
 If  $\pa{\mc{A}_{non-deg}}$, $\pa{\mc{A}_{Cent}}$, $(\mc{A}_{Cent}^*)$,  $\pa{\mc{A}_{Mmt}^{*}}$ and $\pa{\mc{A}_{Cont}}$ hold,
 then $\Gamma(q^*)$ satisfies both $\pa{\mc{P}_{size}}$ and $\pa{\mc{P}_{consist.}}$.
\end{thm}
\emph{Comments.} In the \emph{Linear case} where $h$ is equal to $h_\varphi$ defined by \eqref{hphi}, \\
$\int h(x,x') dP(x) dP({x'}) = \int \varphi(x^1,x^2)\left [dP(x^1,x^2)-dP^1(x^1)dP^2(x^2)\right].$
This means that under the assumptions of Theorem  \ref{thniveauconsistance}, the two-tailed test of $\Gamma(q^*)$ is consistent against any alternative such that $\int \varphi(x^1,x^2) dP(x^1,x^2)$ differs from what is expected under $(H_0)$, i.e. $\int \varphi(x^1,x^2) dP^1(x^1)dP^2(x^2)$.

$(i)$ In particular, in the \emph{Coincidence case} where $h$ is equal to $h_{\Nbc_\delta}$ defined by \eqref{hphicoinc},  the assumptions of  Theorem  \ref{thniveauconsistance} are fulfilled for instance if $X^1$ and $X^2$ are discretized at resolution $r$, with $\delta=kr+r/2$ for some integer $k$, or if $X^1$ and $X^2$ have bounded conditional intensities, with $\delta$ large enough so that $\Nbc_\delta(X^1,X^2)$ is not a.s. null. Theorem \ref{thniveauconsistance} means in such cases that the corresponding two-tailed test is asymptotically of power $1$, for any  alternative $P$ such that $\int \1{|v-u|\leq \delta}\esp{dN_{X^1}(u)dN_{X^2}(v)}\neq \int \1{|v-u|\leq \delta} \esp{dN_{X^1}(u)}\esp{dN_{X^2}(v)}].$
Note that no $\delta$ ensuring this condition can be found if heuristically, the repartition of the delays $|v-u|$ between points of $X^1$ and $X^2$ is the same under $(H_0)$ and under $(H_1)$. For neuroscientists, it means that the cross-correlogram (histogram of the delays, classically represented as a first description of the data) does not show different behaviors in the dependent and independent cases. This would only occur if the dependence could not be measured in terms of delay between points. 

$(ii)$  Furthermore, when $\varphi$ is equal to $\varphi^w$ defined by \eqref{defw} with a continuous integrable function $w$ (see Proposition \ref{hinHw}), Theorem \ref{thniveauconsistance} means that the corresponding two-tailed test is consistent against any alternative such that 
$\beta_{w}=\int w(u,v) \left(\esp{dN_{X^1}(u)dN_{X^2}(v)}-\esp{dN_{X^1}(u)}\esp{dN_{X^2}(v)}\right)\neq 0.$
For the function $w$ chosen in \cite{LaureCh} and under specific Poisson assumptions, $\beta_{w}$  is linked to a coefficient in the Haar basis of the so-called interaction function, which measures the dependence between both processes $X^1$ and $X^2$. Working non asymptotically, one of the main result of \cite{LaureCh} states, after reformulation in the present setting, that if $\beta_{w}$ is larger than an explicit lower bound,  then the second kind error rate of the upper-tailed test is less than a prescribed $\beta$ in $(0,1)$. Theorem  \ref{thniveauconsistance} thus generalizes the result of \cite{LaureCh} to a set-up with much less reductive assumptions on the underlying stochastic models, but in an asymptotic way.

\medskip

Whereas the above family of bootstrap tests $\Gamma(q^*)$ involves an exact computation of the conditional quantiles $q_{\eta,n}^*\pa{\X_n}$, in practice, these quantiles are approximated by a Monte Carlo method. More precisely, let $(B_n)_{n\geq 2}$ be a sequence of possible numbers of Monte Carlo iterations, such that $B_n\to_{n\to+\infty} +\infty$. For $n\geq 1$, let $\pa{\X_n^{*1},\ldots,\X_n^{*B_n}}$ be $B_n$ independent bootstrap samples from $\X_n$. Set $\pa{U^{*1},\ldots,U^{*B_n}}=\pa{U_{n}\pa{\X_n^{*1}},\ldots,U_{n}\pa{\X_n^{*B_n}}}$, and introduce its corresponding order statistic $\pa{U^{*(1)},\ldots,U^{*(B_n)}}$.  The considered family of Monte Carlo bootstrap tests  is then defined from \eqref{deftest} by $\Gamma(q^*_{MC})$, with
\begin{equation}\label{qetoileB}
q^*_{MC}=\left(\sqrt{n}U^{*(\lceil(1-\alpha)B_n\rceil)},\sqrt{n} U^{*(\lfloor\alpha B_n\rfloor +1)}\right)_{n\geq 2}.
\end{equation}

\begin{prop}\label{bootMonteCarlo}
Let $\Gamma(q^*_{MC})$ be the family of Monte Carlo boostrap tests defined by \eqref{deftest} and $q^*_{MC}$ in \eqref{qetoileB}. Under the same assumptions as in Theorem~\ref{thniveauconsistance}, then $\Gamma(q^*_{MC})$ also satisfies both $\pa{\mc{P}_{size}}$ and $\pa{\mc{P}_{consist.}}$.
\end{prop}

% ---------------------------------------------------------------------------------------------------- %
\section{Permutation tests of independence}
\label{permsec}
% ---------------------------------------------------------------------------------------------------- %
%All along this section,  $(X_i)_i$ still denotes a sequence of i.i.d. pairs of point processes, with $X_i=(X^1_i,X^2_i)$ of distribution $P$ on $\calX^2$, whose marginals  are respectively $P^1$ and $P^2$ on $\calX$. For every $n\geq 2$, $\X_n=(X_1,\ldots,X_n)$. 
%\comMag{As specified above, $U_n(\X_n)$ still stands for $U_n(\X_n)$, for a fixed symmetric kernel $h$.}

%\comMag{J'ai repoussé l'hypothèse que $h$ est de la forme $h_\varphi$ car pour les résultats non-asymptotiques, on n'en a pas besoin !}

% ------------------------------------------------------------ %
\subsection{The permutation approach and its known non-asymptotic properties}

Consider a random permutation  $\rperm$, uniformly distributed on the set $\Sn{n}$ of permutations of $\{1,\ldots,n\}$, and independent of $\X_n$. Then a permuted sample from $\X_n$ is defined by $\X_n^{\rperm}=\pa{X_1^{\rperm},\dots,X_n^{\rperm}}$  with $X_i^{\rperm}=\pa{X_i^1,X_{{\rperm}(i)}^2}.$
In the same formalism as for the bootstrap approach, for 
 $n\geq 2$ and $\eta$ in $(0,1)$, let $q_{\eta,n}^\star\pa{\X_n}$ denote the $\eta$-quantile of $\loi{\sqrt{n}U_{n},P_n^\star\middle|\X_n}$,  where $P_n^\star$ stands for the conditional distribution of $\X_n^\rperm$ given $\X_n$.
The family of permutation tests is then defined  by $\Gamma\pa{q^\star}$ (see \eqref{deftest}), with
\begin{equation}\label{qstar}
q^\star=\left(q_{1-\alpha,n}^\star\pa{\X_n}, q_{\alpha,n}^\star\pa{\X_n}\right)_{n\geq 2}.
\end{equation}
As for the bootstrap approach, in practice, the sequence of quantiles $q^\star$ is approximated by a Monte Carlo method.
So, let $(B_n)_{n\geq 2}$ be a sequence of numbers of Monte Carlo iterations, such that $B_n\to_{n\to+\infty} +\infty$. For $n\geq 1$,  let $\pa{\Pi_n^1,\ldots,\Pi_n^{B_n}}$ be a sample of $B_n$ i.i.d. random permutations uniformly distributed on $\Sn{n}$. Set $\pa{U^{\star 1},\ldots,U^{\star B_n}}=\pa{U_{n}\pa{\X_n^{\Pi_n^1}},\ldots,U_{n}\pa{\X_n^{\Pi_n^{B_n}}}}$ and $U^{\star B_n+1}=U_{n}\pa{\X_n}$, the $U$-statistic computed on the original sample $\X_n$. The order statistic associated with  $\pa{U^{\star 1},\ldots,U^{\star B_n+1}}$ is denoted as usual by
$\pa{U^{\star(1)},\ldots,U^{\star(B_n+1)}}$. The considered family of Monte Carlo permutation tests  is then defined from \eqref{deftest} by $\Gamma\pa{q^\star_{MC}}$, with
\begin{equation}\label{qstarB}
q^\star_{MC}=\pa{\sqrt{n}U^{\star(\lceil (1-\alpha)(B_n+1)\rceil)},
\sqrt{n} U^{\star(\lfloor \alpha(B_n+1)\rfloor +1)}}_{n\geq 2}.
\end{equation}

The main advantage of the above families of permutation tests is that any test $\Delta_{\alpha,n}$ from either $\Gamma(q^\star)$ or $\Gamma(q^\star_{MC})$ is exactly of the desired level $\alpha$ i.e. 
\begin{equation}\label{nonasymp}
\mbox{if }P=P^1\otimes P^2, \quad \PP(\Delta_{\alpha,n}(\X_n)=1)\leq \alpha.
\end{equation}
Such non-asymptotic results for the permutation tests are well-known (see for instance  \cite[Lemma 1]{RomanoWolf} and \cite{PesarinSalmaso}). Though similar results are since recently available for bootstrap tests in other settings \cite{DurotRoz, Arlot, FLRB}, there is no known exact counterpart for the bootstrap in the present context.

% ------------------------------------------------------------ %
\subsection{Consistency of the permutation approach}
% ------------------------------------------------------------ %
In this section, we focus on the \emph{Linear case} where $h$ is of the form $h_\varphi$  for some integrable function $\varphi$, as defined in \eqref{hphi}. Indeed, it is the most general case for which we are able to prove a combinatorial Central Limit Theorem under any alternative as well as under the null hypothesis (Theorem \ref{consistencyperm}). Hence in this section, $U_n$ refers to $U_{n,h_\varphi}$. Notice that the centering assumption $(\mc{A}_{Cent})$ is then always satisfied by $U_n(\X_n)$. We here only need the following moment assumption:

\smallskip

\noindent$\pa{\mc{A}_{\varphi,Mmt}}\  \textrm{\begin{tabular}{| l}  For $(X^1,X^2)$ with distribution $P$ or $P^1\otimes P^2$ on $\calX^2$,\\ $\ds{E}\!\cro{\varphi^4\!\!\pa{X^1\!,X^2}\!}\!<\!\infty.$\end{tabular}}$

\smallskip

Though we have no exact counterpart of Theorem \ref{Convergence} for our 
permutation approach, 
the following result combined with  Proposition \ref{TCLUstat} gives a similar result.

\begin{thm}
\label{consistencyperm}
For all $n\geq 2$, let $P_n^\star$ be the conditional distribution of a permuted sample given $\X_n$. In the  \emph{Linear case} where the kernel $h$ is of the form \eqref{hphi} for an integrable function $\varphi$, under $\pa{\mc{A}_{non-deg}}$ and  $(\mc{A}_{\varphi,Mmt})$, with the notations of Section \ref{bootsec},
\begin{equation}
\label{CVstatpermH0}
d_2\pa{\loi{\sqrt{n}U_{n},P_n^\star\middle|\X_n},\mc{N}\pa{0,\sigma_{P^1\otimes P^2}^2}} \cvproba{n\to+\infty} 0,
\end{equation}
\vspace{-0.5cm}\\
where $\cvproba{}$ stands for the usual convergence in $\PP$-probability.
\end{thm}

\emph{Comments.} As pointed out above, unlike the bootstrap approach, the conditional permutation distribution of the test statistic is not here directly compared to the initial distribution of the test statistic under the null hypothesis. It is in fact compared to the Gaussian limit distribution of the test statistic under the null hypothesis, when the non-degeneracy assumption $\pa{\mc{A}_{non-deg}}$ holds. Moreover, the convergence occurs here in probability and not almost surely, but note that no continuity assumption for the kernel $h_\varphi$ is used anymore. The price to pay is that the moment assumption is stronger than the one used for the bootstrap. This assumption, due to our choice to use an existing Central Limit Theorem for martingale difference arrays in the proof, is probably merely technical and maybe dispensable. Indeed, the result of Theorem \ref{consistencyperm} is close to asymptotic results for permutation known as {\it combinatorial Central Limit Theorems} \cite{Hoeff51, PesarinSalmaso}, where this kind of higher moment assumption can be replaced by some Lindeberg conditions \cite{Hajek,Motoo,HoChen}. However, all these existing  results can only be applied directly in our case either when $(X_i)_i$ is deterministic or under the null hypothesis. Up to our knowledge, no combinatorial Central Limit Theorem has been proved for non deterministic and non exchangeable variables, like here under any alternative.

%\comRef{[Ref2, Min8]: je rajoute une explication}. 
%This stronger assumption is the consequence of the following technical fact: integer moments of sums of permuted variables are easily computable \comMel{je trouve que c'est bizarre comme explication, on ne voit pas trop ce que ça vient faire...}. Therefore to show for instance that the reminder terms in the variance (moment of order 2) vanish, one needs to go to moments of order $2\times 2=4$.
%\comMel{Autre formulation: based on Theorem \ref{diffmart}, we need to show the convergence in probability to a constant of a sum of conditional second order moments, and in particular, we prove that its variance vanishes, thus $2\times 2=4$ moments are needed.}
%This assumption is probably merely technical. \comRef{[Ref1 point1] ref qui manque du TCL combinatorial + [point4] dire ce qui est nouveau } Indeed this result is very close to asymptotic results for permutations known as {\it combinatorial Central Limit Theorem} \cite{Hoeff51, PesarinSalmaso} where this kind of higher moment assumption can be replaced by some Lindeberg conditions \cite{Motoo,HoChen}. However all these results can be applied directly in our case either when $(X_i)_i$ \comMel{J'ai changé les $(X_n)_n$ en $(X_i)_i$} is deterministic or under $(H_0)$. Up to our knowledge, no combinatorial Central Limit Theorem has been proved for non deterministic non exchangeable variables, which is what happens under $(H_1)$.

 The above result is thus one of the newest results presented here and its scope is well beyond the only generalization to the point processes setting. Indeed, because it holds not only under $(H_0)$ but also under $(H_1)$, it goes further than any existing one for independence test statistics such as the ones of Romano \cite{Romano89}. 
 The behavior under $(H_1)$ of the permuted test statistic of Van der Vaart and Wellner was also left as an open question  in \cite{vdvwellner96}. 
 
 The proof is presented in the supplementary material.

\smallskip

From Theorem \ref{consistencyperm}, we deduce the following corollary.

\begin{coro}\label{coroPerm}
Under the assumptions of Theorem~\ref{consistencyperm} and with the notations of Proposition \ref{TCLUstat}, for $\eta$ in $(0,1)$, 
$$q_{\eta,n}^\star\pa{\X_n} \cvproba{n\to+\infty} \Phi^{-1}_{0,\sigma^2_{P^1 \otimes P^2}}(\eta).$$
\end{coro}

% ------------------------------------------------------------ %
\subsection{Asymptotic properties of the permutation tests}
% ------------------------------------------------------------ %

As for the bootstrap tests, we obtain the following result.
%Like for the bootstrap approach, the critical values used in the tests defined by \eqref{deftestperm} are random, depending on the observed sample $\X_n$. Since $\rperm$ is  uniformly distributed on the set $\Sn{n}$ of permutations of $\{1,\ldots,n\}$ and independent of $\X_n$, given $\X_n$, the conditional distribution of $\sqrt{n}U_{n,h_\varphi}\!\pa{\X_n^\rperm}$ is discrete and takes the values $\ac{\sqrt{n}U_{n,h_\varphi}\!\pa{\X_n^\perm}}_{\perm\in\Sn{n}}$. Let $\sqrt{n}U_{n,h_\varphi}^{(1)}\!\pa{\X_n}\leq \ldots\leq \sqrt{n}U_{n,h_\varphi}^{(n!)}\!\pa{\X_n}$ be the ordered values.
%Then, $q^\star_{\varphi,\eta,n}\pa{\X_n}=\sqrt{n}U_{n,h_\varphi}^{(\lceil n!\eta\rceil)}\!\pa{\X_n}$, and it is therefore possible to compute it exactly. 

\begin{thm}
\label{NivConsistencyPerm}
Let $\Gamma(q^\star)$ and $\Gamma(q^\star_{MC})$ be the families of permutation and Monte Carlo permutation tests defined by \eqref{deftest} combined with \eqref{qstar} and \eqref{qstarB} respectively. In the \emph{Linear case}, if $\pa{\mc{A}_{non-deg}}$ and $(\mc{A}_{\varphi,Mmt})$ hold, then $\Gamma(q^\star)$  and  $\Gamma(q^\star_{MC})$ both satisfy $\pa{\mc{P}_{size}}$ and $\pa{\mc{P}_{consist.}}$.
\end{thm}

\section{\label{simu}Simulation study}
\label{simusec}
% ---------------------------------------------------------------------------------------------------- %

In this section, we study our testing procedures from a practical point of view, by giving estimations of the size and the power for various underlying distributions that are coherent with real neuronal data. This allows to verify the usability of these new methods in practice, and to compare them with existing classical methods. A real data sets study and a more operational and complete method for neuroscientists derived from the present ones is the subject of  \cite{nous-appli}. % an ongoing work.
 The programs have been optimized, parallelized in \texttt{C++} and interfaced with \texttt{R}. The code is available at  \texttt{\url{https://github.com/ybouret/neuro-stat}}.

\subsection{Presentation of the study}
All along the study, $h$ is taken equal to $h_{\varphi^{coinc}_\delta}$ (see \eqref{hphicoinc}), where $\varphi^{coinc}_\delta$ is defined in \eqref{defnbc} and $\alpha=0.05$. We only present the results for upper-tailed tests, but an analogous study has been performed for lower-tailed tests with similar results. Five different testing procedures are compared. 
\subsubsection{Testing procedures}
\begin{description}
\item[(\texttt{CLT})] Test based on the Central Limit Theorem for $U$-statistics (see Proposition \ref{TCLUstat}) which rejects $(H_0)$ when the test statistic $S_{n}$ in \eqref{stat_S_n}  is larger than the $(1-\alpha)$-quantile of the standard normal distribution.
\item[(\texttt{B})] Monte Carlo bootstrap upper-tailed test of $\Gamma(q^*_{MC})$ (\eqref{deftest} and \eqref{qetoileB}). 
\item[(\texttt{P})] Monte Carlo permutation upper-tailed test of $\Gamma(q^\star_{MC})$ (\eqref{deftest} and \eqref{qstarB}).  
\item[(\texttt{GA})] Upper-tailed tests introduced in \cite[Definition 3]{MTGAUE} under the notation $\Delta^+_{GAUE}(\alpha)$, based on a Gaussian approximation of the total number of coincidences. 
\item[(\texttt{TS})] Trial-shuffling test based on a Monte Carlo approximation of the $p$-value introduced in \cite[Equation (3)]{Pipaet2003}, but adapted to the present notion of coincidences. This test is the reference distribution-free method for neuroscientists. More precisely, let $C(\X_n) = \sum_{i=1}^n \varphi^{coinc}_\delta\pa{X_i^1,X_i^2}$
be the total number of coincidences of  $\X_n$.
The trial-shuffling method consists in uniformly drawing with replacement $n$ i.i.d. pairs of indices $\ac{\pa{i^*(k),j^*(k)}}_{1\leq k\leq n}$ in $\ac{(i,j), 1\!\leq \!i\neq \!j\leq\! n}$, and considering the associated ${TS}$-sample $\X_n^{TS}\!=\! \pa{\!\pa{\!X_{i^*(k)}^1,X_{j^*(k)}^2\!}\!}_{1\leq k\leq n}$. 
The Monte Carlo  $p$-value is defined by $\alpha^{TS}_B = \frac{1}{B}\sum_{b=1}^B \1{C\pa{\X_n^{TS,b}} \geq C\pa{\X_n}},$ where $\X_n^{TS,1},\ldots, \X_n^{TS,B}$ are $B$ independent ${TS}$-samples, and the test rejects $(H_0)$ if $\alpha^{TS}_B\leq \alpha$. 
This procedure is therefore close in spirit to our bootstrap procedure except that it is applied on a non-centered quantity under $(H_0)$, namely $C(\X_n)$.
\end{description}
The number $B$ of steps in the Monte Carlo methods is taken equal to $10 000$.

\subsubsection{Simulated data}\label{simdat}

Various types of point processes are simulated here to check the distribution-free character of our approaches and to investigate their limits. Of course, each of the considered point processes satisfies the moment assumptions on the number of points so that the theorems in this article can be applied. 
From now on and to be coherent with the neuroscience application which originally motivated this work, the point processes are simulated on $[0,0.1]$. Indeed the following experiments have been done to match neurophysiological parameters \cite{MTGAUE,Grunt} and the classical necessary window for detection is usually of duration $0.1$ seconds.

\paragraph{Estimation of the size} The three data sets simulated under $(H_0)$ consist of i.i.d. samples of pairs of independent point processes. For simplicity, both processes have the same distribution, though this is not required. \vspace{-0.1cm}

\begin{description}
\item[Exp. A] Homogeneous Poisson processes on $[0,0.1]$ with intensity $\lambda=60$. 

\item[Exp. B] Inhomogeneous Poisson processes with intensity $f_\lambda : t\in [0,0.1] \mapsto \lambda t $ and $\lambda = 60$. 

\item[Exp. C] Hawkes Processes as detailed in \cite{MTGAUE}
i.e. point processes with conditional intensity $\lambda(t) = \max\left(0,\mu - \int_0^t\nu\ \1{]0,r]}(t-s)\, dN_{X}(s)\right),$
for $t$ in $[0,0.1]$, with spontaneous intensity $\mu=60$, refractory period $r=0.001$, and $\nu>\mu$
such that for all point $T$ in $X$ and $t$ in $]T,T+r]$, $\lambda(t)=0$. 
%Intuitively, $\lambda(t)dt$ gives the probability that a new point appears in $]t, t+dt]$ given the past. In particular, t
This choice of $\nu$ prevents two points to occur at a distance less than the refractory period $r$ to  reflect typical neuronal behavior. This model is also sometimes called Poisson process with dead time.
\end{description}

\paragraph{Study of the power} The three data sets simulated under $(H_1)$ are such that the number of coincidences is larger than  expected under $(H_0)$. The models (injection or Hawkes) are classical in neuroscience and already used in \cite{MTGAUE,GrunB}. 
\begin{description}
\item[Exp. D]  Homogeneous injection model. $X^1 =X_{ind}^1 \cup X_{com}$ and $X^2=X_{ind}^2 \cup X_{com}$,  $X_{ind}^1$ and $X_{ind}^2$ being two independent homogeneous Poisson processes with intensity $\lambda_{ind}=54$, $X_{com}$ being a common homogeneous Poisson process with intensity $\lambda_{com}=6$, independent of $X_{ind}^1$ and $X_{ind}^2$. 
\item[Exp. E]  Inhomogeneous injection model. Similar to {\bf Exp. D}, $X_{ind}^1$ and $X_{ind}^2$ being two independent inhomogeneous Poisson processes with intensity $f_{\lambda_{ind}}$ (see {\bf Exp. B}), $\lambda_{ind}=54$, $X_{com}$ being a homogeneous Poisson process with intensity $\lambda_{com}=6$, independent of $X_{ind}^1$ and $X_{ind}^2$.
\item[Exp. F] Dependent bivariate Hawkes processes. 
The coordinates $X^1$ and $X^2$ of a same pair respectively have the conditional intensities:\\
$\lambda^1(t)\!=\! \max\Big\{0,\mu -\!\int_0^t\!\!\nu \1{]0,r]}(t-s)\, dN_{X^1}(s) + \!\int_0^t\!\!\eta\1{]0,u]}(t-s)\, dN_{X^2}(s) \Big\},$\\
$\lambda^2(t)\!=\! \max\Big\{0,\mu -\! \int_0^t\!\!\nu \1{]0,r]}(t-s)\, dN_{X^2}(s) +\! \int_0^t\!\!\eta\1{]0,u]}(t-s)\, dN_{X^1}(s)\Big\},$
with the spontaneous intensity $\mu= 54$, the interaction intensity $\eta~=~6$ in the period designated by $u=0.005$ and the refractory period designated by $r=0.001$ with $\nu \gg \mu + \eta u$ such that once again, $\lambda^j(t)$ is null on each $]T,T+r]$, for $T$ in $X^j$. We arbitrarily took $\nu=50(2\mu + \eta)$. 
\end{description}

\subsection{Results}

\begin{figure}[h!]
\begin{tabular}{c|c}
{\bf Exp. A} & {\bf Exp. D}\vspace{-8pt}\\
\includegraphics[angle=-90,scale=0.22]{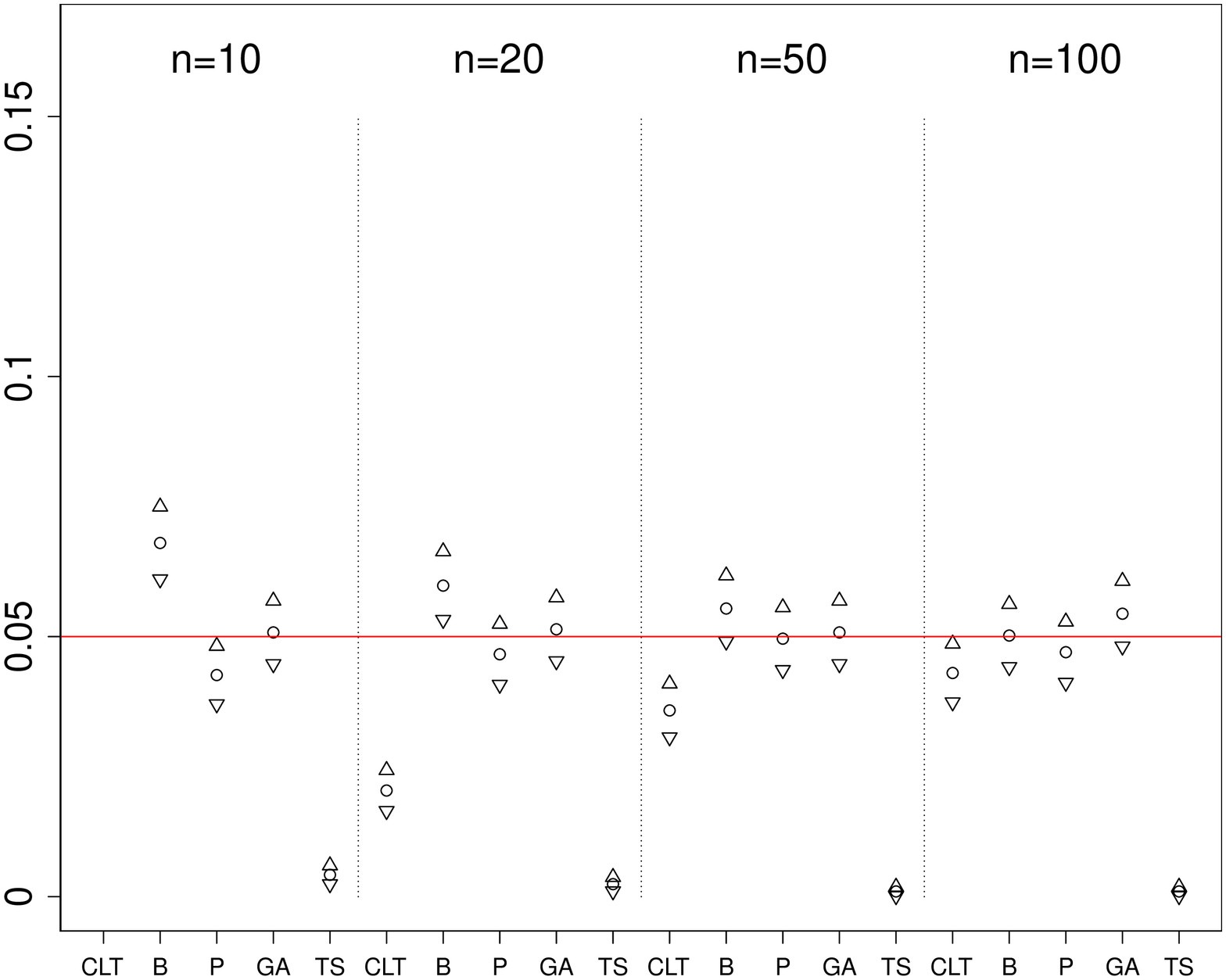}& \includegraphics[angle=-90,scale=0.22]{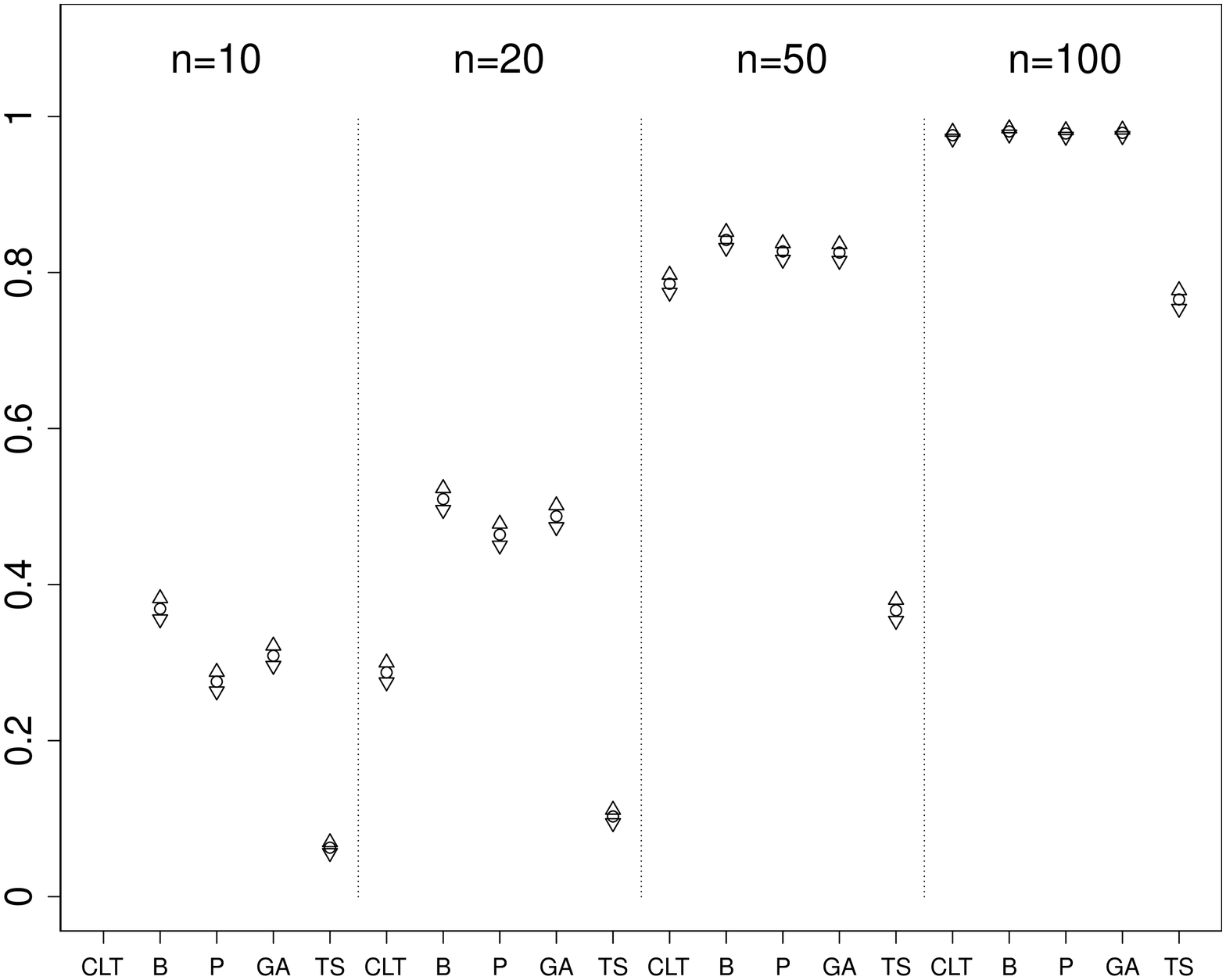}\vspace{5pt}\\
\hline \vspace{-10pt}\\
{\bf Exp. B} & {\bf Exp. E}\vspace{-8pt}\\ 
\includegraphics[angle=-90,scale=0.22]{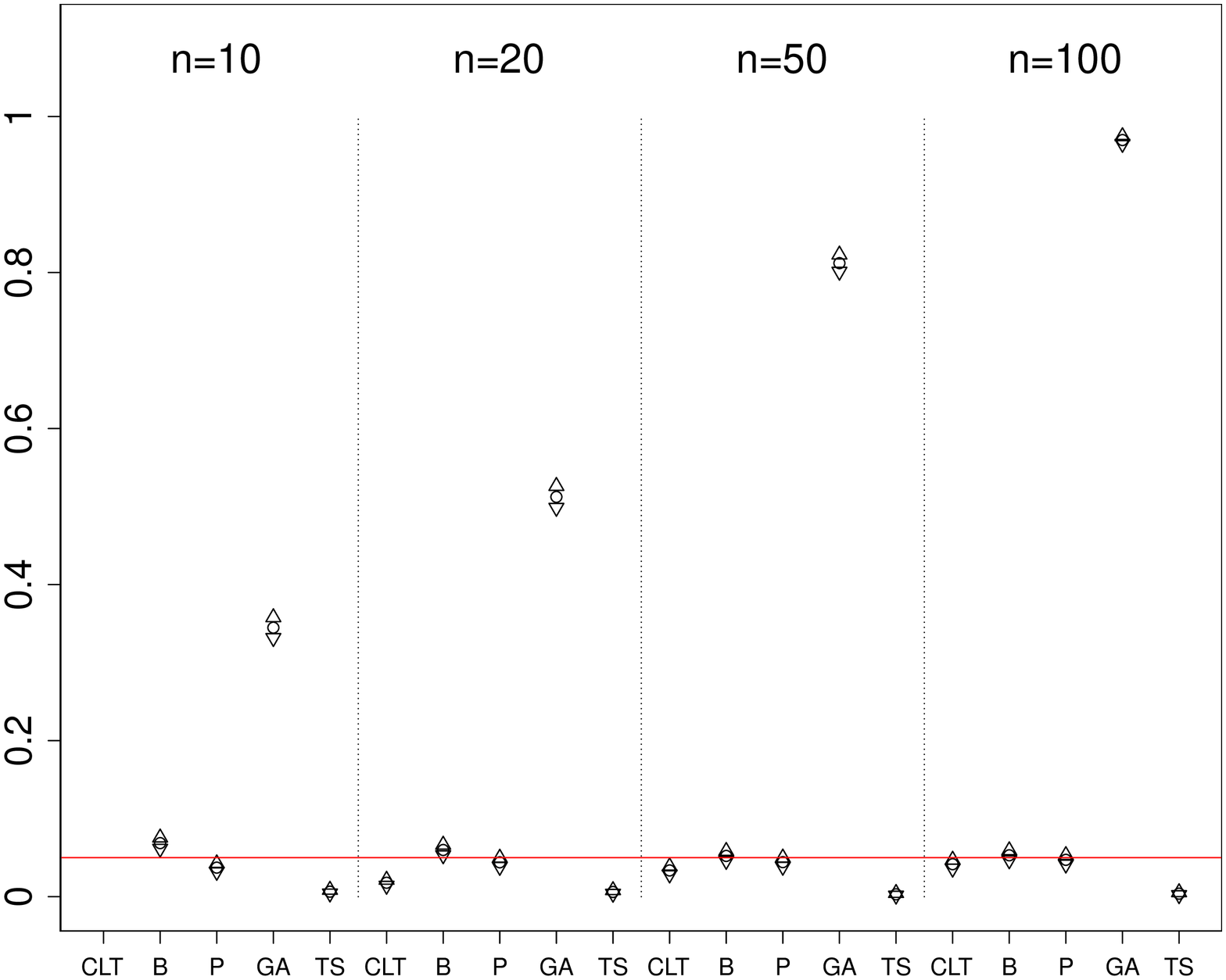}& \includegraphics[angle=-90,scale=0.22]{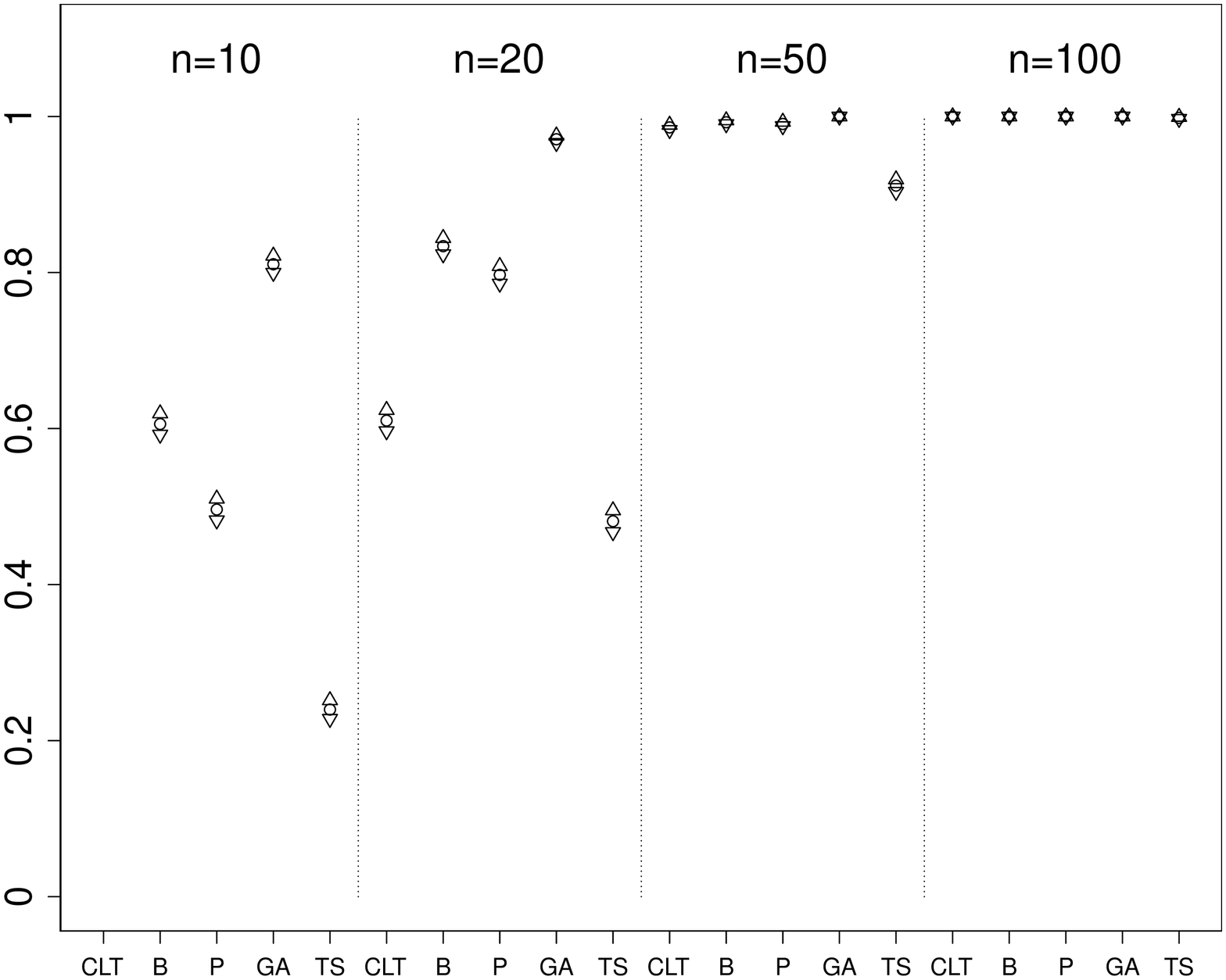}\vspace{5pt}\\
\hline \vspace{-10pt}\\
{\bf Exp. C} & {\bf Exp. F}\vspace{-8pt}\\
 \includegraphics[angle=-90,scale=0.22]{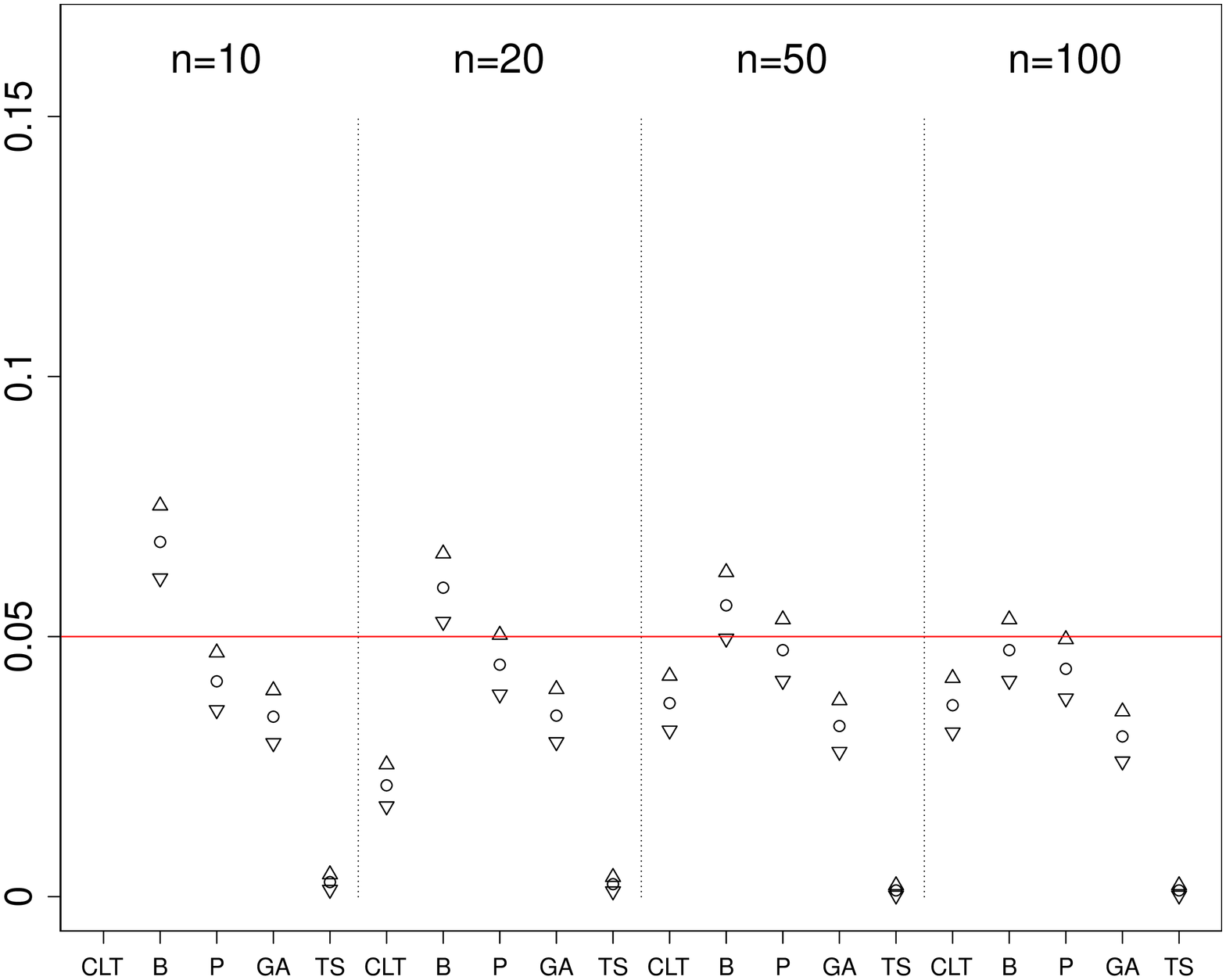}& \includegraphics[angle=-90,scale=0.22]{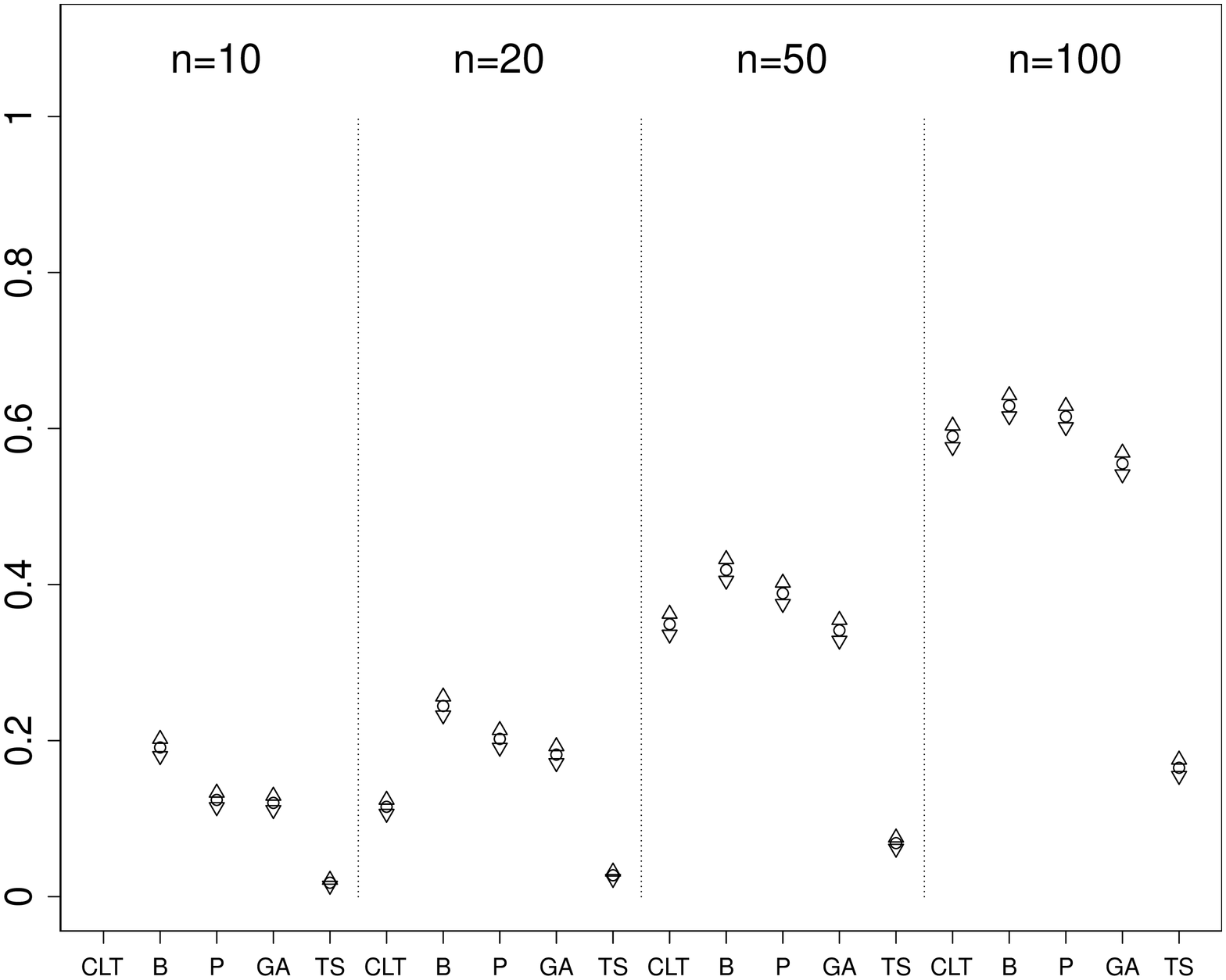}\vspace{5pt}\\
\end{tabular}
\caption{\label{nessai}Estimated sizes and powers for various numbers of trials $n$, all the tests being performed with a level $\alpha=0.05$.  The circles represent the percentage of rejection on 5000 simulations for each method, the triangles represent the corresponding endpoints of a 95\% confidence interval. The corresponding experiments are described in Section \ref{simdat}.}
\end{figure}

\paragraph{Varying number of trials $n$}
 In Figure \ref{nessai}, the delay is fixed at $\delta = 0.01$ and the number $n$ of trials varies in \ac{10, 20, 50, 100}. Note that when the number of trials is too small ($n=10$), the estimated variance in (\texttt{CLT}) is sometimes negative, therefore, the test cannot be implemented. 

The left hand side of  Figure \ref{nessai} corresponds to estimated sizes. On the one hand, one can see in the case of homogeneous Poisson processes ({\bf Exp. A}) and in the case of refractory Hawkes processes ({\bf Exp. C}) that the methods (\texttt{CLT}), (\texttt{B}), (\texttt{P}) and (\texttt{GA}) are quite equivalent, but the size (first kind error rate) seems less  controlled in the bootstrap approach (\texttt{B}) especially for small numbers of trials. Yet, one can see the convergence of the size of the bootstrap test towards $\alpha$ as the number of trials goes to infinity, which  illustrates Proposition \ref{bootMonteCarlo}.  Note that the (\texttt{CLT}) test  also has a well controlled size even if it cannot be used for very small $n$.
On the other hand, in the case of inhomogeneous Poisson processes ({\bf Exp. B}), one can see that the (\texttt{GA}) test has a huge size and  is thus inadequate here. Indeed it is based on the strong assumption that the data are homogeneous Poisson processes though they are in fact strongly nonstationary. The test tends thus to reject the independence null hypothesis even when the data are independent. 
Finally, in the three considered cases, the (\texttt{TS}) approach has a very small size, and is thus too conservative as one can see in the power study. The study of \cite{nous-appli} shows that this lack of performance is 
due to the fact that the (\texttt{TS}) approach is applied here on a non correctly centered quantity.

 The right hand side of Figure \ref{nessai}  corresponds to  estimated powers, which increase as $n$ grows. This is in line with the consistency of the tests. 
Now, as it could be expected when looking at its estimated sizes, for the (\texttt{TS}) approach, the estimated powers are distinctly lower than the ones for the other methods, which confirms its conservative behavior. 
The other approaches are more similar in {\bf Exp. D} or {\bf Exp.~F} though (\texttt{B}) clearly seems to outperform all  tests, but at the price of a less controlled size. 
Note that in the inhomogeneous case ({\bf Exp. E}), (\texttt{GA}) seems to have the best power, but this time, at the price of a totally uncontrolled size.

This part of the simulation study illustrates the convergences of the size and the power of the bootstrap and permutation tests introduced here. The permutation approach seems to actually guarantee the best control of the size as expected, as compared with the bootstrap approach. Nevertheless both approaches are quite effective for any considered kind of point processes  and any sample size, unlike the (\texttt{GA}) test which has very restrictive assumptions. The reference method (\texttt{TS}) for neuroscientists is clearly too conservative. Moreover, the (\texttt{CLT}) test seems to have also satisfying results, but with a slower convergence than the (\texttt{B}) and (\texttt{P}) tests. 
This seems to illustrate that the conditional bootstrap and permutation distributions give better approximations of the original one under independence than a simple Central Limit Theorem. This phenomenon is well-known as the second order accuracy of the bootstrap in more classical frameworks.

\begin{figure}[h!]
\begin{tabular}{c|c}
{\bf Exp. D} & {\bf Exp. E} \vspace{-8pt}\\
\includegraphics[angle=-90,scale=0.23]{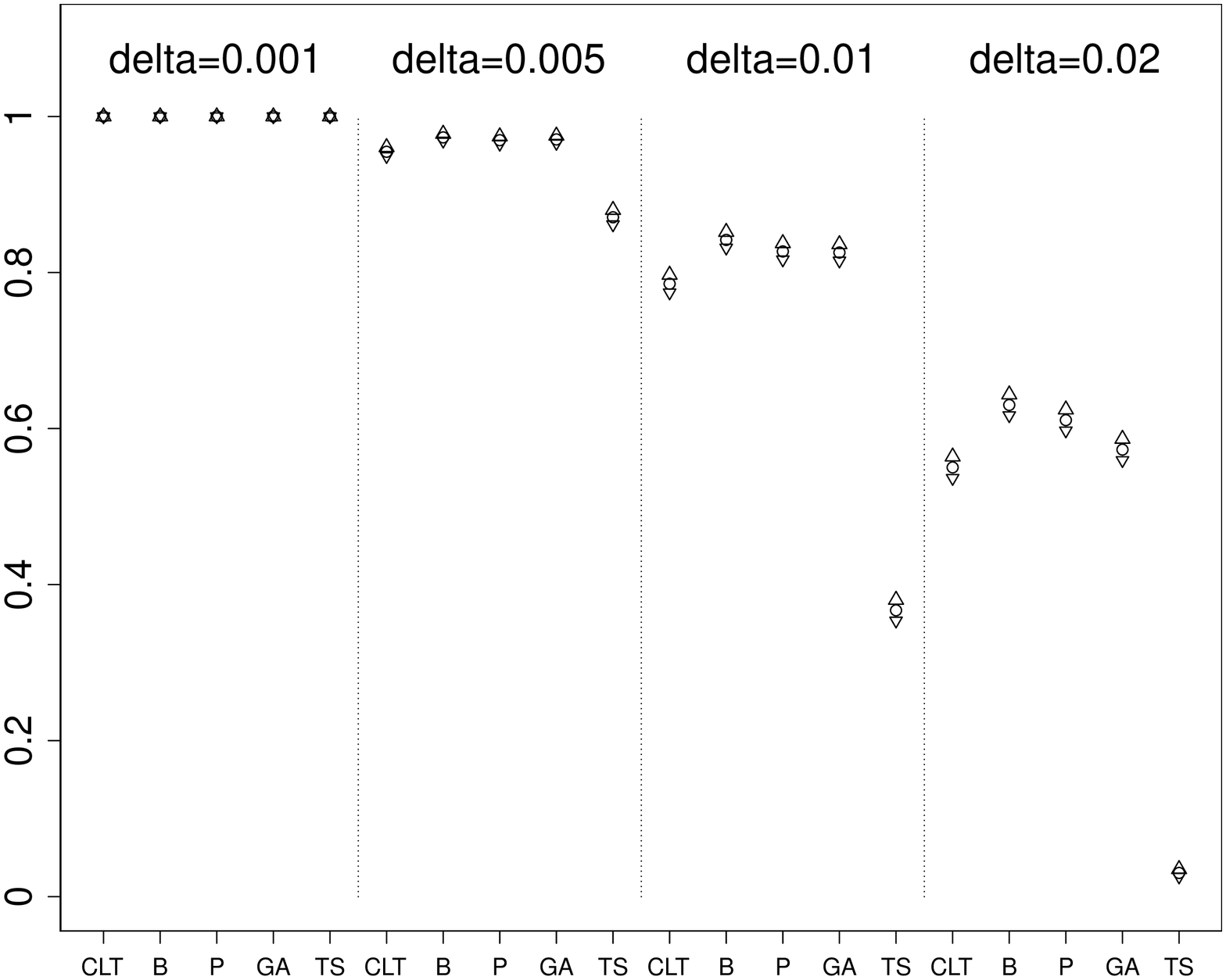}&
\includegraphics[angle=-90,scale=0.23]{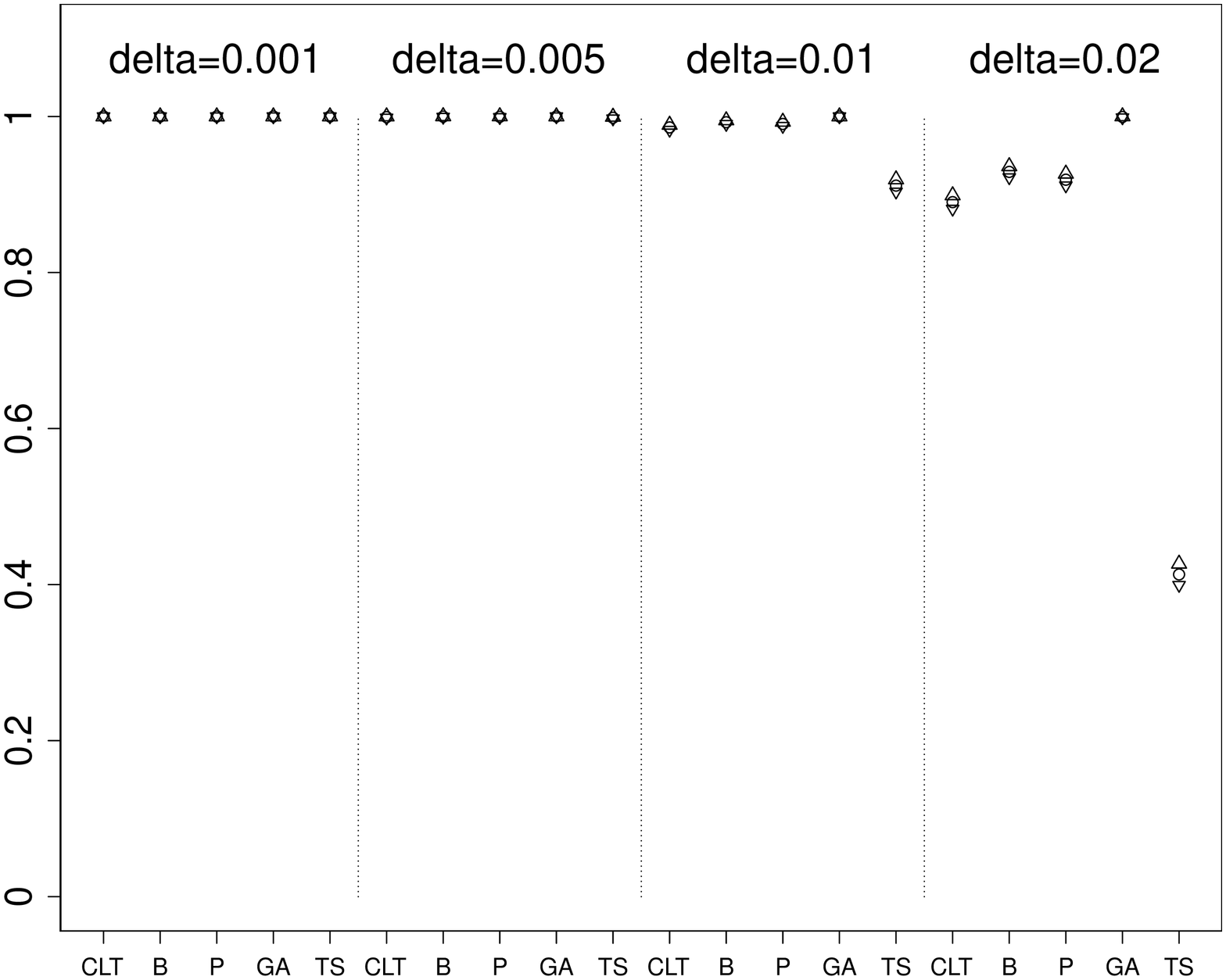}\vspace{5pt}\\
\hline \multicolumn{2}{c}{ } \vspace{-10pt}\\
\multicolumn{2}{c}{{\bf Exp. F}}\vspace{-8pt}\\
\multicolumn{2}{c}{\includegraphics[angle=-90,scale=0.23]{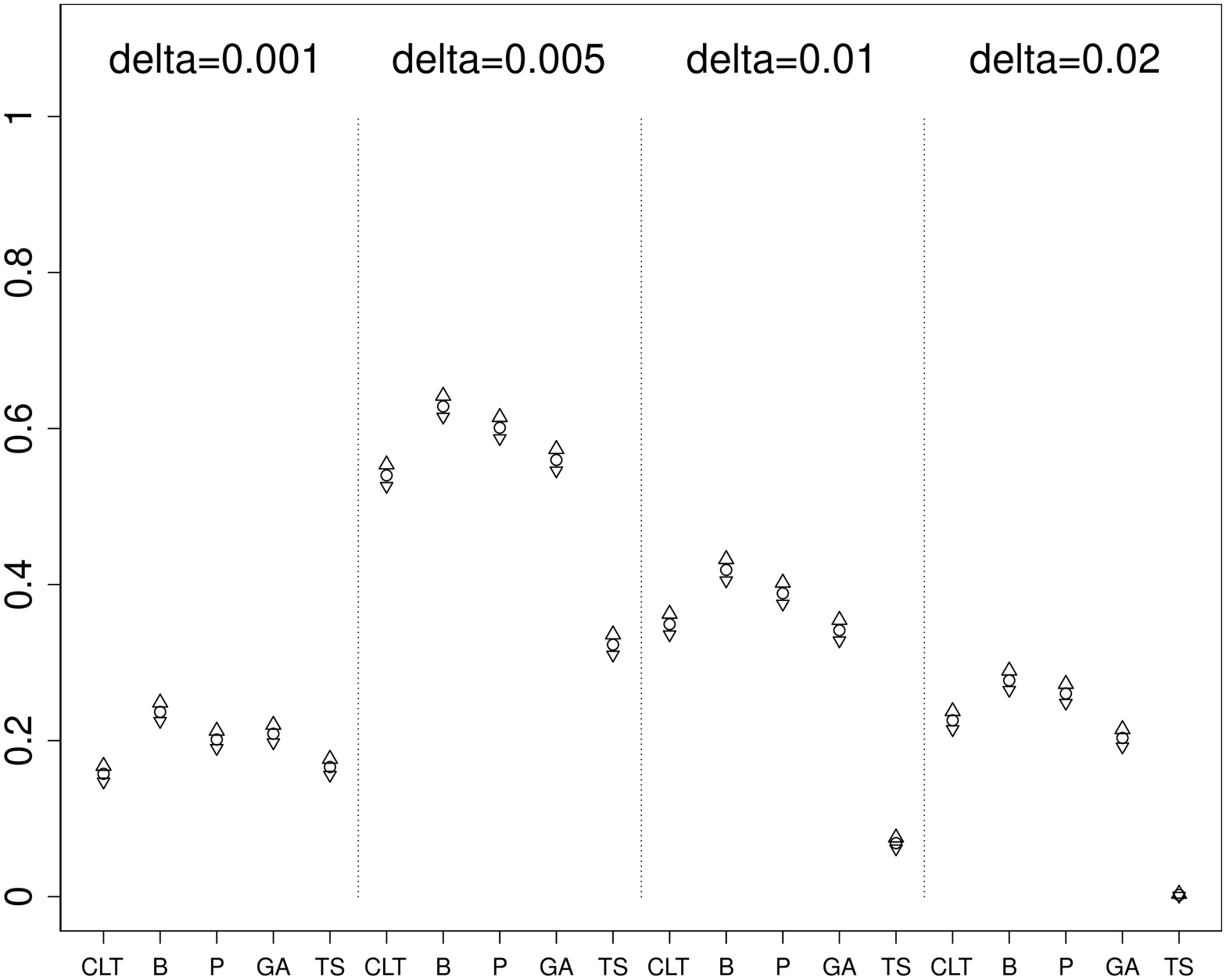}}
\end{tabular}
\caption{\label{puiss_delta} Estimated powers for different $\delta$. Same convention as in Figure \ref{nessai}.}
\end{figure}

\paragraph{Varying delay $\delta$}
We now investigate the impact of the choice for the delay $\delta$ by making $\delta$ vary in \ac{0.001, 0.005, 0.01, 0.02} for a fixed number of trials $n = 50$. 
The results for the sizes being similar to the previous study, only the estimated powers are presented in Figure \ref{puiss_delta}. 

On the top row of Figure \ref{puiss_delta}, the same process is injected in both coordinates: the coincidences are exact in the sense that they have no delay. Therefore, the best choice for the delay parameter $\delta$ is the smallest possible value:  the obtained power is $1$ for very small $\delta$'s (e.g. $\delta=0.001$) and then decreases as $\delta$ increases. On the contrary on the bottom row, it can be noticed that the highest power is for $\delta =  0.005$ which is the exact length of the interaction period~$u$. 
Once again, the (\texttt{TS}) method performs poorly, as does the (\texttt{CLT}) method. The three other methods seem to be quite equivalent except in the inhomogeneous case ({\bf Exp. E}) where the (\texttt{GA}) method has a power always equal to $1$, but at the price of an uncontrolled size.

\section{Conclusion}\label{conclusion}

In the present paper, we have introduced non-parametric independence tests between point processes based on $U$-statistics. The proposed critical values are obtained either by bootstrap or permutation approaches. We have shown that both methods share the same asymptotic properties under the null hypothesis as well as under the alternative. From a theoretical point of view, the main asymptotic results (Theorem \ref{Convergence} and Theorem~\ref{consistencyperm}) have almost the same flavor. However, there are additional assumptions in the permutation case which make the bootstrap results more general (despite the additional continuity assumption, which is very mild). From a more concrete point of view, it is acknowledged (see e.g. \cite{EfronTib}) that permutation should be preferred because of its very general non-asymptotic properties \eqref{nonasymp}. This is confirmed by the experimental study, where clearly permutation leads to a better first kind error rate control. However, both approaches perform much better than a naive procedure, based on a basic application of a Central Limit Theorem, when the number of observation is small. They also outperform existing procedures of the neuroscience literature, namely \cite{MTGAUE}, which  assumes the point processes to be homogeneous Poisson processes and the trial-shuffling procedures \cite{Pipa2003,Pipaet2003}, which are biased bootstrap variants applied on a non-centered  quantity.

One of the main open question with respect to the existing literature is whether our results can be extended to test statistics as $\sup_{h} U_{n,h}$. A first obstacle to this question lies in the nature of the observed random variables (point processes) and the fact that controlling such a supremum leads to controlling the whole $U$-process. This difficulty can probably be overcome, since the asymptotic Gaussian behavior of similar statistics has already been proved in general spaces under $(H_0)$ for product type kernels (see \cite{bouzebda}). The study of such behavior under $(H_1)$ is surely much more complex. A second obstacle comes from a more practical aspect. In neuroscience, and in the particular case of coincidence count, the use of $\sup_{\delta} U_{n,h_{\Nbc_\delta}}$ leads to the following fundamental problems. On the one hand, such a statistic may not be computable if $\delta$ varies in a too large space, typically $[0,1]$. On the other (more important) hand, neuroscientists are especially interested in the value of $\delta$ which leads to a rejection, since it actually provides the delay of interaction (see also Section \ref{simusec}). In this respect, our work in \cite{nous-appli} involves multiple testing aspects, which may answer this issue.

\section*{Acknowledgments}
% ---------------------------------------------------------------------------------------------------- %
We are grateful to both referees, whose discussions and comments allowed us to improve the present article. We  thank F. Grammont  for fruitful discussions. This work was granted access to the HPC and visualization resources of "Centre de Calcul Interactif" hosted by "Universit\'e Nice Sophia Antipolis".
This research was partly supported by the french Agence Nationale de la Recherche (ANR 2011 BS01 010 01 projet Calibration), by the PEPS BMI 2012-2013 {\it Estimation of dependence graphs for thalamo-cortical neurons and multivariate Hawkes processes} and by the interdisciplanary axis MTC-NSC of the University of Nice Sophia-Antipolis. The PhD grant of M. Albert is funded by the PACA french region. 
\medskip

\begin{center}
{\Large\textsc{Supplement}}
\end{center}
\textbf{Supplement: Technical Results and Proofs of "Bootstrap and Permutation tests of independence for Point Processes"}\\
This Supplement consists of all the proofs. It also contains some additional results about non-degeneracy and the empirical centering assumption.

\newpage

\begin{center}
{\Large \textsc{Technical Results and Proofs of \\
\medskip
"Bootstrap and permutation tests of \\
\medskip
independence for point processes"}}
\end{center}

The references of Equations, Theorems, Propositions, etc, that use only numbers such as (3.1) for instance, refer to the main article {\it Bootstrap and permutation tests of independence for point processes}.

\appendix
\section{Complete Proofs}

All along this section, $C$ and $C'$ denote positive constants, that may vary from one line to another one.

% ---------------------------------------------------------------------------------------------------- %
\subsection{Proof of Proposition \ref{Prop:non_deg_coinc}} 
\label{nondegcoinc}
% ---------------------------------------------------------------------------------------------------- %
We focus on the \emph{Coincidence case}. According to the comment following the definition of $\pa{\mc{A}_{non-deg}}$, $U_{n}(\X_n)$ is non-degenerate under $(H_0)$ if one can find some borelian set $\mc{B}$ of $\calX^2$ such that $P^1\otimes P^2 (\mc{B})>0$ and such that for all $x$ in $\mc{B}$, 
$\esp{h_{\varphi_{\delta}^{coinc}} (x,X) } \neq 0$, where $X$ has distribution $P^1\otimes P^2$. 

Consider $\mc{B}=\ac{\pa{\emptyset,\emptyset}}$.
Then $P^1\otimes P^2(\mc{B}) = P^1\pa{\ac{\emptyset}} P^2\pa{\ac{\emptyset}} >0$. \\
Moreover, as $\varphi_{\delta}^{coinc}\pa{\cdot,\emptyset}$ and $\varphi_{\delta}^{coinc}\pa{\emptyset,\cdot}$ are both the zero function, under $(H_0)$, 
$$\esp{h_{\varphi_{\delta}^{coinc}} \pa{(\emptyset,\emptyset),X} } = \frac12 \esp{\varphi_{\delta}^{coinc}\pa{X^1,X^2}}>0,$$ 
as $\varphi_{\delta}^{coinc}\pa{X^1,X^2}$ is non-negative and not almost surely null under $(H_0)$. 

\medskip
See also appendix \ref{Apdx:dege} for further results on the non-degeneracy of the $U$-statistic in more general cases.

% ---------------------------------------------------------------------------------------------------- %
\subsection{Proof of Proposition \ref{hinHw}} 
\label{contHw}
% ---------------------------------------------------------------------------------------------------- %
Consider $w: [0,1]^2\rightarrow \R$ a continuous integrable function.
Let us prove that $h=h_{\varphi^w}$ given by \eqref{defw} and \eqref{hphi} is continuous for the topology induced by $d$ (see \eqref{distance}).
Recall that for $x_1=(x_1^1,x_1^2)$ and $x_2=(x_2^1,x_2^2)$ in $(\calX^2)$, 
$$h_{\varphi^w}(x_1,x_2) = \frac{1}{2}\pa{\varphi^w (x_1^1,x_1^2) + \varphi^w (x_2^1,x_2^2) - \varphi^w (x_1^1,x_2^2) - \varphi^w (x_2^1,x_1^2)}.$$
The first step is to show that for each $i,j$ in \ac{1,2}, the projection defined  by
$$p_{i,j}:\fct{\pa{(\calX^2)^2,d}}{\pa{\calX^2,d_{\calX^2}}}{\pa{\pa{x_1^1,x_1^2},\pa{x_2^1,x_2^2}}}{\pa{x_i^1,x_j^2}},$$
is continuous. 
Let %$i$ and $j$ in \ac{1,2}, 
${\bf x}=\pa{\pa{x_1^2,x_1^2},\pa{x_2^1,x_2^2}}$ and ${\bf x}'=\pa{\pa{{x'}_1^1,{'}_1^2},\pa{{x'}_2^1,{x'}_2^2}}$ in $(\calX^2)^2$. Then, 
$$d_{\calX^2}\!\pa{p_{i,j}({\bf x}),p_{i,j}({\bf x}')} = d_{\calX^2}\!\pa{\pa{x_i^1,x_j^2},\pa{{x'}_i^1,{x'}_j^2}} \leq d\!\pa{{\bf x},{\bf x}'}.$$
Hence, $p_{i,j}$ is $1$-Lipschitz and therefore continuous.

\smallskip

The second step is to show that if $w$ is continuous on $\pa{[0,1]^2,\norm{\cdot}_{\infty}}$, with $\norm{(u,v)-(u',v')}_{\infty} = \max\ac{|u-u'|,|v-v'|}$, then $\varphi^w$ is also continuous. \\
Let $\varepsilon > 0$ and for $z$ in $\calX$, recall that $N_z$ is the counting process associated with $z$, defined by
$$N_z(t)=\int_0^1 \1{u\leq t} dN_{z}(u).$$ 
First notice that, $w$ being continuous on the compact set $[0,1]^2$, $w$ is uniformly continuous. 
Thus one can find some $\eta$ in $(0,1)$ such that, for all $(u,v)$, $(u',v')$ in $[0,1]^2$, 
\begin{equation}
\label{UCw}
\norm{(u,v)-(u',v')}_{\infty} \leq \eta \quad \mbox{implies} \quad \left|w(u,v)-w(u',v')\right|\leq \varepsilon.
\end{equation}
Consider such $\eta$. \\
Let $\ac{x_n}_{n\geq 0}$ be a sequence in $\calX^2$ such that $d_{\calX^2}\pa{x_n,x_0}\cv{n\to+\infty}0$ and let us show that 
$\varphi^w(x_n) \cv{n\to+\infty} \varphi^w(x_0)$.
There exists $n_0$ in $\N$ such that for all $n\geq n_0$, $d_{\calX^2}\pa{x_n,x_0}\leq \eta$. 
Then, for such $n$, by definition of $d_{\calX^2}$, we have that $d_\calD\!\pa{\!N_{x_n^1},N_{x_0^1}\!}\leq \eta$ and $d_\calD\pa{N_{x_n^2},N_{x_0^2}}\leq \eta$. Thus, by definition of $d_\calD$, 
$$ \exists \lambda_n^1 \in \Lambda \tq 
\left\{ \begin{array}{lr} \sup_{t\in [0,1]} \left|\lambda_n^1(t) - t\right| \leq \eta, & (1\mbox{-}i)\\
\sup_{t\in [0,1]}\left|N_{x_n^1}(t) - N_{x_0^1}\pa{\lambda_n^1(t)} \right| \leq \eta, & (1\mbox{-}ii)\end{array}\right. $$
%and 
$$ \exists \lambda_n^2 \in \Lambda \tq 
\left\{ \begin{array}{lr} \sup_{t\in [0,1]} \left|\lambda_n^2(t) - t\right| \leq \eta, & (2\mbox{-}i)\\
\sup_{t\in [0,1]}\left|N_{x_n^2}(t) - N_{x_0^2}\pa{\lambda_n^2(t)} \right| \leq \eta. & (2\mbox{-}ii)\end{array}\right. $$

In particular, as $\eta$ is chosen strictly smaller than 1 and as the $N_{x_n^j}$'s ($n\geq 0$, $j=1,2$) are counting processes with values in $\N$, 
$(1\mbox{-}ii)$ implies that $\forall t\in[0,1]$, $N_{x_n^1}(t) = N_{x_0^1}\pa{\lambda_n^1(t)}$ and thus, 
$$u_0\in x_0^1 \quad \Leftrightarrow \quad u_n=\lambda_n^1(u_0) \in x_n^1.$$ 
Similarly, $(2\mbox{-}ii)$ implies that
$$v_0\in x_0^2 \quad \Leftrightarrow \quad v_n=\lambda_n^2(v_0) \in x_n^2.$$ 
Therefore, 
\begin{eqnarray*}
\varphi^w(x_n) &=& \iint_{[0,1]^2} w(u,v) dN_{x_n^1}(u) dN_{x_n^2}(v) \\
&=& \sum_{(u_n,v_n)\in x_n^1\times x_n^2} w(u_n,v_n) \\
&=& \sum_{(u_0,v_0)\in x_0^1\times x_0^2} w\pa{\lambda_n^1(u_0),\lambda_n^2(v_0)}. 
\end{eqnarray*}
Hence, 
$$\left|\varphi^w(x_n)-\varphi^w(x_0)\right| \leq  \sum_{(u_0,v_0)\in x_0^1\times x_0^2} \left| w\pa{\lambda_n^1(u_0),\lambda_n^2(v_0)} - w\pa{u_0,v_0}\right|.$$
Yet, by $(1\mbox{-}i)$ and $(2\mbox{-}i)$, for each $(u_0,v_0)$ in $x_0^1\times x_0^2$, we have 
$$\norm{\pa{\lambda_n^1(u_0),\lambda_n^2(v_0)}-\pa{u_0,v_0}}_\infty \leq \eta,$$ and thus, applying \eqref{UCw}, we obtain
$$\left|\varphi^w(x_n)-\varphi^w(x_0)\right| \leq  \#x_0^1 \#x_0^2\varepsilon,$$
and this for all $n\geq n_0$, which ends the proof of Proposition \ref{hinHw}.

% ---------------------------------------------------------------------------------------------------- %
\subsection{Proof of Proposition \ref{hinH}} 
\label{contcoinc}
% ---------------------------------------------------------------------------------------------------- %

Let us prove that in the \emph{Coincidence case}, the kernel $h=h_{\varphi_\delta^{coinc}}$ given by \eqref{defnbc} and \eqref{hphicoinc} is continuous for the topology induced by the metric $d$ (defined in \eqref{distance}) in any $\pa{x_0,y_0}$ in $\mathcal{C}_\delta$ satisfying 
\begin{equation*}%\label{1indep2}
\pa{\ac{x_0^1}\cup\ac{y_0^1}} \cap \pa{\ac{x_0^2 \pm \delta}\cup\ac{y_0^2 \pm\delta}} = \emptyset.
\end{equation*}
As in the proof of Proposition \ref{hinHw}, denote by $N_z$ the counting process associated with $z$:
$$N_z(t)=\int_0^1 \1{u\leq t} dN_{z}(u).$$
Consider a sequence $\ac{\pa{x_n, y_n}}_{n\in\N}$ of elements in $\calX^2\times\calX^2$, where $x_n = \pa{x_n^1,x_n^2}$ and $y_n = \pa{y_n^1,y_n^2}$ such that $d \pa{\pa{x_n,y_n},\pa{x_0,y_0}} \cv{n\to+\infty} 0$ and $\pa{x_0,y_0}$ belongs to $\mathcal{C}_\delta$.\\
We want to show that $\left| h\pa{x_n,y_n} - h\pa{x_0,y_0}\right| \cv{n\to+\infty} 0$.\\
Since $\pa{x_0,y_0}$ is in $\mathcal{C}_\delta$, for any $t_0$ in $\ac{x_0^2\pm\delta}\cup\ac{y_0^2\pm\delta}$, 
$t_0\notin x_0^1$, which means that $N_{x_0^1}$ is continuous in $t_0$ and therefore constant in a neighborhood:
$$\exists\eta_{t_0} >0 \tq \forall t\in [0,1],\quad |t-t_0|\leq\eta_{t_0} \Rightarrow N_{x_0^1}(t) =N_{x_0^1}(t_0).$$
As $\ac{x_0^2\pm\delta}\cup\ac{y_0^2\pm\delta}$ is finite, $\eta_{x_0^1} = \min_{t_0\in \ac{x_0^2\pm\delta}\cup\ac{y_0^2\pm\delta}} \eta_{t_0} >0$ is well defined, and satisfies 
$$\forall u\in \ac{x_0^2\pm\delta}\cup \ac{y_0^2\pm\delta},\ \forall t\in [0,1],\quad |t-u|\leq\eta_{x_0^1} \Rightarrow N_{x_0^1}(t) =N_{x_0^1}(u).$$
By the same argument using continuity of $N_{y_0^1}$ over $\ac{x_0^2\pm\delta}\cup\ac{y_0^2\pm\delta}$, one can find $\eta_{y_0^1}>0$ such that 
$$\forall u\in \ac{x_0^2\pm\delta}\cup \ac{y_0^2\pm\delta},\ \forall t\in [0,1],\quad |t-u|\leq\eta_{y_0^1} \Rightarrow N_{y_0^1}(t) = N_{y_0^1}(u).$$
Since $(x_0,y_0)\in\mathcal{C}_\delta\Leftrightarrow \pa{\ac{x_0^2}\cup\ac{y_0^2}} \cap \pa{\ac{x_0^1 \pm \delta}\cup\ac{y_0^1 \pm\delta}} = \emptyset$, one can construct $\eta_{x_0^2}$ and $\eta_{y_0^2}$ satisfying 
$$\forall u\in \ac{x_0^1\pm\delta}\cup \ac{y_0^1\pm\delta},\ \forall t\in [0,1],
\left\{\begin{array}{l}
|t-u|\leq\eta_{x_0^2}  \Rightarrow N_{x_0^2}(t)=N_{x_0^2}(u) \\ 
|t-u|\leq\eta_{y_0^2} \Rightarrow N_{y_0^2}(t)=N_{y_0^2}(u) \end{array}\right. .$$
Finally, if $\eta =\min\ac{\eta_{x_0^1},\eta_{y_0^1},\eta_{x_0^2},\eta_{y_0^2}}>0$,
\begin{equation}\label{etaconst1}
\forall s\in\ac{x_0^2 \pm \delta}\cup\ac{y_0^2\pm\delta},\forall t\in [0,1],
|t-s|\leq\eta \Rightarrow
\left\{\begin{array}{l}
N_{x_0^1}(t)= N_{x_0^1}(s)\\ 
N_{y_0^1}(t) = N_{y_0^1}(s)\end{array}\right. ,
\end{equation}
\begin{equation}\label{etaconst2}
\forall s\in\ac{x_0^1 \pm \delta}\cup\ac{y_0^1\pm\delta}, \forall t\in [0,1],
|t-s|\leq\eta\Rightarrow
\left\{\begin{array}{l}
N_{x_0^2}(t) = N_{x_0^2}(s) \\ 
N_{y_0^2}(t) = N_{y_0^2}(s) \end{array}\right. .
\end{equation}
As $d\pa{\pa{x_n,y_n},\pa{x_0,y_0}} \to_{n\to+\infty} 0$, there exists $n_0\geq 0$ such that for $n\geq n_0$, 
$d\pa{\pa{x_n,y_n},\pa{x_0,y_0}}\leq \eta /4.$
From the definition of $d$, we deduce that
$$ \exists \lambda_n^1 \in \Lambda \tq 
\left\{ \begin{array}{lr} \sup_{t\in [0,1]} \left|\lambda_n^1(t) - t\right| \leq \frac{\eta}{4} & (1\mbox{-}i)\\
\sup_{t\in [0,1]}\left|N_{x_n^1}(t) - N_{x_0^1}\pa{\lambda_n^1(t)} \right| \leq \frac{\eta}{4} & (1\mbox{-}ii)\end{array}\right., $$
and $$ \exists \lambda_n^2 \in \Lambda \tq 
\left\{ \begin{array}{lr} \sup_{t\in [0,1]} \left|\lambda_n^2(t) - t\right| \leq \frac{\eta}{4} & (2\mbox{-}i)\\
\sup_{t\in [0,1]}\left|N_{x_n^2}(t) - N_{x_0^2}\pa{\lambda_n^2(t)} \right| \leq \frac{\eta}{4} & (2\mbox{-}ii)\end{array}\right. .$$
Notice that similar results occur for $y_n$ and $y_0$, but there are not detailed here since we do not use them explicitly.\\
By definition of $h$, 
\begin{align}\label{controlh}
 &h(x_n,y_n) - h\pa{x_0,y_0} \\
\nonumber &= \frac12\iint \1{|u-v|\leq\delta} \ac{dN_{x_n^1}dN_{x_n^2} + dN_{y_n^1}dN_{y_n^2} - dN_{x_n^1}dN_{y_n^2} - dN_{y_n^1}dN_{x_n^2}}(u,v) \\
 \nonumber&\quad -\frac12\iint \1{|u-v|\leq\delta} \ac{dN_{x_0^1}dN_{x_0^2} + dN_{y_0^1}dN_{y_0^2} - dN_{x_0^1}dN_{y_0^2} - dN_{y_0^1}dN_{x_0^2}}(u,v) \\
\nonumber &=\frac12 \iint \1{|u-v|\leq\delta} \Big(dN_{x_n^1}(u)\pa{dN_{x_n^2}-dN_{x_0^2}}(v) + dN_{y_n^1}(u)\pa{dN_{y_n^2}-dN_{y_0^2}}(v) \\
\nonumber &\quad - dN_{x_n^1}(u)\pa{dN_{y_n^2}-dN_{y_0^2}}(v) - dN_{y_n^1}(u)\pa{dN_{x_n^2}-dN_{x_0^2}}(v) \\
\nonumber &\quad + \pa{dN_{x_n^1}-dN_{x_0^1}}(u)\ dN_{x_0^2}(v) + \pa{dN_{y_n^1}-dN_{y_0^1}}(u)\ dN_{y_0^2}(v) \\
\nonumber &\quad -\pa{dN_{x_n^1}-dN_{x_0^1}}(u)\ dN_{y_0^2}(v) + \pa{dN_{y_n^1}-dN_{y_0^1}}(u)\ dN_{x_0^2}(v) \Big).
\end{align}
By symmetry of the problem, we just need to study the terms 
$$A_n=\iint \1{|u-v|\leq\delta}\pa{dN_{x_n^1}-dN_{x_0^1}}(u)\ dN_{x_0^2}(v),$$ and
$$B_n=\iint \1{|u-v|\leq\delta} dN_{x_n^1}(u)\pa{dN_{x_n^2}-dN_{x_0^2}}(v).$$

% --------------------- study of A_n --------------------- %

\paragraph{Study of $A_n$}

\begin{eqnarray*}
A_n &=& \iint \1{|u-v|\leq\delta}\pa{dN_{x_n^1}-dN_{x_0^1}}(u)\ dN_{x_0^2}(v) \\
&=& \iint \1{u\leq v+\delta} \pa{dN_{x_n^1}-dN_{x_0^1}}(u)\ dN_{x_0^2}(v)\\
&& - \iint \1{u< v-\delta} \pa{dN_{x_n^1}-dN_{x_0^1}}(u)\ dN_{x_0^2}(v).
\end{eqnarray*}
We have that
 \begin{align*}
 \Big|\iint \1{u\leq v+\delta} \Big(dN_{x_n^1}-&dN_{x_0^1}\Big)(u)\ dN_{x_0^2}(v) \Big|\\
 &= \abs{\int\pa{N_{x_n^1}(v+\delta)-N_{x_0^1}(v+\delta)}dN_{x_0^2}(v)}  \\
 &\leq \sum_{T\in x_0^2} \abs{N_{x_n^1}(T+\delta)-N_{x_0^1}(T+\delta)}\\
 &\leq \sum_{T\in x_0^2} \left| N_{x_n^1}(T+\delta)-N_{x_0^1}\pa{\lambda_n^1(T+\delta)} \right| \\
 &+\sum_{T\in x_0^2}\left| N_{x_0^1}\pa{\lambda_n^1(T+\delta)}-N_{x_0^1}(T+\delta) \right|.
\end{align*}
Now, using the notation $N_{x_i^1}^-(t)=\int \1{u< t} dN_{x_i^1}(u)$,
$$\abs{ \iint \1{u< v-\delta} \pa{dN_{x_n^1}-dN_{x_0^1}}(u)\ dN_{x_0^2}(v)} 
 \leq \sum_{T\in x_0^2} \abs{N_{x_n^1}^-(T-\delta)-N_{x_0^1}^-(T-\delta)}.$$
Therefore,
 \begin{eqnarray}\label{controlAn}
 |A_n|&\leq &\sum_{T\in x_0^2} \Bigg(\left| N_{x_n^1}(T+\delta)-N_{x_0^1}\pa{\lambda_n^1(T+\delta)} \right|\\
\nonumber&&+\left| N_{x_0^1}\pa{\lambda_n^1(T+\delta)}-N_{x_0^1}(T+\delta) \right|\\
\nonumber &&+  \abs{N_{x_n^1}^-(T-\delta)-N_{x_0^1}^-(T-\delta)}\Bigg).
\end{eqnarray}
Let us study individually each term in the sum.\\
Fix $T$ in $x_0^2$.  By $(1\mbox{-}ii)$,
\begin{equation}\label{An1}
\left| N_{x_n^1}(T+\delta)-N_{x_0^1}\pa{\lambda_n^1(T+\delta)} \right| \leq \frac{\eta}{4} \leq \varepsilon.
\end{equation}
From $(1\mbox{-}i)$, one has $|\lambda_n^1(T+\delta)-(T+\delta)|\leq\frac{\eta}{2} \leq \eta$ which, with \eqref{etaconst1}, implies 
\begin{equation}\label{An2}
\left|N_{x_0^1}\pa{\lambda_n^1(T+\delta)}-N_{x_0^1}(T+\delta)\right|=0.
\end{equation}
As $N_{x_n^1}^-(T-\delta) =\lim\limits_{\substack{u\to T-\delta \\ u<T-\delta}} N_{x_n^1}(u)$, there exists 
$u_T$ in $[T-\delta-\eta/4, T-\delta[$ such that 
$$\left| N_{x_n^1}^-(T-\delta)-N_{x_n^1}\pa{u_T} \right| \leq \varepsilon,$$ so
\begin{multline}\label{decoup-}
 \left| N_{x_n^1}^-(T-\delta)-N_{x_0^1}^-(T-\delta) \right|\leq\varepsilon+ \left| N_{x_n^1}\pa{u_T}-N_{x_0^1}\pa{\lambda_n^1(u_T)} \right|\\
+\left| N_{x_0^1}\pa{\lambda_n^1(u_T)}-N_{x_0^1}^-(T-\delta) \right|.
\end{multline}
From $(1\mbox{-}ii)$, one has $\left| N_{x_n^1}\pa{u_T}-N_{x_0^1}\pa{\lambda_n^1(u_T)} \right|\leq {\eta}/{4} \leq \varepsilon$.\\
Then, by continuity of $N_{x_0^1}$ at $T-\delta$, first remark that $N_{x_0^1}^-(T-\delta)=N_{x_0^1}(T-\delta)$. Moreover, by $(1\mbox{-}i)$ and construction of $u_T$,
 $$\left|\lambda_n^1(u_T)-(T-\delta)\right| \leq \left|\lambda_n^1(u_T)-u_T\right| + \left|u_T-(T-\delta)\right|\leq \frac{\eta}{4}+\frac{\eta}{4} < \eta,$$
 hence, using \eqref{etaconst1},
 $\left| N_{x_0^1}\pa{\lambda_n^1(u_T)}-N_{x_0^1}^-(T-\delta) \right| =0.
$ So finally, \eqref{decoup-} gives
\begin{equation}\label{An3}
 \left| N_{x_n^1}^-(T-\delta)-N_{x_0^1}^-(T-\delta) \right| \leq 2\varepsilon.
\end{equation}
Combining \eqref{controlAn}, \eqref{An1}, \eqref{An2}, and \eqref{An3}, we obtain that for any $n\geq n_0$:
\begin{equation}\label{majoAn}
|A_n|\leq 3\varepsilon \#x_0^2.
\end{equation}

\paragraph{Study of $B_n$}

Recall that $B_n=\iint \1{|u-v|\leq\delta} dN_{x_n^1}(u)\pa{dN_{x_n^2}-dN_{x_0^2}}(v)$. \\
As for $A_n$,  $B_n$ is upper bounded by a sum of several terms, that we study separately.
$$B_n= \sum_{T\in x_n^1} \pa{N_{x_n^2}\pa{T+\delta}-N_{x_0^2}\pa{T+\delta}} -\sum_{T\in x_n^1} \pa{N_{x_n^2}^-\pa{T-\delta}-N_{x_0^2}^-\pa{T-\delta}}.$$
So
\begin{equation}\label{controlBn}
 B_n \leq |B_{n,1}|+|B_{n,2}|+|B_{n,3}|+|B_{n,4}|,
\end{equation}
with
\begin{eqnarray*}
B_{n,1}&=&\sum_{T\in x_n^1} \pa{ N_{x_n^2}\pa{T+\delta}-N_{x_0^2}\pa{\lambda_n^2\pa{T+\delta}} },\\
B_{n,2}&=& \sum_{T\in x_n^1} \pa{N_{x_0^2}\pa{\lambda_n^2\pa{T+\delta}} - N_{x_0^2}\pa{T+\delta}},\\
B_{n,3}&=&\sum_{T\in x_n^1}  \pa{N_{x_n^2}^-\pa{T-\delta}-N_{x_0^2}^-\pa{\lambda_n^2\pa{T-\delta}}},\\
B_{n,4}&=&\sum_{T\in x_n^1} \pa{N_{x_0^2}^-\pa{\lambda_n^2\pa{T-\delta}} + N_{x_0^2}^-\pa{T-\delta}}.
\end{eqnarray*}
The control of $B_n$ is quite similar to the one of $A_n$ except that the sums are over $T$ in $x_n^1$ instead of $T$ in $x_0^1$, which prevents us to use \eqref{etaconst2} and \eqref{etaconst1} directly.

\smallskip

\noindent Control of $B_{n,1}$. Due to $(2\mbox{-}ii)$, 
$\left| N_{x_n^2}\pa{T+\delta}-N_{x_0^2}\pa{\lambda_n^2\pa{T+\delta}} \right|\leq \varepsilon$, so
\begin{equation}\label{Bn1}
|B_{n,1}|\leq \varepsilon \#x_n^1.
\end{equation}

\smallskip

\noindent Control of $B_{n,2}$. One can easily see that
\begin{eqnarray*}
B_{n,2}&=& \iint \pa{\1{v\leq \lambda_n^2\pa{u+\delta}} - \1 {v\leq u+\delta}} dN_{x_0^2}(v)\ dN_{x_n^1}(u) \\
&=& \iint \cro{\pa{1- \1{u<\pa{\lambda_n^2}^{-1}(v) - \delta}} - \pa{1-\1{u< v - \delta}}}dN_{x_0^2}(v)\ dN_{x_n^1}(u) \\
&=& \sum_{T\in x_0^2} \pa{N_{x_n^1}^-(T-\delta) - N_{x_n^1}^-(\pa{\lambda_n^2}^{-1}(T)-\delta)}.
\end{eqnarray*}
Fix now $T$ in $x_0^2$. 
\begin{multline*}
\left|N_{x_n^1}^-(T-\delta) - N_{x_n^1}^-(\pa{\lambda_n^2}^{-1}(T)-\delta) \right|
\leq  \left|N_{x_n^1}^-(T-\delta) - N_{x_0^1}(T-\delta) \right|\\
+ \left| N_{x_0^1}(T-\delta) - N_{x_n^1}^-(\pa{\lambda_n^2}^{-1}(T)-\delta) \right|.
\end{multline*}
As shown in \eqref{decoup-}, $\left|N_{x_n^1}^-(T-\delta) - N_{x_0^1}(T-\delta) \right|\leq 2\varepsilon$.\\
Furthermore, 
take $v_T$ in $\left[\pa{\lambda_n^2}^{-1}(T)-\delta - \eta/4, \pa{\lambda_n^2}^{-1}(T)-\delta\right[$  such that 
$$\left| N_{x_n^1}^-(\pa{\lambda_n^2}^{-1}(T)-\delta) -N_{x_n^1}(v_T)\right|\leq \varepsilon.$$
 So,
\begin{multline*}
\left| N_{x_0^1}(T-\delta) - N_{x_n^1}^-(\pa{\lambda_n^2}^{-1}(T)-\delta) \right| 
\leq \varepsilon+ \left|N_{x_n^1}(v_T)-N_{x_0^1}\pa{\lambda_n^1(v_T)}\right|\\
+\left|N_{x_0^1}\pa{\lambda_n^1(v_T)} -  N_{x_0^1}(T-\delta) \right|.\end{multline*}
By construction of $v_T$ and $\lambda_n^1$ (see $(1\mbox{-}ii)$), $\left|N_{x_n^1}(v_T)-N_{x_0^1}\pa{\lambda_n^1(v_T)}\right| \leq \varepsilon$.\\
Because of \eqref{etaconst1} which is true as 
$$\left|\lambda_n^1(v_T) - (T-\delta)\right| \leq |\lambda_n^1(v_T) - v_T| + |v_T - (T-\delta)|\leq \frac{\eta}{4}+\frac{\eta}{4}<\eta$$ by $(1\mbox{-}i)$,
$\left|N_{x_0^1}\pa{\lambda_n^1(v_T)} -  N_{x_0^1}(T-\delta) \right|=0.$ Hence,
$$\left| N_{x_0^1}(T-\delta) - N_{x_n^1}^-(\pa{\lambda_n^2}^{-1}(T)-\delta) \right|\leq 2\varepsilon.$$ 
Finally, 
$$\left|N_{x_n^1}^-(T-\delta) - N_{x_n^1}^-(\pa{\lambda_n^2}^{-1}(T)-\delta) \right| \leq 4\varepsilon,$$ 
and 
\begin{equation}\label{Bn2}
|B_{n,2}| \leq 4\varepsilon \#x_0^2.
\end{equation}
Control of $B_{n,3}$. First, for all $T$ in $x_n^1$, we find some $\nu_{n,T}$ in $]0, \eta/4]$ such that 
$$\forall u \in [T-\delta - \nu_{n,T}, T-\delta[,\quad \left|N_{x_n^2}^-\pa{T-\delta} - N_{x_n^2}\pa{u}\right|\leq \varepsilon.$$ Setting $\nu_n=\min_{T\in x_n^1} \nu_{n,T}$,
\begin{eqnarray*}
\left|B_{n,3}\right|& \leq& \sum_{T\in x_n^1} \left|N_{x_n^2}^-\pa{T-\delta} - N_{x_n^2}\pa{T-\delta-\nu_n}\right|\\
&& +  \sum_{T\in x_n^1} \left|N_{x_n^2}\pa{T-\delta-\nu_n} - N_{x_0^2}\pa{\lambda_n^2\pa{T-\delta-\nu_n}}\right| \\
&  &+  \left|\sum_{T\in x_n^1} \pa{N_{x_0^2}\pa{\lambda_n^2\pa{T-\delta-\nu_n}}-N_{x_0^2}^-\pa{\lambda_n^2\pa{T-\delta}}}\right|.
\end{eqnarray*}
For each $T$ in $x_n^1$, $\left|N_{x_n^2}^-\pa{T-\delta} - N_{x_n^2}\pa{T-\delta-\nu_n}\right| \leq \varepsilon $ and\\
$\left|N_{x_n^2}\pa{T-\delta-\nu_n} - N_{x_0^2}\pa{\lambda_n^2\pa{T-\delta-\nu_n}}\right| \leq \varepsilon $ by $(2\mbox{-}ii)$. Therefore,
$$\left|B_{n,3}\right| \leq2\varepsilon \#x_{n}^1+  \left|\sum_{T\in x_n^1} \pa{N_{x_0^2}\pa{\lambda_n^2\pa{T-\delta-\nu_n}}-N_{x_0^2}^-\pa{\lambda_n^2\pa{T-\delta}}}\right|.$$
Now,
\begin{align*}
\sum_{T\in x_n^1} \Big(N_{x_0^2}\Big(\lambda_n^2&\pa{T-\delta-\nu_n}\Big)-N_{x_0^2}^-\pa{\lambda_n^2\pa{T-\delta}} \Big)\\
&=\iint \1{v\leq \lambda_n^2\pa{u-\delta-\nu_n}} - \1{v<\lambda_n^2\pa{u-\delta}} dN_{X_n^1}(u)\ dN_{X_0^2}(v) \\
&= \sum_{T\in x_0^2}\pa{N_{x_n^1}\pa{\pa{\lambda_n^2}^{-1}(T)+\delta} - N_{x_n^1}^-\pa{\pa{\lambda_n^2}^{-1}(T)+\delta+\nu_n}}.
\end{align*}
For each $T$ in $x_0^2$,
\begin{align*}
\Big|N_{x_n^1}&\Big(\pa{\lambda_n^2}^{-1}(T)+\delta\Big)- N_{x_n^1}^-\pa{\pa{\lambda_n^2}^{-1}(T)+\delta+\nu_n}\Big|\\
&\leq  \left|N_{x_n^1}\pa{\pa{\lambda_n^2}^{-1}(T)+\delta} - N_{x_0^1}\pa{\lambda_n^1\pa{\pa{\lambda_n^2}^{-1}(T)+\delta}} \right| \\
&\quad  + \left|N_{x_0^1}\pa{\lambda_n^1\pa{\pa{\lambda_n^2}^{-1}(T)+\delta}} - N_{x_0^1}\pa{\lambda_n^1\pa{\pa{\lambda_n^2}^{-1}(T)+\delta+\nu_n}}\right| \\
& \quad +  \left| N_{x_0^1}\pa{\lambda_n^1\pa{\pa{\lambda_n^2}^{-1}(T)+\delta+\nu_n}}- N_{x_n^1}^-\pa{\pa{\lambda_n^2}^{-1}(T)+\delta+\nu_n}\right|\\
&\leq 2\varepsilon +\left| N_{x_0^1}\pa{\lambda_n^1\pa{\pa{\lambda_n^2}^{-1}(T)+\delta+\nu_n}}- N_{x_n^1}^-\pa{\pa{\lambda_n^2}^{-1}(T)+\delta+\nu_n}\right|,
\end{align*}
where the last line comes from $(1\mbox{-}ii)$, and \eqref{etaconst1}.\\
We now find some $w_T$ in $\left[\pa{\lambda_n^2}^{-1}(T)+\delta+\nu_n-\eta/4\ ,\ \pa{\lambda_n^2}^{-1}(T)+\delta+\nu_n\right[$ such that 
$$\left| N_{x_n^1}^-\pa{\pa{\lambda_n^2}^{-1}(T)+\delta+\nu_n}- N_{x_n^1}\pa{w_T}\right|\leq\varepsilon,$$ 
so
\begin{multline*}
 \left| N_{x_0^1}\pa{\lambda_n^1\pa{\pa{\lambda_n^2}^{-1}(T)+\delta+\nu_n}}- N_{x_n^1}^-\pa{\pa{\lambda_n^2}^{-1}(T)+\delta+\nu_n}\right|\\
 \leq \left| N_{x_0^1}\pa{\lambda_n^1\pa{\pa{\lambda_n^2}^{-1}(T)+\delta+\nu_n}}- N_{x_0^1}\pa{\lambda_n^1\pa{w_T}}\right|\\
  + \left| N_{x_0^1}\pa{\lambda_n^1\pa{w_T}} - N_{x_n^1}\pa{w_T} \right| +\varepsilon. 
\end{multline*}
From $(1\mbox{-}ii)$, we deduce that $\left| N_{x_0^1}\pa{\lambda_n^1\pa{w_T}} - N_{x_n^1}\pa{w_T} \right|\leq \varepsilon.$ Due to \eqref{etaconst1}, $(1\mbox{-}i)$, and the construction of $w_T$,
$$\left|\pa{\lambda_n^1\pa{\pa{\lambda_n^2}^{-1}(T)+\delta+\nu_n}} - (T-\delta)\right| \leq\frac{3\eta}{4} < \eta,$$
and 
$$\left|\pa{\lambda_n^1\pa{w_T}} - (T-\delta)\right| 
\leq \left|\pa{\lambda_n^1\pa{w_T} - w_T}\right|+\left|w_T - (T-\delta)\right|<\eta.$$So
$\left| N_{x_0^1}\pa{\lambda_n^1\pa{\pa{\lambda_n^2}^{-1}(T)+\delta+\nu_n}}- N_{x_0^1}\pa{\lambda_n^1\pa{w_T}}\right|=0.$
As a consequence,
$$\Big|N_{x_n^1}\Big(\pa{\lambda_n^2}^{-1}(T)+\delta\Big)- N_{x_n^1}^-\pa{\pa{\lambda_n^2}^{-1}(T)+\delta+\nu_n}\Big|\leq 4\varepsilon,$$
 and
 \begin{equation}\label{Bn3}
\left|B_{n,3}\right| \leq2\varepsilon \#x_{n}^1+  4\varepsilon  \#x_0^2.
\end{equation}

\smallskip

\noindent Control of $B_{n,4}$.
 \begin{eqnarray*}
B_{n,4}
&=& \iint \pa{\1{v<\lambda_n^2\pa{u-\delta}} - \1{v < u-\delta} }dN_{x_0^2}(v)\ dN_{x_n^1}(u) \\
&=& \sum_{T\in x_0^2} \pa{N_{x_n^1}(T+\delta) - N_{x_n^1}\pa{\pa{\lambda_n^2}^{-1}(T)+\delta}}.
\end{eqnarray*}
Let us fix  $T$ in $x_0^2$. We have 
\begin{align*}
\Big|N_{x_n^1}(T+\delta) &- N_{x_n^1}\pa{\pa{\lambda_n^2}^{-1}(T)+\delta}\Big|\\
&\leq  \left| N_{x_n^1}(T+\delta) - N_{x_0^1}\pa{\lambda_n^1\pa{T+\delta}}\right| \\
& \quad +   \left| N_{x_0^1}\pa{\lambda_n^1\pa{T+\delta}} - N_{x_0^1}\pa{\lambda_n^1\pa{\pa{\lambda_n^2}^{-1}(T)+\delta}}\right| \\
& \quad +   \left| N_{x_0^1}\pa{\lambda_n^1\pa{\pa{\lambda_n^2}^{-1}(T)+\delta}} - N_{x_n^1}\pa{\pa{\lambda_n^2}^{-1}(T)+\delta}\right|.
\end{align*}
The first and the last terms are upper bounded by $\varepsilon$ due to $(1\mbox{-}ii)$.
Furthermore, since $N_{x_0^1}\pa{\lambda_n^1\pa{\pa{\lambda_n^2}^{-1}(T)+\delta}} = N_{x_0^1}\pa{T+\delta} = N_{x_0^1}\pa{\lambda_n^1\pa{T+\delta}}$
by applying \eqref{etaconst1} and using $(1\mbox{-}i)$ and $(2\mbox{-}i)$,
$$\left| N_{x_0^1}\pa{\lambda_n^1\pa{T+\delta}} - N_{x_0^1}\pa{\lambda_n^1\pa{\pa{\lambda_n^2}^{-1}(T)+\delta}}\right|=0.$$
So finally, 
\begin{equation}\label{Bn4}
\left|B_{n,4}\right| \leq 2\varepsilon \#{x_0^2}.
\end{equation}

\smallskip

\noindent Combining \eqref{controlBn}, \eqref{Bn1}, \eqref{Bn2}, \eqref{Bn3}, and \eqref{Bn4}, we can conclude that
\begin{equation}\label{majoBn}
|B_n|\leq 3\varepsilon \#x_{n}^1+  10\varepsilon  \#x_0^2.
\end{equation}
We now just remark that $\pa{\#{x_n^1}}_{n\geq n_0}$ is bounded because it converges to $\#{x_0^1}$. 
Indeed, since $\#{x_n^1}=N_{x_n^1}(1)$, $\#{x_0^1}=N_{x_0^1}(1)$ and for every $n$, $\lambda_n^1(1) = 1$, 
\begin{eqnarray*}
\left|\#{x_n^1} - \#{x_0^1} \right| &=& \left| N_{x_n^1}(1) - N_{x_0^1}(1)\right| \\
									&=& \left| N_{x_n^1}(1) - N_{x_0^1}\pa{\lambda_n^1(1)}\right| \\
									&\cv{n\to+\infty}& 0.
\end{eqnarray*}
With \eqref{controlh}, \eqref{majoAn}, and \eqref{majoBn}, this concludes the proof of Proposition \ref{hinH}.

% ------------------------------------------------------------ %
\subsection{Proof of Theorem \ref{Convergence}}
% ------------------------------------------------------------ %
By Proposition \ref{inegalite_Wasserstein},  for all $n\geq 2$, 
\begin{multline*}
d_2\Big(\loi{\sqrt{n}U_{n},P_n^1\otimes P_n^2\middle | \X_n},\loi{\sqrt{n}U_{n},P^1\otimes P^2}\Big) \\
\leq 
C \inf\limits_{\substack{\pa{Y_{n,a}^*,Y_a}, \pa{Y_{n,b}^*,Y_b}\ i.i.d \tq\\
Y_{n,a}^*,Y_{n,b}^* \thicksim P_n^1\otimes P_n^2,\ 
Y_a, Y_b \thicksim P^1\otimes P^2}}
\Estar{\Big(h\pa{Y_{n,a}^*,Y_{n,b}^*}-h\pa{Y_a,Y_b}\Big)^2}.
\end{multline*}
Our goal is  to construct, for almost all $\omega$ in $\Omega$, a sequence of random variables $\pa{\bar{Y}_{n,\omega,a}^*}_{n\geq 1}$ such that for every $n\geq 1$, $\bar{Y}_{n,\omega,a}^*\sim P_{n,\omega}^1\otimes P_{n,\omega}^2$, where $P_{n,\omega}^j=n^{-1}\sum_{i=1}^n \delta_{X_i^j(\omega)}$ is the $j$th marginal empirical measure corresponding to the realization $\X_n(\omega)$,
a random variable $\bar{Y}_{\omega,a} \sim P^1\otimes P^2$, and $\ac{\pa{\bar{Y}_{n,\omega,b}^*}_{n\geq 1},\bar{Y}_{\omega,b}}$ an independent copy of $\ac{\pa{\bar{Y}_{n,\omega,a}^*}_{n\geq 1},\bar{Y}_{\omega,a}}$ on some probability space $\pa{\Omega_\omega', \mc{A}_\omega',\PP_\omega'}$ depending on $\omega$ such that
\begin{equation}
\label{CVL2}
\mathds{E}'_{\omega}\cro{\Big(h\pa{\bar{Y}_{n,\omega,a}^*,\bar{Y}_{n,\omega,b}^*} - h\pa{\bar{Y}_{\omega,a},\bar{Y}_{\omega,b}}\Big)^2} \cv{n\to+\infty} 0,
\end{equation}
where $\mathds{E}'_{\omega}$ denotes the expectation corresponding to $\PP_\omega'$.
Then from \eqref{CVL2}, we can conclude by noting that, for almost all $\omega$ in $\Omega$, 
\begin{multline*}
\inf\limits_{\substack{\pa{Y_{n,a}^*,Y_a}, \pa{Y_{n,b}^*,Y_b}\ i.i.d \tq\\
Y_{n,a}^*,Y_{n,b}^* \thicksim P_{n,\omega}^1\otimes P_{n,\omega}^2, \ Y_a, Y_b \thicksim P^1\otimes P^2}}
\!\!\!\!\!\Estar{\Big(h\pa{Y_{n,a}^*,Y_{n,b}^*}-h\pa{Y_a,Y_b}\Big)^2}(\omega) \\
\leq \quad \mathds{E}'_{\omega}\cro{\Big(h\pa{\bar{Y}_{n,\omega,a}^*,\bar{Y}_{n,\omega,b}^*} - h\pa{\bar{Y}_{\omega,a},\bar{Y}_{\omega,b}}\Big)^2} \cv{n\to+\infty} 0.
\end{multline*}
%which concludes the proof of the theorem.

To prove \eqref{CVL2}, consider $(\Omega,\mc{A},\PP)$ the probability space on which all the $X_i$'s are defined. In what follows, one can keep in mind that $\Omega$ represents the randomness in the original sequence $(X_i)_i$. 
Thus, a given $\omega$ in $\Omega$ represents a given realization of $(X_i)_i$.\\
As a preliminary step, from Proposition \ref{LGNhcarre}, there exists some subset $\Omega_1$ of $\Omega$ such that $\PP(\Omega_1)=1$ and for every $\omega$ in $\Omega_1$, 
\begin{multline}
\label{LLNh2}
\frac{1}{n^4} \sum_{i,j,k,l=1}^n h^2\pa{\pa{X_i^1(\omega),X_j^2(\omega)},\pa{X_k^1(\omega),X_l^2(\omega)}} \\
\cv{n\to+\infty} \esp{h^2\pa{\pa{X_1^1,X_2^2},\pa{X_3^1,X_4^2}}}.
\end{multline}
Applying Theorem 3 in \cite{Varadarajan58suppl}, since $\pa{\calX,d_{\calX}}$ defined by \eqref{defmetric} is  separable,  $P$-a.s. in $(X_i)_i$, $P_n^1 \cvf{n\to+\infty} P^1$ and $P_n^2 \cvf{n\to+\infty} P^2$. Hence
there exists some subset $ \Omega_2$ of $\Omega$ such that $ \PP(\Omega_2)=1$ and for every $\omega$ in $\Omega_2,$ 
\begin{equation}
\label{CVfaible}
P_{n,\omega}^1\otimes P_{n,\omega}^2 \cvf{n\to+\infty} P^1\otimes P^2,
\end{equation} 
Now, consider $\Omega_0=\Omega_1\cap\Omega_2$, and fix $\omega$ in $\Omega_0$. \\
Following the proof of Skorokhod's representation theorem in \cite[Theorem~11.7.2,  p.~415]{Dudleysuppl}, since $(\calX^2, d_{\calX^2})$ is a separable space, it is possible to construct
\begin{itemize}
\item some probability space $\pa{\Omega_\omega', \mc{A}_\omega',\PP_\omega'}$, 
\item some random variables 
$\bar{Y}_{n,\omega,a}^*:\Omega_\omega'\to \calX^2$, $\bar{Y}_{n,\omega,b}^*:\Omega_\omega'\to \calX^2$ with distribution  $P_{n,\omega}^1\otimes P_{n,\omega}^2$, 
\item 
$\bar{Y}_{\omega,a}:\Omega_\omega'\to \calX^2$, $\bar{Y}_{\omega,b}:\Omega_\omega'\to \calX^2$  with distribution $P^1\otimes P^2,$ 
\end{itemize}
satisfying: 
\begin{itemize}
\item $\PP_\omega'\mbox{-a.s.}$, 
 $\bar{Y}_{n,\omega,a}^* \cv{n\to+\infty} \bar{Y}_{\omega,a}$ and $\bar{Y}_{n,\omega,b}^* \cv{n\to+\infty} \bar{Y}_{\omega,b}$, 
 \item $\left\{\pa{\bar{Y}_{n,\omega,a}^*}_{n\geq 1},  \bar{Y}_{\omega,a}\right\}$ and $\left\{\pa{\bar{Y}_{n,\omega,b}^*}_{n\geq 1},  \bar{Y}_{\omega,b}\right\}$ are independent, 
 \end{itemize} so that w.r.t. the metric $d$ (see \eqref{distance}),
\begin{equation}
\label{cvpsbar}
\PP_\omega'\mbox{-a.s.},\quad \pa{\bar{Y}_{n,\omega,a}^*,\bar{Y}_{n,\omega,b}^*} \cv{n\to+\infty} \pa{\bar{Y}_{\omega,a},\bar{Y}_{\omega,b}}.
\end{equation}
But under $\pa{\mc{A}_{Cont}}$, $h$ is continuous on a set $\mathcal{C}$ s. t.  $\PP_\omega'\pa{\pa{\bar{Y}_{\omega,a},\bar{Y}_{\omega,b}}\in \mathcal{C}}=\pa{P^1\otimes P^2}^{\otimes 2}(\mathcal{C})=1$,  hence
$$\PP_\omega'\mbox{-a.s.},\ h\pa{\bar{Y}_{n,\omega,a}^*,\bar{Y}_{n,\omega,b}^*} \to_{n\to+\infty} h\pa{\bar{Y}_{\omega,a},\bar{Y}_{\omega,b}}.$$
As $\PP_\omega'$-a.s. convergence implies convergence in probability, to obtain \eqref{CVL2},  we only need to prove that the sequence $\pa{ h^2\pa{\bar{Y}_{n,\omega,a}^*,\bar{Y}_{n,\omega,b}^*}}_{n\geq 1}$ is uniformly integrable, according to Theorem 16.6  p. 165 of \cite{Schilling05suppl}. We therefore conclude since \eqref{LLNh2} is equivalent to
\begin{align*}
\mathds{E}'_{\omega}\big[h^2\big(\bar{Y}_{n,\omega,a}^*,&\bar{Y}_{n,\omega,b}^*\big)\big]=\frac{1}{n^4}\sum_{i,j,k,l=1}^n h^2\pa{(X_i^1(\omega),X_j^2(\omega)),(X_k^1(\omega),X_l^2(\omega))} \\
&\cv{n\to+\infty} \esp{h^2\pa{\pa{X_1^1,X_2^2},\pa{X_3^1,X_4^2}}} 
= \mathds{E}'_{\omega}\left[ h^2\pa{\bar{Y}_{\omega,a},\bar{Y}_{\omega,b}}\right].
\end{align*}
\eqref{CVL2} is thus obtained for any $\omega$ in $\Omega_0$, with $\PP(\Omega_0)=1$. This ends the proof.

% ------------------------------------------------------------ %
\subsection{Proof of Proposition \ref{inegalite_Wasserstein}}
% ------------------------------------------------------------ %

Fix some integer $n\geq 2$ and recall that the $P_n^j$'s ($j=1,2$) are the marginal empirical measures associated with $\X_n$.\\ 
Let $\pa{Y_{n,i}^*,Y_i}_{1\leq i\leq n}$ be an i.i.d. sample such that for every $i\!=\!1\ldots n$, $Y_{n,i}^*\sim P_n^1\!\otimes\! P_n^2$, $Y_i\sim P^1\otimes P^2$, and such that, from the definition of Wasserstein's metric $d_2$ recalled in \eqref{Wasserstein}, 
\begin{multline*}d_2^2\left(\calL\left( \left. \sqrt{n} U_{n}, P_n^1\otimes P_n^2 \right|\X_n\right) , \calL\left( \sqrt{n}U_{n},P^1\otimes P^2\right)  \right) \\
\leq  \frac{1}{n(n-1)^2}\ \Estar{\pa{\sum_{i\neq i'} \pa{h\pa{ Y_{n,i}^*, Y_{n,i'}^*} - 
 h\pa{ Y_i, Y_{i'}}}}^2}.
\end{multline*}
Notice that the upper bound is finite under $\pa{\mc{A}_{Mmt}^{*}}$.\\ Introducing  for $(i,i',j,j')$ in $\ac{1,2,\dots,n}^4$, and $m$ in $\{2,3,4\}$, 
$$\text{E}_{(i,i',j,j')}=
\Estar{ \Big(h\big(Y_{n,i}^*, Y_{n,i'}^*\big) - h( Y_i, Y_{i'})\Big)\Big( h\big(Y_{n,j}^*, Y_{n,j'}^*\big) - h(Y_j, Y_{j'})\Big)},$$
$$I_{m}=\ac{(i,i',j,j')\in\ac{1,2,\dots,n}^4 ; i\neq i',\ j\neq j',\ \#\ac{i,i',j,j'}=m},$$ 
where $\#\ac{i,i',j,j'}$ denotes the number of different elements in $\ac{i,i',j,j'}$,  one has:
\begin{multline*}\Estar{\Big(\sum_{i\neq i'} \!\!\pa{ h\pa{ Y_{n,i}^*, Y_{n,i'}^*} - 
h\pa{ Y_i, Y_{i'}}}\Big)^2}\\
=\sum_{(i,i',j,j')\in I_4}\text{E}_{(i,i',j,j')}+\sum_{(i,i',j,j')\in I_3} \text{E}_{(i,i'\!,j,j')} +\sum_{(i,i',j,j')\in I_2} \text{E}_{(i,i',j,j')}.
\end{multline*}
Let us now upper bound each term of this sum separately.\\
If $(i,i',j,j')$ is in $I_4$, then by independence, 
\begin{multline*}
\text{E}_{(i,i',j,j')}
=\big(\Estar{h\pa{ Y_{n,i}^*, Y_{n,i'}^*}} - \esp{h\pa{ Y_i, Y_{i'}}}\big)\times\\
\big(\Estar{h\pa{ Y_{n,j}^*, Y_{n,j'}^*}} - \esp{h\pa{ Y_j, Y_{j'}}}\big).
\end{multline*}
Under $\pa{\mc{A}_{Cent}}$ and $(\mc{A}_{Cent}^*)$, $\esp{h\!\pa{ Y_i, Y_{i'}}}=\Estar{h\!\pa{ Y_{n,i}^*, Y_{n,i'}^*}}=0$, so\\ $\text{E}_{(i,i',j,j')}=0$.\\
If $(i,i',j,j')$ is in $I_3$, by the Cauchy-Schwarz inequality, 
$$
\text{E}_{(i,i',j,j')} \leq \Estar{\Big(h\big(Y_{n,a}^*, Y_{n,b}^*\big) - h(Y_a, Y_b)\Big)^2},$$
where $\pa{Y_{n,a}^*,Y_a}$ and $\pa{Y_{n,b}^*,Y_b}$ are independent copies of the $\pa{Y_{n,i},Y_i}$'s.
If $(i,i',j,j')$ is in $I_2$, then $\text{E}_{(i,i',j,j')} = \Estar{\pa{h\pa{ Y_{n,a}^*,Y_{n,b}^*} - h\pa{ Y_a, Y_b}}^2}$ is immediate.\\
But $\#I_3 = 4 n(n-1)(n-2)$ and $\#I_2 = 2 n(n-1)$, so
\begin{multline*}
d_2^2\left(\calL \left( \left. \sqrt{n} U_{n}, P_n^1\otimes P_n^2 \right|\X_n\right) , \calL\left( \sqrt{n}U_{n},P^1 \otimes P^2\right)  \right)\\
\leq 4\Estar{\Big(h\big(Y_{n,a}^*,Y_{n,b}^*\big) - h\big(Y_a,Y_b\big)\Big)^2 } .\end{multline*}
Since  $\pa{Y_{n,a}^*,Y_a}$ and $\pa{Y_{n,b}^*,Y_b}$ may be arbitrarily chosen, Proposition \ref{inegalite_Wasserstein} follows.

% ------------------------------------------------------------ %
\subsection{Proof of Proposition \ref{LGNhcarre}}
% ------------------------------------------------------------ %

Let us first notice that  \eqref{cvh} is a direct application of the strong law of large numbers for $U$-statistics, proved by Hoeffding \cite{Hoeffding61suppl}. 

Next, for $m$ in $\ac{1,\ldots,4}$, introduce
$$g_{m}\pa{X_{i_1},\dots,X_{i_{m}}} = \sum_{\pa{i,j,k,l}\in I_{\{i_1,\ldots i_m\}}}h^2\pa{\pa{X_i^1,X_j^2},\pa{X_k^1,X_l^2}},$$
where 
$I_{\{i_1,\ldots i_m\}}$ is the set $\ac{\pa{i,j,k,l}\in \ac{i_1,\ldots i_m}^4 ; \# \ac{i,j,k,l}=m}.$\\
Then,
\begin{multline*}
\frac{1}{n^4}\sum_{i,j,k,l=1}^n h^2\!\pa{\pa{X_i^1,X_j^2},\pa{X_k^1,X_l^2}} \\
= \sum_{m=1}^4 \frac{1}{m!}\ \Bigg(\frac{1}{n^4} \sum\limits_{\substack{\pa{i_1,\dots,i_m}\in \{1,\ldots,n\}^m\\ i_1,\dots,i_m \textrm{ all different} }} g_{m}\pa{X_{i_1},\dots,X_{i_{m}}}\!\Bigg).
\end{multline*}
Each of the four terms in the right hand side of the above decomposition being, up to a multiplicative factor, a classical $U$-statistic, and since under $\pa{\mc{A}_{Mmt}^{*}}$, $\esp{|g_{m}(X_{i_1},\dots,X_{i_m})|}<+\infty$, we can now apply the strong law of large numbers for $U$-statistics again. Therefore $P$-a.s. in $(X_i)_i$, 
\begin{equation*}
\frac{1}{n(n\!-\!1)\dots(n\!-\!m\!+\!1)}\!\!\sum_{\pa{i_1,\dots,i_m}}\!\!\! g_{m}\pa{X_{i_1},\dots,X_{i_{m}}}\cv{n\to+\infty} \esp{g_{m}\pa{X_1,\dots,X_{m}}}.
\end{equation*}
In particular, $P$-a.s. in $(X_i)_i$, 
$n^{-4}\sum_{\pa{i_1,\dots,i_m}}g_{m}\pa{X_{i_1},\dots,X_{i_{m}}}$ converges towards $0$ for $m$ in $\ac{1,2,3}$, and
towards $\esp{g_4\pa{X_1,X_2,X_3,X_4}}$ for $m=4$. Finally noticing that $\esp{g_4\pa{X_1,X_2,X_3,X_4}}=4!\esp{h^2\pa{\pa{X_{1}^1,X_{2}^2}\!,\!\pa{X_{3}^1,X_{4}^2}}}$ allows to conclude.

% ------------------------------------------------------------ %
\subsection{Proof of Proposition \ref{TCLUstat}}
% ------------------------------------------------------------ %

Let $(X_i)_i$ be a sequence of i.i.d pairs of point processes with distribution $P^1\otimes P^2$ on $\calX^2$. 
According to $\pa{\mc{A}_{Cent}}$, for $i\neq j$, $\esp{h(X_i,X_j)}=0$. For a better readability, we set $\esp{h\middle | X_i} = \esp{h(X_i,X)\middle| X_i} = \esp{h(X,X_i)\middle| X_i}$ for some $X$ with distribution $P^1\otimes P^2$, and independent of $X_i$.
By Hoeffding's decomposition for non-degenerate $U$-statistics, which also holds when the $X_i$'s are non necessarily real-valued (see \cite{RubinVitale80suppl}) we obtain that
$$\sqrt{n} U_{n}\!\pa{\X_n} =\frac{2}{\sqrt{n}(n-1)} \pa{T_n+M_n},$$
where $T_n=\sum_{i<j}\pa{\esp{h\middle | X_i}+\esp{h\middle | X_j}}$, and $M_n=\sum_{i<j}g(X_i,X_j)$, with $g(X_i,X_j) = h(X_i,X_j) - \esp{h\middle | X_i}- \esp{h\middle | X_j}$. \\
Firstly, we have that $\esp{M_n^2} = \sum_{i<j}\sum_{k<l} \esp{g(X_i,X_j)g(X_k,X_l)}.$ But if $\ac{i,j}\cap \ac{k,l}=\emptyset$, $i<j$, $k<l$, $\esp{g(X_i,X_j)g(X_k,X_l)}=\left(\esp{g(X_i,X_j)}\right)^2=0$. If $\#(\ac{i,j}\cap \ac{k,l})=1$, with for instance $k=i,\ j\neq l$, ($i<j$, $i<l$) (the other cases may be treated similarly), then $$\esp{g(X_i,X_j)g(X_i,X_l)}=\esp{\esp{g(X_i,X_j)\middle | X_i} \esp{  g(X_i,X_l)\middle | X_i}}=0.$$
Therefore,
$\esp{M_n^2} = \sum_{i<j} \esp{g^2(X_i,X_j)} = n(n-1)\esp{g^2(X_1,X_2)}/2,$ and since
$\esp{g^2(X_i,X_j)}<+\infty$, from Chebychev's inequality, we deduce that
\begin{equation}
\label{EqCVM_n}
\frac{2}{\sqrt{n}(n-1)}M_n \cvproba{n\to+\infty}0.
\end{equation}
Secondly, we have that $T_n =(n-1) \sum_{i=1}^n \esp{h\middle|X_i}. $
Since the $\esp{h\middle|X_i}$'s are i.i.d, with $\esp{\esp{h\middle|X_i}}=0$ and $\var{\esp{h\middle|X_i}}= {\sigma_{P^1\otimes P^2}^2}/4,$ thanks to $\pa{\mc{A}_{Mmt}}$, the Central Limit Theorem leads to 
\begin{equation}
\label{EqCVT_n}
\frac{2}{\sqrt{n}(n-1)}T_n \cvloi{n\to+\infty} \mc N\pa{0,\sigma_{P^1\otimes P^2}^2}.
\end{equation}
Thus, combining \eqref{EqCVM_n} and \eqref{EqCVT_n}, Slutsky's lemma ensures the convergence in distribution of $\sqrt{n}U_{n}\!\pa{\X_n}$ towards $\mc N\pa{0,\sigma_{P^1\otimes P^2}^2}$.\\
Now, in order to obtain the convergence in the Wasserstein metric, one needs to check the convergence of the second order moments. Notice that
$$\esp{\pa{\sqrt{n}U_{n}(\X_n)}^2} = \frac{1}{n(n-1)^2} \sum_{i\neq i'}\sum_{j\neq j'} \esp{h(X_i,X_{i'})h(X_j,X_{j'})}.$$
Let us consider all the cases where $i\neq i'$ and $j\neq j'$.\\
If $\#\{i,i',j,j'\} = 4$, $\esp{h(X_i,X_{i'})h(X_j,X_{j'})} =0$,  by independence and $\pa{\mc{A}_{Cent}}$.\\
If $\#\{i,i',j,j'\} \!= \!3$, 
$\esp{h(X_i,X_{i'})h(X_j,X_{j'})}\! = \!\sigma^{2}_{P^1 \otimes P^2}/4$, by symmetry of~$h$.\\
If $\#\{i,i',j,j'\} = 2$,
$\esp{h(X_i,X_{i'})h(X_j,X_{j'})} = \esp{\pa{h(X_1,X_2)}^2}$. Therefore,
$$\esp{\pa{\sqrt{n}U_{n}(\X_n)}^2} =\frac{n-2}{n-1} \sigma^{2}_{P^1 \otimes P^2} + \frac{2}{n-1}\esp{\pa{h(X_1,X_2)}^2}\!\cv{n\to+\infty}\!\!\sigma^{2}_{P^1 \otimes P^2},$$
which ends the proof of Proposition \ref{TCLUstat}.

% ------------------------------------------------------------ %
\subsection{Proof of Corollary \ref{coroBoot}}
% ------------------------------------------------------------ %
By Proposition \ref{TCLUstat}, we have that
\begin{equation}\label{macvenloi}
\calL\left(\sqrt{n}U_{n},P^1\otimes P^2 \right)\cvf{n\to+\infty} \mathcal{N}(0, \sigma^2_{P^1 \otimes P^2}),
\end{equation}
where $\mathcal{N}(0, \sigma^2_{P^1 \otimes P^2})$ has a continuous c.d.f. Therefore, by \cite[Lemma 2.11]{vandervaartsuppl}, 
\begin{equation}\label{cvsupCondGauss0}
\sup_{z\in \R} \left|\proba{\sqrt{n}U_{n}(\X_n^{\independent})\leq z}-\Phi_{0,\sigma^2_{P^1 \otimes P^2}}(z)\right|\cv{n\to+\infty} 0.
\end{equation}
Furthermore, since convergence w.r.t the $d_2$ distance implies weak convergence, Theorem \ref{Convergence} combined with \eqref{macvenloi} leads to
\begin{equation}\label{macvenloistar}
\calL\left(\left.\sqrt{n}U_{n},P_n^1\otimes P_n^2 \right|\X_n\right)\cvf{n\to+\infty} \mathcal{N}(0, \sigma^2_{P^1 \otimes P^2})\ \ \textrm{$P$-a.s. in $(X_i)_i$}.
\end{equation}
Hence, 
\begin{equation}
\label{cvsupCondGauss}
\sup_{z\in \R}\! \left|\proba{\sqrt{n}U_{n}(\X_n^*)\!\leq \!z|\X_n}\!-\!\Phi_{0,\sigma^2_{P^1 \otimes P^2}}\!(z)\right|\!\!\cv{n\to+\infty} \!\!\!0\ \textrm{$P$-a.s. in $(X_i)_i$},
\end{equation}
and the first part of the corollary is obtained. 

\medskip

Moreover, \cite[Lemma 21.2]{vandervaartsuppl} can then be applied to both \eqref{macvenloi} and \eqref{macvenloistar}, to obtain that on the event where \eqref{macvenloistar} holds :
 \begin{equation}\label{cvqstar}
q_{\eta,n}^*\pa{\X_n}\cv{n\to+\infty}\Phi^{-1}_{0,\sigma^2_{P^1 \otimes P^2}}(\eta)\quad \textrm{$P$-a.s. in $(X_i)_i$},
\end{equation}
 and that $q_{\eta,n}^{\independent}$ also converges to $\Phi^{-1}_{0,\sigma^2_{P^1 \otimes P^2}}(\eta)$.

% ------------------------------------------------------------ %
\subsection{Proof of Theorem \ref{thniveauconsistance}}
% ------------------------------------------------------------ %

Let us focus on the sequence of upper-tailed tests in $\Gamma(q^*)$,  the proof for the other tests being similar. 
%Given $\alpha$ in $(0,1)$, we have proved for Corollary \ref{coroBoot} that
%\begin{equation}\label{cvqstar}
%q_{1-\alpha,n}^*\pa{\X_n}\cv{n\to+\infty}\Phi^{-1}_{0,\sigma^2_{P^1 \otimes P^2}}(1-\alpha)\quad \textrm{$P$-a.s. in $(X_i)_i$}.
%\end{equation}

Under $(H_0)$, from Proposition~\ref{TCLUstat} and \eqref{cvqstar}, by Slutsky's lemma,\\ $(\sqrt{n} U_{n}(\X_n), q_{1-\alpha,n}^*\!\pa{\X_n})$ converges in distribution to $(Z,\Phi^{-1}_{0,\sigma^2_{P^1 \otimes P^2}}(1-\alpha))$, where $Z\sim \mathcal{N}(0,\sigma^2_{P^1 \otimes P^2})$. Therefore, under $(H_0)$, 
$$\mathds{P}(\sqrt{n} U_{n}(\X_n)>q_{1-\alpha,n}^*\pa{\X_n})\to_{n\to+\infty}\alpha,$$
 which proves $\pa{\mc{P}_{size}}$.

Under any alternative such that $\int h(x,x')dP(x)dP(x') >0$, by Proposition \ref{LGNhcarre},
$$U_{n}(\X_n)\cv{n\to+\infty} \int h(x,x') dP(x) dP(x') > 0, \ \textrm{$P$-a.s. in $(X_i)_i$}.$$
Furthermore, due to \eqref{cvqstar}, $q_{1-\alpha,n}^*\pa{\X_n}/\sqrt{n}\to_{n\to+\infty} 0$ $\textrm{$P$-a.s. in $(X_i)_i$}$. Hence, $\mathds{P}(\sqrt{n} U_{n}(\X_n)\leq q_{1-\alpha,n}^*(\X_n))\to_{n\to+\infty} 0,$ and thus $\pa{\mc{P}_{consist.}}$ is proved.

% ------------------------------------------------------------ %
\subsection{Proof of Proposition \ref{bootMonteCarlo}}
% ------------------------------------------------------------ %

 As above, we focus on the sequence of upper-tailed tests in $\Gamma(q^*_{MC})$. Let $Z\sim\Nm$ and
define for $z$ in $\mathbb{R},$
$$F^*_{n,\X_n}(z) = \proba{\sqrt{n}U_{n}(\X_n^*)\leq z|\X_n}, \ F^{* B_n}_{n,\X_n}(z) = \frac{1}{B_n}\sum_{b=1}^{B_n} \1{\sqrt{n}U_{n}\!\pa{\X_n^{*b}}\leq z}.$$
By the Dvoretzky-Kiefer-Wolfowitz inequality (see \cite{vandervaartsuppl}), for $n\geq 2$ and $\varepsilon >0$, 
\begin{align*}
\PP\!\bigg(\!\sup_{z\in\R}\Big|F^{* B_n}_{n,\X_n}(z)\!-\! F^*_{n,\X_n}(z)\Big| \!>\! \varepsilon\!\bigg)&\!= \!\esp{\proba{\!\sup_{z\in\R}\left|F^{*,B_n}_{n,\X_n}(z) \!- \!F^*_{n,\X_n}(z)\right| \!>\! \varepsilon \middle| \X_n}} \\
&\leq  2 e^{-2 B_n \varepsilon^2} \cv{n\to+\infty}0,
\end{align*}
that is $\sup_{z\in\R}|F^{* B_n}_{n,\X_n}(z) - F^*_{n,\X_n}(z)| \cvproba{n\to+\infty} 0.$ With \eqref{cvsupCondGauss}, this leads to
\begin{equation}
\label{cvProbaSupEmpVraieboot}
\sup_{z\in\R}\left|F^{*B_n}_{n,\X_n}(z) - \Phi_{0,\sigma^2_{P^1 \otimes P^2}}(z)\right| \cvproba{n\to+\infty} 0.
\end{equation}
We finish the proof using similar arguments as in \cite[Lemma 21.2]{vandervaartsuppl}, combined with a subsequence argument \cite[Theorem 9.2.1]{Dudleysuppl}. % since the convergence in \eqref{cvProbaSupEmpVraieboot} occurs in probability, and not almost surely.
Let $\phi_0$ be an extraction. 
Then, by \eqref{cvProbaSupEmpVraieboot}, there exists an extraction $\phi_1$, and some $\Omega_0\subset\Omega$ such that $\proba{\Omega_0}= 1$, and for every $ \omega$ in $\Omega_0$, 
$$\sup_{z\in\R}\left|F^{*B_{\phi_1\circ\phi_0(n)}}_{{\phi_1\circ\phi_0(n)},\X_{\phi_1\circ\phi_0(n)}}(\omega)(z) - \Phi_{0,\sigma^2_{P^1 \otimes P^2}}(z)\right| \cv{n\to+\infty} 0.$$
From now on, fix $\omega$ in $\Omega_0$. In particular, this fixes a realisation of $\X_n$, and a realisation of $\pa{\X_n^{*1},\ldots,\X_n^{*B_n}}$ and thus, $F^{* B_n}_{n,\X_n}(\omega)$ is deterministic. \\
Hence, $F^{*B_{\phi_1\circ\phi_0(n)}}_{{\phi_1\circ\phi_0(n)},\X_{\phi_1\circ\phi_0(n)}}(\omega)(Z) \!\cvps{n\to+\infty} \!\Phi_{0,\sigma^2_{P^1 \otimes P^2}}(Z),$ and for $\eta$ in $(0,1)$, 
\begin{align*}
&\Phi_{0,1}\pa{\pa{F^{*B_{\phi_1\circ\phi_0(n)}}_{{\phi_1\circ\phi_0(n)},\X_{\phi_1\circ\phi_0(n)}}(\omega)}^{-1}(\eta)} = \proba{F^{*B_{\phi_1\circ\phi_0(n)}}_{{\phi_1\circ\phi_0(n)},\X_{\phi_1\circ\phi_0(n)}}(\omega)(Z) < \eta} \\
&\cv{n\to+\infty} \proba{\Phi_{0,\sigma^2_{P^1 \otimes P^2}}(Z)<\eta}= \Phi_{0,1} \pa{\pa{\Phi_{0,\sigma^2_{P^1 \otimes P^2}}}^{-1}(\eta)}.
\end{align*}
Finally, as $\Phi_{0,1}$ is a one-to-one function and $\Phi_{0,1}^{-1}$ is continuous, 
\begin{multline}
\sqrt{\phi_1\!\circ\!\phi_0(n)}U^{* (\lceil \eta(B_{\phi_1\circ\phi_0(n)}) \rceil)} (\omega) \\
=\pa{ F^{*B_{\phi_1\circ\phi_0(n)}}_{{\phi_1\circ\phi_0(n)},\X_{\phi_1\circ\phi_0(n)}}(\omega)}^{-1}(\eta) 
\cv{n\to+\infty} \Phi^{-1}_{0,\sigma^2_{P^1 \otimes P^2}}(\eta),
\end{multline}
%$$\!\sqrt{\!\phi_1\!\circ\!\phi_0(n)}U^{* (\lceil \eta(B_{\phi_1\circ\phi_0(n)}) \rceil)} \! (\omega)
%\!=\! \pa{ \!F^{*B_{\phi_1\circ\phi_0(n)}}_{{\phi_1\circ\phi_0(n)},\X_{\phi_1\circ\phi_0(n)}}\!\!(\omega)}^{\!-1}\!\!\!\!(\eta) 
% \!\!\cv{n\to+\infty} \!\!\Phi^{-1}_{0,\sigma^2_{P^1 \otimes P^2}}\!\! \!\! \!\!(\eta),$$ 
and this for all $\omega$ in $\Omega_0$, and any initial extraction $\phi_0$. Therefore, we obtain that $\sqrt{n}U^{* (\lceil \eta B_n \rceil)}\cvproba{n\to+\infty}\Phi^{-1}_{0,\sigma^2_{P^1 \otimes P^2}}(\eta).$
 We conclude as for Theorem \ref{thniveauconsistance}.

\subsection{Proof of Theorem \ref{consistencyperm}}
% ------------------------------------------------------------ %
For the sake of clarity and a better readability, we first present a sketch of the proof of this Theorem in Subsection \ref{sketchPerm}. A complete version is detailed in Subsection \ref{complPerm}.

% ------------------------------------------------------------ %
\subsubsection{Sketch of proof of Theorem \ref{consistencyperm}}\label{sketchPerm}
% ------------------------------------------------------------ %
Let $d_{BL}$ denote the bounded Lipschitz metric, which metrizes the weak convergence \cite[Prop. 11.3.2 and Th. 11.3.3]{Dudleysuppl}. For  any variable $Z_n$ depending on $\X_n$ and $\rperm$, $\loi{ Z_n \middle |\X_n}$ denotes  the conditional distribution of $Z_n$ given $\X_n$  and  for any integrable function $f$, $\mathds{E}_{P^1\otimes P^2}[f]=\ds{E}[f(X_1^1,X_2^2)]$.

\smallskip

\noindent $\bullet$ The first step of the proof consists in decomposing $\sqrt{n}U_{n}\pa{\X_n^\rperm}$ in
$$\sqrt{n}U_{n}\pa{\X_n^\rperm}=\frac{n}{n-1} \pa{M_n^\rperm\pa{\X_n}+\frac{R_n^\rperm\pa{\X_n}}{\sqrt{n}}-\frac{T_n\pa{\X_n}}{\sqrt{n}} },$$
where
\begin{itemize}
\item $\displaystyle M_n^\rperm\pa{\X_n}=\frac{1}{\sqrt{n}}\sum_{i\neq j} \1{\rperm(i)=j} C_{i,j},$ 
\item $\displaystyle R_n^\rperm\pa{\X_n}=\sum_{i=1}^n \pa{\1{\rperm(i)=i}-\frac{1}{n}} C_{i,i},$ 
\item $\displaystyle T_n\pa{\X_n}=\frac{1}{n} \sum_{i\neq j}  C_{i,j},$ 
\end{itemize} with 
$$C_{i,j}=\varphi\pa{X_i^1,X_j^2}- \esp{\varphi\pa{X_i^1,X^2}\middle| X_i^1} -\esp{\varphi\pa{X^1,X_j^2}\middle| X_j^2} +  \mathds{E}_{P^1\otimes P^2}\cro{\varphi},$$
$X=(X^1,X^2)$ being $P$-distributed and independent of $(X_i)_i$.\\
We then prove from Cauchy-Schwarz inequality that 
$$\esp{\pa{\esp{\frac{\left|R_n^\rperm\pa{\X_n}\right|}{\sqrt{n}} \middle| \X_n}}^2} \cv{n\to+\infty} 0 \textrm{ and }\esp{\pa{\frac{T_n\pa{\X_n}}{\sqrt{n}} }^2}\cv{n\to+\infty} 0,$$
therefore from Markov's inequality,
$$\esp{\frac{\left|R_n^\rperm\pa{\X_n}\right|}{\sqrt{n}} \middle| \X_n}\cvproba{n\to+\infty} 0 \textrm{ and }  \frac{T_n\pa{\X_n}}{\sqrt{n}} \cvproba{n\to+\infty} 0.$$
From the definition of $d_{BL}$, this allows us to derive that
\begin{equation}\label{dbl}
d_{BL}\pa{ \loi{ \sqrt{n}U_{n}\pa{\X_n^\rperm} \middle | \X_n}, \loi{ \frac{n}{n-1}M_n^\rperm\pa{\X_n}\middle | \X_n }}\cvproba{n\to+\infty} 0.
\end{equation}

\noindent $\bullet$ The second, and most difficult, step of the proof consists in proving that
\begin{equation}\label{secstep}
d_{BL}\pa{\loi{M_n^\rperm\pa{\X_n}\middle| \X_n},\mc{N}\pa{0,\sigma^2_{P^1 \otimes P^2}}} \cvproba{n\to+\infty} 0. 
\end{equation}

\noindent Consider  
\begin{equation}\label{YM}
Y_{n,i}= \frac{1}{\sqrt{n}}\sum_{j=1}^{i-1} \pa{\1{\rperm(i)=j} C_{i,j}+\1{\rperm(j)=i} C_{j,i} },%\quad M_{n,k}=\sum_{i=2}^k Y_{n,i},
\end{equation}

\noindent and for ${\Pi_n'}$ another uniformly distributed random permutation with values in $\Sn{n}$, independent of $\rperm$ and $\X_n$, define accordingly $Y_{n,i}'$ %and $M_{n,k}'$ 
by replacing $\rperm$ by $\Pi_n'$ in \eqref{YM}, 
% $$Y_{n,i}'=\frac{1}{\sqrt{n}}\sum_{j=1}^{i-1} \pa{\1{{\Pi_n'}(i)=j} C_{i,j}+\1{{\Pi_n'}(j)=i} C_{j,i} },\quad M_{n,k}'=\sum_{i=2}^k Y_{n,i}',$$ 
so that $M_n^\rperm\pa{\X_n}=\sum_{i=1}^n Y_{n,i}$ and similarly for $M_n^{{\Pi_n'}}\pa{\X_n}$.\\
Setting $\mc{F}_{n,i}=\sigma\pa{\rperm,{\Pi_n'},X_1,X_2,\ldots,X_i}$ for $n\geq i \geq 2$,  we prove through  technical computations that for $a,b$ in $\R$, $\pa{aY_{n,i}+bY_{n,i}',\mc{F}_{n,i}}_{2\leq i\leq n}$ is a martingale difference array which  satisfies the assumptions of the following result, commonly attributed to Brown \cite{Brownsuppl}.

\begin{thm}
\label{diffmart}
Let $(X_{n,k})_{k\in\{1,...,p_n\},n\in \mathds{N}^* }$ be a martingale difference array, i.e. such that
there exists an array of $\sigma$-algebra $(\mathcal{F}_{n,k})_{k\in\{1,...,p_n\},n\in \mathds{N}^* }$ that is increasing w.r.t. $k$ such that for all $k=1,...,p_n$, $ \esp{X_{n,k}|\mathcal{F}_{n,k-1}}=0. $\\
Let $A_n= \sum_{k=1}^{p_n} \esp{X_{n,k}^2|\mathcal{F}_{n,k-1}}$, and assume that
\begin{itemize}
\item $A_n\cvproba{n\to+\infty} \sigma^2>0$, 
\item $\forall\varepsilon>0$,\ $\displaystyle\sum_{k=1}^{p_n} \esp{X_{n,k}^2 \1{|X_{n,k}|>\varepsilon}}\to_{n\to+\infty} 0.$
\end{itemize}
Then $Z_n=\sum_{k=1}^{p_n} X_{n,k}$  converges in distribution towards $\mathcal{N}(0,\sigma^2)$.
\end{thm}
Thus, given $a,b$ in $\R$, we obtain that  
$$\loi{a M_n^\rperm\pa{\X_n}+b M_n^{{\Pi_n'}} \pa{\X_n}}\cvf{n\to+\infty}  \mc{N}\pa{0, \pa{a^2+b^2}\sigma^{2}_{P^1\otimes P^2}},$$ which, according to the Cramér-Wold device,  leads to Lemma \ref{lm1} below.

\begin{lm}\label{lm1}
Considering the above notation, 
$$\loi{\pa{M_n^\rperm\pa{\X_n},M_n^{{\Pi_n'}} \pa{\X_n}}'}\cvf{n\to+\infty}  \mc{N}_2\pa{0, \pa{\begin{array}{cc}\sigma^{2}_{P^1\otimes P^2} & 0\\
0 & \sigma^{2}_{P^1\otimes P^2} \end{array}}},$$ 
where $\mc{N}_2\pa{M,V}$ denotes the $2$-dimensional Gaussian distribution with mean vector $M$ and variance-covariance matrix $V$.
\end{lm}

From Lemma \ref{lm1},  we deduce that for every $t$ in $\R$,
\begin{equation*}
\left\{ \begin{array}{l}
\proba{M_n^\rperm\pa{\X_n}  \leq t} \cv{n\to+\infty} \Phi_{0,\sigma^2_{P^1\otimes P^2}}(t), \\
\proba{M_n^\rperm\pa{\X_n}  \leq t, M_n^{\Pi_n'}\pa{\X_n}  \leq t} \cv{n\to+\infty} \Phi_{0,\sigma^2_{P^1\otimes P^2}}^2(t).
\end{array}\right.
\end{equation*}
Using Chebychev's inequality,  with the fact (see \cite[Th. 9.2.1]{Dudleysuppl} for instance) that in a separable metric space, convergence in probability is metrizable, and therefore is equivalent to almost sure convergence of a subsequence of any initial subsequence,  we prove that this leads to \eqref{secstep}, and therefore,
$$d_{BL}\pa{ \loi{ \sqrt{n}U_{n}\pa{\X_n^\rperm} \middle | \X_n},   \mc{N}\pa{0,\sigma^2_{P^1 \otimes P^2}} }\cvproba{n\to+\infty} 0.$$

\noindent $\bullet$ The third, and final, step of the proof consists in deriving, by direct computations and the strong law of large numbers of H{\oe}ffding \cite{Hoeffding61suppl}, the convergence of the conditional second order moments
$$\esp{\pa{\sqrt{n} U_{n}\pa{\X_n^\rperm} }^2\middle |\X_n}\cvps{n\to+\infty} \sigma^2_{P^1 \otimes P^2},$$
which ends the proof.

% ------------------------------------------------------------ %
\subsubsection{Complete proof of Theorem \ref{consistencyperm}}\label{complPerm}
% ------------------------------------------------------------ %

Recall that $d_{BL}$ denotes the bounded Lipschitz metric which metrizes the weak convergence, defined by
$$d_{BL}(\mu,\nu) = \sup_{f\in BL,\ \norm{f}_{BL}\leq 1} \abs{\int_\R f\ (d\mu-d\nu)},$$
where, as defined in \cite{Dudleysuppl}, $BL$ is the set of bounded Lipschitz functions on \R, and 
$\displaystyle\norm{f}_{BL} = \norm{f}_\infty + \sup_{x\neq y}\frac{|f(x)-f(y)|}{|x-y|}$.\\
Recall that the proof consists of three steps presented in Section \ref{sketchPerm}. We give below a complete proof for each of these steps.

% -------------------------------------------------- %
\paragraph{First step: decomposition of $\sqrt{n} U_{n}\!\pa{\X_n^\rperm}$ in the {\it Linear case}} It is obvious that by the definition \eqref{hphi} of $h_\varphi$, 
\begin{equation}\label{simpleU}
U_{n}\!\pa{\X_n^\rperm} =\frac{1}{n-1} U_{n}^\rperm,
\end{equation}
 where $\displaystyle U_{n}^\rperm=\sum_{i=1}^n \varphi\!\pa{X_i^1,X_{\rperm(i)}^2}-\frac{1}{n}\sum_{i,j=1}^n  \varphi\!\pa{X_i^1,X_j^2}$, so,
\begin{eqnarray*}
U_{n}^\rperm&=&\sum_{i=1}^n \varphi\!\pa{X_i^1,X_{\rperm(i)}^2}-\frac{1}{n}\sum_{i,j=1}^n \esp{\varphi\!\pa{X_i^1,X_j^2}\middle| X_i^1} \\
&&-\frac{1}{n}\sum_{i,j=1}^n  \esp{\varphi\!\pa{X_i^1,X_j^2}\middle| X_j^2} + \frac{1}{n}\sum_{i,j=1}^n \esp{\varphi\!\pa{X_i^1,X_j^2}}\\
&&-\frac{1}{n}\sum_{i,j=1}^n \big( \varphi\!\pa{X_i^1,X_j^2}- \esp{\varphi\!\pa{X_i^1,X_j^2}\middle| X_i^1} \\
&&-\esp{\varphi\!\pa{X_i^1,X_j^2}\middle| X_j^2} +   \esp{\varphi\!\pa{X_i^1,X_j^2}}\big).
\end{eqnarray*}
On the one hand, if $\mathds{E}_P\cro{f}$ and $\mathds{E}_{P^1\otimes P^2}\cro{f}$ respectively denote $\esp{f\pa{X_1^1,X_1^2}}$, and $\esp{f\pa{X_1^1,X_2^2}}$, for any integrable function $f$, then
\begin{align*}
 &\frac{1}{n}\sum_{i,j=1}^n \esp{\varphi\!\pa{X_i^1,X_j^2}} \\
 &= \sum_{i,j=1}^n \1{\rperm(i)=j}   \esp{\varphi\!\pa{X_i^1,X_j^2}}  \! - \!\sum_{i,j=1}^n \pa{\1{\rperm(i)=j}-\frac{1}{n}}\esp{\varphi\!\pa{X_i^1,X_j^2}}\\
 &=\sum_{i,j=1}^n \1{\rperm(i)=j} \esp{\varphi\!\pa{X_i^1,X_j^2}}\!-\!\pa{\mathds{E}_P\cro{\varphi}\!-\!\mathds{E}_{P^1\otimes P^2}\cro{\varphi}} 
\sum_{i=1}^n \!\!\pa{\1{\rperm(i)=i}-\frac{1}{n}}.
\end{align*}
On the other hand,
\begin{multline*}
 \frac{1}{n}\sum_{i,j=1}^n  \esp{\varphi\!\pa{X_i^1,X_j^2}\middle| X_i^1}=\sum_{i,j=1}^n \1{\rperm(i)=j}  \esp{\varphi\!\pa{X_i^1,X_j^2}\middle| X_i^1}\\
  -\sum_{i=1}^n \pa{\1{\rperm(i)=i}-\frac{1}{n}}  \pa{ \esp{\varphi\!\pa{X_i^1,X_i^2}\middle| X_i^1}-\esp{\varphi\!\pa{X_i^1,X^2}\middle| X_i^1}} ,
 \end{multline*}
 where $X=(X^1,X^2)$ is assumed to be $P$-distributed and independent of $(X_i)_i$, and in the same way,
 \begin{multline*}
 \frac{1}{n}\sum_{i,j=1}^n  \esp{\varphi\!\pa{X_i^1,X_j^2}\middle| X_j^2}=\sum_{i,j=1}^n \1{\rperm(i)=j}  \esp{\varphi\!\pa{X_i^1,X_j^2}\middle| X_j^2}\\
 -\sum_{j=1}^n \pa{\1{\rperm(j)=j}-\frac{1}{n}}  \pa{ \esp{\varphi\!\pa{X_j^1,X_j^2}\middle| X_j^2}-\esp{\varphi\!\pa{X^1,X_j^2}\middle| X_j^2}}.
 \end{multline*}
 Therefore, $U_{n}^\rperm$ is equal to
\begin{align*}
&\sum_{i,j=1}^n \1{\rperm(i)=j} \Big(\varphi\!\pa{X_i^1,X_j^2}-\esp{\varphi\!\pa{X_i^1,X_j^2}\middle| X_i^1}\\
&- \esp{\varphi\!\pa{X_i^1,X_j^2}\middle| X_j^2}+  \esp{\varphi\!\pa{X_i^1,X_j^2}}\Big)\\
&+\sum_{i=1}^n \pa{\1{\rperm(i)=i}-\frac{1}{n}} \Big( \esp{\varphi\!\pa{X_i^1,X_i^2}\middle| X_i^1}+\esp{\varphi\!\pa{X_i^1,X_i^2}\middle| X_i^2}\\
&-\esp{\varphi\!\pa{X_i^1,X^2}\middle| X_i^1}-\esp{\varphi\!\pa{X^1,X_i^2}\middle| X_i^2}-
\mathds{E}_P\cro{\varphi}+\mathds{E}_{P^1\otimes P^2}\cro{\varphi} \Big)\\
&-\frac{1}{n}\sum_{i,j=1}^n \Big( \varphi\!\pa{X_i^1,X_j^2}- \esp{\varphi\!\pa{X_i^1,X_j^2}\middle| X_i^1} -\esp{\varphi\!\pa{X_i^1,X_j^2}\middle| X_j^2} +   \esp{\varphi\!\pa{X_i^1,X_j^2}}\Big).
\end{align*}
As a consequence, setting
$$C_{i,j}=\varphi\!\pa{X_i^1,X_j^2}- \esp{\varphi\!\pa{X_i^1,X^2}\middle| X_i^1} -\esp{\varphi\!\pa{X^1,X_j^2}\middle| X_j^2} +  \mathds{E}_{P^1\otimes P^2}\cro{\varphi},$$ 
\begin{equation}\label{decomp}
\sqrt{n}U_{n}\!\pa{\X_n^\rperm}= \frac{n}{n-1} \pa{M_n^\rperm\!\pa{\X_n}+\frac{R_n^\rperm\!\pa{\X_n}}{\sqrt{n}}-\frac{T_n\!\pa{\X_n}}{\sqrt{n}}},
\end{equation}
with
\begin{eqnarray*}
M_n^\rperm\!\pa{\X_n}&=&\frac{1}{\sqrt{n}}\sum_{i\neq j} \1{\rperm(i)=j} C_{i,j},\\
R_n^\rperm\!\pa{\X_n}&=&\sum_{i=1}^n \pa{\1{\rperm(i)=i}-\frac{1}{n}} C_{i,i},\\
T_n\!\pa{\X_n}&=&\frac{1}{n} \sum_{i\neq j} C_{i,j}.
\end{eqnarray*}
Let us now prove that 
\begin{equation}\label{step1Perm}
d_{BL}\pa{ \loi{ \sqrt{n}U_{n}\!\pa{\X_n^\rperm} \middle | \X_n}, \loi{ \frac{n}{n-1}M_n^\rperm\!\pa{\X_n}\middle | \X_n }}\cvproba{n\to+\infty} 0.
\end{equation}
To do this, first notice that for every function $f$ in $BL$ such that $\norm{f}_{BL} \leq 1$,
\begin{align*}
\bigg| \esp{f\pa{\sqrt{n}U_{n}\!\pa{\X_n^\rperm}} \middle| \X_n} & -\esp{f\pa{ \frac{n}{n-1} M_n^\rperm\!\pa{\X_n}}\middle| \X_n} \bigg|\\
&\leq  \esp{\left| \sqrt{n}U_{n}\!\pa{\X_n^\rperm}-\frac{n}{n-1} M_n^\rperm\!\pa{\X_n}\right| \middle| \X_n}\\
&\leq \frac{n}{n-1} \pa{ \esp{\frac{\left|R_n^\rperm\!\pa{\X_n}\right|}{\sqrt{n}} \middle| \X_n}+ \frac{\left|T_n\!\pa{\X_n}\right|}{\sqrt{n}} }.
\end{align*}
Hence, taking the supremum over \ac{f\in BL ; \norm{f}_{BL} \leq 1},
\begin{multline}\label{majdbl}
d_{BL}\pa{ \loi{ \sqrt{n}U_{n}\!\pa{\X_n^\rperm} \middle | \X_n}, \loi{ \frac{n}{n-1}M_n^\rperm\!\pa{\X_n}\middle | \X_n }}\\
\leq \frac{n}{n-1} \pa{ \esp{\frac{\left|R_n^\rperm\!\pa{\X_n}\right|}{\sqrt{n}} \middle| \X_n}+ \frac{\left|T_n\!\pa{\X_n}\right|}{\sqrt{n}} }.
\end{multline}
Moreover, on the one hand, since $\rperm$ is independent of $(X_i)_i$, by Cauchy-Schwarz inequality,
$$\esp{\pa{\esp{\frac{\left|R_n^\rperm\!\pa{\X_n}\right|}{\sqrt{n}} \middle| \X_n}}^2} \leq \frac{1}{n}\esp{\pa{R_n^\rperm\!\pa{\X_n}}^2},$$
and
\begin{align*}
&\esp{\pa{R_n^\rperm\!\pa{\X_n}}^2}\\
&\leq \sum_{i,j=1}^n \esp{\pa{\1{\rperm(i)=i}-\frac{1}{n}}\pa{\1{\rperm(j)=j}-\frac{1}{n}}} \mathds{E}\left[ C_{i,i}C_{j,j}\right]\\
&\leq C\pa{\esps{P}{\varphi^2}+\esps{P^1 \otimes P^2}{\varphi^2}}\sum_{i,j=1}^n \pa{\esp{\1{\rperm(i)=i}\1{\rperm(j)=j}}-\frac{1}{n^2}}\\
&\leq C\pa{\esps{P}{\varphi^2}+\esps{P^1 \otimes P^2}{\varphi^2}}\pa{\sum_{i=1}^n \pa{\frac{1}{n}-\frac{1}{n^2}}+\sum_{i\neq j} \pa{\frac{1}{n(n-1)}-\frac{1}{n^2}}}\\
&\leq C\pa{\esps{P}{\varphi^2}+\esps{P^1 \otimes P^2}{\varphi^2}}<+\infty.
\end{align*}
Therefore, from Markov's inequality, we deduce that 
 $$\esp{\frac{\left|R_n^\rperm\!\pa{\X_n}\right|}{\sqrt{n}} \middle| \X_n}\cvproba{n\to+\infty} 0.$$
On the other hand, 
$$\esp{\pa{\frac{T_n\!\pa{\X_n}}{\sqrt{n}} }^2}= \frac{1}{n^3}\sum_{i\neq j}\sum_{k\neq l}\esp{C_{i,j} C_{k,l}}.$$
Notice that for $i\neq j$, $\esp{C_{i,j} \middle|X_i}=\esp{C_{i,j} \middle|X_j}=0$.\\
If $\# \ac{i,j,k,l}=4$, then $\esp{C_{i,j} C_{k,l}}=\pa{\esp{C_{i,j}}}^2=0$.\\
If $i,j,l$ are all different,  then 
\begin{eqnarray*}
\esp{C_{i,j}C _{i,l}}&=&\esp{\esp{C_{i,j}C_{i,l} \middle| X_i,X_l}}\\
&=&\esp{\esp{C_{i,j} \middle| X_i}C_{i,l}}\\
&=&0.
\end{eqnarray*}
In the same way, for $i,j,k$ all different, then $\esp{C_{i,j} C_{k,i}}=0$.\\
If $i\neq j$, 
\begin{equation}\label{sig}
\esp{C_{i,j}^2}=\sigma^{2}_{P^1\otimes P^2},\textrm{ and } \esp{C_{i,j}C_{j,i}}\leq \sigma^{2}_{P^1\otimes P^2},
\end{equation}
by the Cauchy-Schwarz inequality. Combining the above computations, we obtain that
$$\esp{\pa{\frac{T_n\!\pa{\X_n}}{\sqrt{n}} }^2}\leq 2\frac{n(n-1)}{n^3} \sigma^{2}_{P^1\otimes P^2}\cv{n\to+\infty} 0,$$
and therefore, 
 $$\frac{T_n\!\pa{\X_n}}{\sqrt{n}} \cvproba{n\to+\infty} 0.$$
Finally, from \eqref{majdbl}, we derive \eqref{step1Perm}.

\paragraph{Second step: asymptotic normality of $M_n^\rperm\!\pa{\X_n}$ given $\X_n$, in probability} Recall that 
\begin{eqnarray*}
M_n^\rperm\!\pa{\X_n}&=&\frac{1}{\sqrt{n}}\sum_{i\neq j} \1{\rperm(i)=j} C_{i,j}\\
&=&\frac{1}{\sqrt{n}}\sum_{i=2}^n\sum_{j=1}^{i-1} \pa{\1{\rperm(i)=j} C_{i,j}+\1{\rperm(j)=i} C_{j,i}}.
\end{eqnarray*}
Let ${\Pi_n'}$ be another uniformly distributed random permutation with values in $\Sn{n}$, independent of $\rperm$ and $\X_n$, 
and 
\begin{eqnarray*}
M_n^{{\Pi_n'}} \pa{\X_n}&=&\frac{1}{\sqrt{n}}\sum_{i\neq j} \1{{\Pi_n'}(i)=j} C_{i,j}\\
&=&\frac{1}{\sqrt{n}}\sum_{i=2}^n\sum_{j=1}^{i-1} \pa{\1{{\Pi_n'}(i)=j} C_{i,j}+\1{{\Pi_n'}(j)=i} C_{j,i} }.
\end{eqnarray*}
Let us now recall the result of Lemma \ref{lm1}:
$$\loi{\pa{M_n^\rperm\!\pa{\X_n},M_n^{{\Pi_n'}} \pa{\X_n}}'}\cvf{n\to+\infty}  \mc{N}_2\pa{0, \pa{\begin{array}{cc}\sigma^{2}_{P^1\otimes P^2} & 0\\
0 & \sigma^{2}_{P^1\otimes P^2} \end{array}}}.$$ 
\paragraph{Proof of Lemma \ref{lm1}} According to the Cramér-Wold device, given $a,b$ in $\R$, we aim at proving that  
$$\loi{a M_n^\rperm\!\pa{\X_n}+b M_n^{{\Pi_n'}} \pa{\X_n}}\cvf{n\to+\infty}  \mc{N}\pa{0, \pa{a^2+b^2}\sigma^{2}_{P^1\otimes P^2}}.$$ 
In order to deal with simpler mathematical expressions, we introduce below some additional notation.
\begin{itemize}
\item For $n\geq i \geq 2$, $\mc{F}_{n,i}=\sigma\pa{\rperm,{\Pi_n'},X_1,X_2,\ldots,X_i}$.
\item Let
 $$Y_{n,i}= \frac{1}{\sqrt{n}}\sum_{j=1}^{i-1} \pa{\1{\rperm(i)=j} C_{i,j}+\1{\rperm(j)=i} C_{j,i} },$$ 
 $$Y_{n,i}'=\frac{1}{\sqrt{n}}\sum_{j=1}^{i-1} \pa{\1{{\Pi_n'}(i)=j} C_{i,j}+\1{{\Pi_n'}(j)=i} C_{j,i} },$$ %$M_{n,k}=\sum_{i=2}^k Y_{n,i}$ and $M_{n,k}'=\sum_{i=2}^k Y_{n,i}'$, 
 so that  $M_n^\rperm\!\pa{\X_n}=\sum_{i=1}^n Y_{n,i}$ and $M_n^{{\Pi_n'}}\pa{\X_n}=\sum_{i=1}^n Y_{n,i}'$.
 \end{itemize}
 Let us first prove that for a fixed integer $n\geq 2$, $\pa{aY_{n,i}+bY_{n,i}',\mc{F}_{n,i}}_{2\leq i\leq n}$ is a martingale difference array. Note that for $2\leq i\leq n$, 
\begin{eqnarray*}
\esp{Y_{n,i} \middle| \mc{F}_{n,i-1}}&=&\frac{1}{\sqrt{n}} \sum_{j=1}^{i-1} \esp{ \1{\rperm(i)=j} C_{i,j}+\1{\rperm(j)=i} C_{j,i}\middle| \mc{F}_{n,i-1}}\\
&=& \frac{1}{\sqrt{n}} \sum_{j=1}^{i-1}\pa{ \1{\rperm(i)=j} \esp{C_{i,j} \middle|X_{j}}+\1{\rperm(j)=i} \esp{C_{j,i} \middle|X_{j}}}\\
&=&0.
\end{eqnarray*}
In the same way, we have that $\esp{Y_{n,i}' \middle| \mc{F}_{n,i-1}}=0$, so
$\esp{aY_{n,i}+bY_{n,i}' \middle| \mc{F}_{n,i-1}}=0.$
From  Theorem \ref{diffmart}, we thus deduce that if
\begin{eqnarray*}
(i)&& \sum_{i=2}^n\esp{\pa{aY_{n,i}+b Y_{n,i}'}^2 \middle| \mc{F}_{n,i-1}} \cvproba{n\to+\infty} (a^2+b^2) \sigma^{2}_{P^1\otimes P^2},\\
(ii)&&\sum_{i=2}^n\esp{\pa{aY_{n,i}+b Y_{n,i}'}^2\1{\left| aY_{n,i}+b Y_{n,i}'\right|>\varepsilon}}\cv{n\to+\infty} 0 \textrm{ for any $\varepsilon>0$},\end{eqnarray*}
then
$$\loi{a M_n^\rperm\!\pa{\X_n}+b M_n^{{\Pi_n'}} \pa{\X_n}}\cvf{n\to+\infty}  \mc{N}\pa{0, \pa{a^2+b^2}\sigma^{2}_{P^1\otimes P^2}}.$$ 
Let us now check that both $(i)$ and $(ii)$ are satisfied.

\smallskip
\noindent{\emph{Assumption $(i)$}}.  In all the following, only consider $n\geq 4$. Noticing that 
\begin{multline}\label{decompab2}
 \sum_{i=2}^n\esp{\pa{aY_{n,i}+b Y_{n,i}'}^2 \middle| \mc{F}_{n,i-1}} \\
 =(a^2+b^2)\sum_{i=2}^n \esp{Y_{n,i}^2\middle| \mc{F}_{n,i-1}} +2ab \sum_{i=2}^n\esp{Y_{n,i}Y_{n,i}'\middle| \mc{F}_{n,i-1}},
\end{multline}
the proof of $(i)$ can be decomposed into two points.\\
The first point consists in proving that
$$\sum_{i=2}^n \esp{Y_{n,i}^2}\cv{n\to+\infty} \sigma^{2}_{P^1\otimes P^2}\quad \textrm{and}\quad
\var{\sum_{i=2}^n \esp{Y_{n,i}^2\middle| \mc{F}_{n,i-1}}}\cv{n\to+\infty}0,$$
which leads, thanks to Chebychev's inequality, to
$$\sum_{i=2}^n \esp{Y_{n,i}^2\middle| \mc{F}_{n,i-1}}\cvproba{n\to+\infty} \sigma^{2}_{P^1\otimes P^2}.$$
The second point consists in proving that
$$\esp{\pa{\sum_{i=2}^n\esp{Y_{n,i}Y_{n,i}'\middle| \mc{F}_{n,i-1}}}^2} \cv{n\to+\infty} 0,$$
so
$$\sum_{i=2}^n\esp{Y_{n,i}Y_{n,i}'\middle| \mc{F}_{n,i-1}}\cvproba{n\to+\infty} 0.$$
$\bullet$ First point. On the one hand,
\begin{multline*}\sum_{i=2}^{n} \esp{Y_{n,i}^2}=\frac{1}{n} \sum_{i=2}^{n} \sum_{j,k=1}^{i-1} \mathds{E}\big[ \pa{\1{\rperm(i)=j} C_{i,j}+\1{\rperm(j)=i} C_{j,i} }\times\\
\pa{\1{\rperm(i)=k} C_{i,k}+\1{\rperm(k)=i} C_{k,i} }\big].\end{multline*}
Furthermore, if $1\leq j\neq k \leq i-1$, 
\begin{eqnarray*}
&&\esp{\pa{\1{\rperm(i)=j} C_{i,j}+\1{\rperm(j)=i} C_{j,i} }\pa{\1{\rperm(i)=k} C_{i,k}+\1{\rperm(k)=i} C_{k,i} }}\\
&&=\esp{\esp{ \pa{\1{\rperm(i)=j} C_{i,j}+\1{\rperm(j)=i} C_{j,i} }\pa{\1{\rperm(i)=k} C_{i,k}+\1{\rperm(k)=i} C_{k,i} }\middle| X_i,X_j,\rperm}}\\
&&=\esp{\pa{\1{\rperm(i)=j} C_{i,j}+\1{\rperm(j)=i} C_{j,i} }\pa{\1{\rperm(i)=k} \esp{C_{i,k}\middle |X_i} +\1{\rperm(k)=i} \esp{C_{k,i} \middle |X_i}}}\\
&&=0.
\end{eqnarray*}
Thus,
\begin{eqnarray*}
&&\sum_{i=2}^{n} \esp{Y_{n,i}^2}=\frac{1}{n} \sum_{i=2}^{n} \sum_{j=1}^{i-1}  \esp{ \pa{\1{\rperm(i)=j} C_{i,j}+\1{\rperm(j)=i} C_{j,i} }^2}\\
&&=\frac{1}{n} \sum_{i=2}^{n} \sum_{j=1}^{i-1}   \esp{\1{\rperm(i)=j} C_{i,j}^2+\1{\rperm(j)=i} C_{j,i}^2 +2\1{\rperm(i)=j} \1{\rperm(j)=i} C_{i,j}C_{j,i}}\\
&&=\frac{1}{n} \sum_{i=2}^{n} \sum_{j=1}^{i-1}\pa{ \frac{2}{n}\esp{C_{i,j}^2} +\frac{2}{n(n-1) }\esp{C_{i,j}C_{j,i}}}\\
&&= \frac{2}{n^2} \sum_{i=2}^{n} (i-1)\pa{ \esp{C_{1,2}^2} +\frac{1}{n-1 }\esp{C_{1,2}C_{2,1}}},
\end{eqnarray*}
so $\sum_{i=2}^{n} \esp{Y_{n,i}^2}=\frac{n-1}{n}\esp{C_{1,2}^2}+\frac{1}{n} \esp{C_{1,2}C_{2,1}}.$
From \eqref{sig}, we derive that
\begin{equation}\label{part1}
\sum_{i=2}^{n} \esp{Y_{n,i}^2}\cv{n\to+\infty} \sigma^{2}_{P^1\otimes P^2}.
\end{equation}
On the other hand, we have that
\begin{eqnarray*}
\esp{Y_{n,i}^2\middle| \mc{F}_{n,i-1}}&=&\frac{1}{n} \sum_{j=1}^{i-1} \1{\rperm(i)=j} \esp{C_{i,j}^2\middle| X_j}+\frac{1}{n} \sum_{j=1}^{i-1} \1{\rperm(j)=i} \esp{C_{j,i}^2\middle| X_j}\\
&&+\frac{2}{n} \sum_{j=1}^{i-1} \1{\rperm(i)=j} \1{\rperm(j)=i}  \esp{C_{i,j}C_{j,i}\middle| X_j}\\
&& +\frac{2}{n} \sum_{1\leq j\neq k \leq i-1} \1{\rperm(i)=j} \1{\rperm(k)=i}  \esp{C_{i,j}C_{k,i}\middle| X_j,X_k}.
\end{eqnarray*}
Then, 
$$\sum_{i=2}^n \pa{\esp{Y_{n,i}^2\middle| \mc{F}_{n,i-1}}-\esp{Y_{n,i}^2}}=A_{n,1}+A_{n,2}+2A_{n,3}+2A_{n,4},$$
with
\begin{eqnarray*}
A_{n,1}&\!\!\!=\!\!\!&\frac{1}{n}\!\sum_{1\leq j<i\leq n}\!\!\! \pa{ \1{\rperm(i)=j} \esp{C_{i,j}^2\middle| X_j} -\frac{1}{n}\esp{C_{i,j}^2}},\\
A_{n,2}&\!\!\!=\!\!\!&\frac{1}{n}\!\sum_{1\leq j<i\leq n}\!\!\! \pa{ \1{\rperm(j)=i} \esp{C_{j,i}^2\middle| X_j} -\frac{1}{n}\esp{C_{i,j}^2}},\\
A_{n,3}&\!\!\!=\!\!\!&\frac{1}{n}\!\sum_{1\leq j<i\leq n}\!\!\! \pa{ \1{\rperm(i)=j}\1{\rperm(j)=i} \esp{C_{i,j}C_{j,i}\middle| X_j} -\frac{1}{n(n-1)}\esp{C_{i,j}C_{j,i}}},\\
A_{n,4}&\!\!=\!\!&\frac{1}{n}\!\sum_{1\leq j\neq k<i\leq n}\!\!\!\!\! \pa{ \1{\rperm(i)=j}\1{\rperm(k)=i} \esp{C_{i,j}C_{k,i}\middle| X_j,X_k}}.
\end{eqnarray*}
Thus, 
\begin{multline}\label{controlA}
\var{\sum_{i=2}^n \pa{\esp{Y_{n,i}^2\middle| \mc{F}_{n,i-1}}}}\\\leq 4\pa{\esp{A_{n,1}^2}+\esp{A_{n,2}^2}+4\esp{A_{n,3}^2}+4\esp{A_{n,4}^2}}.
\end{multline}
Let us now control each term of the above right-hand side.\\
\emph{Convergence of $\esp{A_{n,1}^2}$ and $\esp{A_{n,2}^2}$}.
\begin{eqnarray*}
\esp{A_{n,1}^2}&=&\frac{1}{n^2} \sum_{1\leq j<i\leq n} \sum_{1\leq l<k\leq n} \Big(\esp{ \1{\rperm(i)=j}\1{\rperm(k)=l}}\times\\
&&\esp{\esp{C_{i,j}^2\middle| X_j}\esp{C_{k,l}^2\middle| X_l} }-\frac{1}{n^2}\pa{\esp{C_{k,l}^2}}^2\Big).
\end{eqnarray*}
Let us now consider all the cases where $1\leq j<i\leq n$, and $1\leq l<k\leq n$.\\
If $i=k$ and $j=l$, then 
$$\esp{ \1{\rperm(i)=j}\1{\rperm(k)=l}}\esp{\esp{C_{i,j}^2\middle| X_j}\esp{C_{k,l}^2\middle| X_l}}=\frac{1}{n}\esp{\pa{\esp{C_{2,1}^2\middle| X_1}}^2}.$$
If $i=k$ and $j\neq l$, or if $i\neq k$ and $j=l$, then
$$\esp{ \1{\rperm(i)=j}\1{\rperm(k)=l}}\esp{\esp{C_{i,j}^2\middle| X_j}\esp{C_{k,l}^2\middle| X_l}}=0.$$
If $i\neq k$ and $j\neq l$, then 
$$\esp{ \1{\rperm(i)=j}\1{\rperm(k)=l}}\esp{\esp{C_{i,j}^2\middle| X_j}\esp{C_{k,l}^2\middle| X_l}}=\frac{1}{n(n-1)}\pa{\esp{C_{2,1}^2}}^2.$$
By combining these results, from \eqref{sig} and under the assumption $\pa{\mc{A}_{\varphi,Mmt}}$, we obtain that
\begin{eqnarray*}
\esp{A_{n,1}^2}&\leq& \frac{n-1}{2n^2}\pa{\esp{\pa{\esp{C_{2,1}^2\middle| X_1}}^2}-\frac{\sigma^4_{P^1\otimes P^2}}{n}}\\
&&+\ C\ n^2 \pa{\frac{1}{n(n-1)}-\frac{1}{n^2}}\sigma^4_{P^1\otimes P^2}\cv{n\to+\infty} 0.
\end{eqnarray*}
One can prove in the same way that $\esp{A_{n,2}^2}\cv{n\to+\infty} 0.$\\
\emph{Convergence of $\esp{A_{n,3}^2}$.} We easily prove that
$$\esp{A_{n,3}^2}=\frac{1}{n^2}\sum_{1\leq j<i\leq n} \sum_{1\leq l<k\leq n}\kappa_{i,j,k,l}-\frac{1}{4n^2} \pa{\esp{C_{1,2}C_{2,1}}}^2,$$
where 
\begin{multline*}
\kappa_{i,j,k,l}=\esp{\1{\rperm(i)=j}\1{\rperm(j)=i}\1{\rperm(k)=l}\1{\rperm(l)=k}}\times\\
\esp{\esp{C_{i,j}C_{j,i}\middle |X_j} \esp{C_{k,l}C_{l,k}\middle |X_l} }.
\end{multline*}
Let us again consider $\kappa_{i,j,k,l}$ in all the cases where $1\leq j<i\leq n$, and $1\leq l<k\leq n$.
If $i=k$ and $j=l$, then 
$$\kappa_{i,j,k,l}=\frac{1}{n(n-1)} \esp{\pa{\esp{C_{2,1}C_{1,2}\middle| X_1}}^2}.$$
If $i=k$ and $j\neq l$, or if $i\neq k$ and $j=l$, then $\kappa_{i,j,k,l}=0.$\\
If $i\neq k$ and $j\neq l$, then 
$$\kappa_{i,j,k,l}=\frac{\pa{\esp{C_{1,2}C_{2,1}}}^2}{n(n-1)(n-2)(n-3)}.$$
Thus, under $\pa{\mc{A}_{\varphi,Mmt}}$, we finally have that
$$\esp{A_{n,3}^2}\leq \frac{1}{2n^2}\esp{\pa{\esp{C_{1,2}C_{2,1}\middle| X_1}}^2}+C\frac{n\pa{\esp{C_{1,2}C_{2,1}}}^2}{(n-1)(n-2)(n-3)}\cv{n\to+\infty} 0.$$
\emph{Convergence of $\esp{A_{n,4}^2}$.}
\begin{multline*}
\esp{A_{n,4}^2}=
\frac{1}{n^2}\!\!\!\sum\limits_{\substack{1\leq j\neq k<i\leq n\\ 1\leq p\neq q<l \leq n}}\!\!\!\!\!\!\Big(\esp{\1{\rperm(i)=j}\1{\rperm(k)=i}\1{\rperm(l)=p}\1{\rperm(q)=l}}\times\\
 \esp{\esp{C_{i,j}C_{k,i}\middle| X_j,X_k} \esp{C_{l,p}C_{q,l}\middle| X_p,X_q}}\Big).
 \end{multline*}
Let us  consider all the cases where $1\leq j\neq k<i\leq n$, and $1\leq p\neq q<l \leq n$.\\
If $\#\{j,k,p,q\}\geq 3$, there exists at least one element in $\{j,k,p,q\}$, $j$ for instance (the other cases are studied in the same way), which differs from the other ones. Then,
\begin{align*}
\ds{E}\big[\ds{E}\big[C_{i,j}C_{k,i}&\big| X_j,X_k\big] \esp{C_{l,p}C_{q,l}\middle| X_p,X_q}\big]\\
&=\esp{\esp{\esp{C_{i,j}C_{k,i}\middle| X_j,X_k} \esp{C_{l,p}C_{q,l}\middle| X_p,X_q}\middle| X_k,X_p,X_q}}\\
&=\esp{\esp{C_{i,j}C_{k,i}\middle| X_k} \esp{C_{l,p}C_{q,l}\middle| X_p,X_q}}\\
&=\esp{\esp{\esp{C_{i,j}C_{k,i}\middle|X_i,X_k} \middle| X_k}\esp{C_{l,p}C_{q,l}\middle| X_p,X_q}}\\
&=\esp{\esp{C_{k,i}\esp{C_{i,j}\middle|X_i} \middle| X_k}\esp{C_{l,p}C_{q,l}\middle| X_p,X_q}}.
\end{align*}
Since $\esp{C_{i,j}\middle|X_i}=0$, this leads to
\begin{equation}\label{eqstar}
\esp{\esp{C_{i,j}C_{k,i}\middle| X_j,X_k} \esp{C_{l,p}C_{q,l}\middle| X_p,X_q}}=0.
\end{equation}
If $j=p$, $k=q$, and $i=l$, then,
$$\esp{\1{\rperm(i)=j}\1{\rperm(k)=i}\1{\rperm(l)=p}\1{\rperm(q)=l}}=\frac{1}{n(n-1)},$$
and 
\begin{eqnarray*}
\left|\esp{\esp{C_{i,j}C_{k,i}\middle| X_j,X_k} \esp{C_{l,p}C_{q,l}\middle| X_p,X_q}}\right|&=&\esp{\pa{\esp{C_{i,j}C_{k,i}\middle| X_j,X_k}}^2}\\
&=&\esp{\pa{\esp{C_{3,1}C_{2,3}\middle| X_1,X_2}}^2}\\
&<&+\infty \textrm{ \quad under $\pa{\mc{A}_{\varphi,Mmt}}$}.
\end{eqnarray*}
If $j=p$, $k=q$, and $i\neq l$, then $\1{\rperm(k)=i}\1{\rperm(q)=l}=0$, so
$$\esp{\1{\rperm(i)=j}\1{\rperm(k)=i}\1{\rperm(l)=p}\1{\rperm(q)=l}}=0.$$
If $j=q$, $k=p$, and $i=l$, then $\1{\rperm(i)=j}\1{\rperm(l)=p}=0$, so
$$\esp{\1{\rperm(i)=j}\1{\rperm(k)=i}\1{\rperm(l)=p}\1{\rperm(q)=l}}=0.$$
If $j=q$, $k=p$, and $i\neq l$, then 
$$\esp{\1{\rperm(i)=j}\1{\rperm(k)=i}\1{\rperm(l)=p}\1{\rperm(q)=l}}=\frac{(n-4)!}{n !},$$
and 
\begin{align*}
\left|\esp{\esp{C_{i,j}C_{k,i}\middle| X_j,X_k} \esp{C_{l,p}C_{q,l}\middle| X_p,X_q}}\right|&=\left|\esp{\esp{C_{i,j}C_{k,i}C_{l,k}C_{j,l}\middle| X_j,X_k}}\right|\\
&\leq  \esp{|C_{3,1}C_{2,3}C_{4,2}C_{1,4}|}\\
&<+\infty \textrm{ \quad under $\pa{\mc{A}_{\varphi,Mmt}}$}.
\end{align*}
By combining these results, we obtain that
\begin{multline*}
\esp{A_{n,4}^2}\leq C\frac{n^3}{n^2} \frac{\esp{\pa{\esp{C_{3,1}C_{2,3}\middle| X_1,X_2}}^2}}{n(n-1)}\\
+C'\ \frac{n^4}{n^2} \frac{(n-4)!}{n !}\esp{C_{3,1}C_{2,3}C_{4,2}C_{1,4}}
 \cv{n\to+\infty} 0.\end{multline*}
From \eqref{controlA}, and the above results of convergence towards $0$ for $\esp{A_{n,1}^2}$, $\esp{A_{n,2}^2}$, $\esp{A_{n,3}^2}$, and $\esp{A_{n,4}^2}$, we firstly derive that
$$\var{\sum_{i=2}^n \pa{\esp{Y_{n,i}^2\middle| \mc{F}_{n,i-1}}}}\cv{n\to+\infty} 0.$$
$\bullet$ Second point. Notice that
$$\esp{Y_{n,i}Y_{n,i}'\middle| \mc{F}_{n,i-1}}=B_{n,1}+B_{n,2}+B_{n,3}+B_{n,4},$$
with
\begin{align*}
&B_{n,1}=\frac{1}{n} \!\sum_{1\leq j<i\leq n}\!\!\1{\rperm(i)=j} \1{{\Pi_n'}(i)=j} \esp{ C_{i,j}^2\middle| X_j} ,\\
&B_{n,2}=\frac{1}{n} \!\sum_{1\leq j<i\leq n}\!\!\1{\rperm(j)=i} \1{{\Pi_n'}(j)=i} \esp{C_{j,i}^2\middle|X_j},\\
&B_{n,3}=\frac{1}{n} \!\sum_{1\leq j<i\leq n}\!\!\!\!\pa{ \1{\rperm(i)=j} \1{{\Pi_n'}(j)=i} +\1{\rperm(j)=i} \1{{\Pi_n'}(i)=j}}\esp{ C_{i,j} C_{j,i}\middle|X_j},
\end{align*}
and
\begin{align*}B_{n,4}=\frac{1}{n} \sum_{1\leq j\neq k<i\leq n}  \Big(&\1{\rperm(i)=j} \1{{\Pi_n'}(i)=k} \esp{ C_{i,j}C_{i,k}\middle| X_j,X_k}\\
&+\1{\rperm(i)=j} \1{{\Pi_n'}(k)=i} \esp{ C_{i,j} C_{k,i}\middle|X_j,X_k}\\
&+\1{\rperm(j)=i} \1{{\Pi_n'}(i)=k} \esp{ C_{j,i}C_{i,k}\middle| X_j,X_k}\\
&+\1{\rperm(j)=i} \1{{\Pi_n'}(k)=i} \esp{C_{j,i}C_{k,i} \middle| X_j,X_k}\Big).
\end{align*}
Thus,
\begin{multline}\label{controlB}
\esp{\pa{\sum_{i=2}^n\esp{Y_{n,i}Y_{n,i}'\middle| \mc{F}_{n,i-1}}}^2} \\
\leq 4 \pa{\esp{B_{n,1}^2}+\esp{B_{n,2}^2}+\esp{B_{n,3}^2}+\esp{B_{n,4}^2}}.
\end{multline}
\emph{Convergence of $\esp{B_{n,1}^2}$ and $\esp{B_{n,2}^2}$.} It can be proved that
\begin{multline*}
\esp{B_{n,1}^2}
\leq \frac{1}{n^3} \sum_{1\leq j<i\leq n}\sum_{1\leq l<k\leq n} \esp{\1{\rperm(i)=j} \1{\rperm(k)=l}}\times\\
\esp{\esp{ C_{i,j}^2\middle| X_j} \esp{ C_{k,l}^2\middle| X_l} }.
\end{multline*}
Then, with the same computations as for the convergence of $ \esp{A_{n,1}^2}$ above, we prove that
$$\esp{B_{n,1}^2}\leq \frac{n-1}{2n^3}\esp{\pa{\esp{C_{1,2}^2\middle| X_2}}^2}+C\frac{\sigma^4_{P^1\otimes P^2}}{n-1}\cv{n\to+\infty}0.$$
In the same way, we also prove that $\esp{B_{n,2}^2}\cv{n\to+\infty}0$.\\
\emph{Convergence of $\esp{B_{n,3}^2}$.} We also have that
\begin{multline*}
\esp{B_{n,3}^2}\leq \frac{4}{n^2} \sum_{1\leq j<i\leq n}\sum_{1\leq l<k\leq n} \esp{\1{\rperm(i)=j} \1{\rperm(k)=l}}\times\\
\esp{\1{{\Pi_n'}(j)=i}\1{{\Pi_n'}(l)=k}} \esp{\esp{ C_{i,j} C_{j,i}\middle|X_j}\esp{ C_{k,l} C_{l,k}\middle|X_l}}.
\end{multline*}
Now, with similar computations as for the convergence of $ \esp{A_{n,1}^2}$ above again, we prove that
$$\esp{B_{n,3}^2}\leq 2\frac{n-1}{n^3}\esp{\pa{\esp{C_{1,2}C_{2,1}\middle| X_2}}^2}+C\frac{ \pa{\esp{C_{1,2}C_{2,1}}}^2}{n-1}\cv{n\to+\infty}0.$$
\emph{Convergence of $\esp{B_{n,4}^2}$.} Setting 
$$B_{n,4,1}=\frac{1}{n} \sum_{1\leq j\neq k<i\leq n} \1{\rperm(i)=j} \1{{\Pi_n'}(i)=k} \esp{ C_{i,j}C_{i,k}\middle| X_j,X_k},$$
$$B_{n,4,2}=\frac{1}{n} \sum_{1\leq j\neq k<i\leq n} \1{\rperm(i)=j} \1{{\Pi_n'}(k)=i} \esp{ C_{i,j} C_{k,i}\middle| X_j,X_k},$$
$$B_{n,4,3}=\frac{1}{n} \sum_{1\leq j\neq k<i\leq n} \1{\rperm(j)=i} \1{{\Pi_n'}(i)=k} \esp{ C_{j,i}C_{i,k}\middle| X_j,X_k},$$
$$B_{n,4,4}=\frac{1}{n} \sum_{1\leq j\neq k<i\leq n} \1{\rperm(j)=i} \1{{\Pi_n'}(k)=i} \esp{C_{j,i}C_{k,i} \middle| X_j,X_k},$$
then $B_{n,4}=B_{n,4,1}+B_{n,4,2}+B_{n,4,3}+B_{n,4,4}$ and in particular, 
$$\esp{{B_{n,4}}^2} \leq 4 \pa{\esp{{B_{n,4}}^2}+\esp{{B_{n,4}}^2}+\esp{{B_{n,4}}^2}+\esp{{B_{n,4}}^2}}.$$
Yet,
\begin{multline*}
\esp{B_{n,4,1}^2}=\frac{1}{n^2} \sum_{1\leq j\neq k<i\leq n} \sum_{1\leq p\neq q<l\leq n} \esp{\1{\rperm(i)=j} \1{\rperm(l)=p}} \times\\
\esp{ \1{{\Pi_n'}(i)=k} \1{{\Pi_n'}(l)=q} }\esp{\esp{ C_{i,j}C_{i,k}\middle| X_j,X_k}\esp{C_{l,p}C_{l,q}\middle| X_p,X_q}}.
\end{multline*}
Now, consider all the cases where $1\leq j\neq k<i\leq n$, $1\leq p\neq q<l\leq n$.\\
If $\#\{j,k,p,q\}\geq 3$, using a similar argument as in \eqref{eqstar}, we obtain that 
$$\esp{\esp{ C_{i,j}C_{i,k}\middle| X_j,X_k}\esp{C_{l,p}C_{l,q}\middle| X_p,X_q}}=0.$$
If $j=p$, $k=q$, and $i=l$, then,
$$\esp{\1{\rperm(i)=j} \1{\rperm(l)=p}}\esp{ \1{{\Pi_n'}(i)=k} \1{{\Pi_n'}(l)=q} } =\frac{1}{n^2},$$
and
\begin{eqnarray*}
\left|\esp{\esp{C_{i,j}C_{i,k}\middle| X_j,X_k} \esp{C_{l,p}C_{l,q}\middle| X_p,X_q}}\right|&=&\esp{\pa{\esp{C_{i,j}C_{i,k}\middle| X_j,X_k}}^2}\\
&=&\esp{\pa{\esp{C_{3,1}C_{3,2}\middle| X_1,X_2}}^2}\\
&<&+\infty \textrm{ \quad under $\pa{\mc{A}_{\varphi,Mmt}}$}.
\end{eqnarray*}
If $j=p$, $k=q$, $i\neq l$, or if $j=q$, $k=p$,  $i=l$, then $\1{\rperm(i)=j}\1{\rperm(l)=p}$ is equal to $0$, so
$$\esp{\1{\rperm(i)=j} \1{\rperm(l)=p}}\esp{ \1{{\Pi_n'}(i)=k} \1{{\Pi_n'}(l)=q} } =0.$$
If $j=q$, $k=p$, and $i\neq l$, then 
$$\esp{\1{\rperm(i)=j} \1{\rperm(l)=p}}\esp{ \1{{\Pi_n'}(i)=k} \1{{\Pi_n'}(l)=q} } =\frac{1}{n^2(n-1)^2},$$
and 
\begin{align*}
\left|\esp{\esp{C_{i,j}C_{i,k}\middle| X_j,X_k} \esp{C_{l,p}C_{l,q}\middle| X_p,X_q}}\right|&=\left|\esp{\esp{C_{i,j}C_{i,k}C_{l,k}C_{l,j}\middle| X_j,X_k}}\right|\\
&= \esp{\left|C_{3,1}C_{2,3}C_{4,2}C_{1,4}\right|}\\
&<+\infty \textrm{ \quad under $\pa{\mc{A}_{\varphi,Mmt}}$}.
\end{align*}
By combining these results, we obtain that
$$\esp{B_{n,4,1}^2}\leq C\frac{\esp{\pa{\esp{C_{3,1}C_{3,2}\middle| X_1,X_2}}^2}   }{n} +
C'\frac{\esp{C_{3,1}C_{2,3}C_{4,2}C_{1,4}}}{(n-1)^2}\cv{n\to+\infty} 0.$$
Following the same lines of proof, we furthermore obtain that $\esp{B_{n,4,2}^2}$, $\esp{B_{n,4,3}^2}$, and $\esp{B_{n,4,4}^2}$ also converge towards $0$. Hence, $\esp{B_{n,4}^2}\cv{n\to+\infty}0.$
From \eqref{controlB}, and the convergence towards $0$ of $\esp{B_{n,1}^2}$, $\esp{B_{n,2}^2}$, $\esp{B_{n,3}^2}$, and $\esp{B_{n,4}^2}$, we derive that
$$\esp{\pa{\sum_{i=2}^n\esp{Y_{n,i}Y_{n,i}'\middle| \mc{F}_{n,i-1}}}^2}\cv{n\to+\infty}0,$$
which finally allows to conclude that assumption $(i)$ is satisfied.

\smallskip

\noindent \emph{Assumption $(ii)$.} Given $\varepsilon>0$, let us prove that 
$$\sum_{i=2}^n\esp{\pa{aY_{n,i}+b Y_{n,i}'}^2\1{\left| aY_{n,i}+b Y_{n,i}'\right|>\varepsilon}}\cv{n\to+\infty} 0.$$
\begin{align*}
\sum_{i=2}^n\esp{\pa{aY_{n,i}+b Y_{n,i}'}^2\1{\left| aY_{n,i}+b Y_{n,i}'\right|>\varepsilon}}&\leq \frac{1}{\varepsilon^2}
\sum_{i=2}^n\esp{\pa{aY_{n,i}+b Y_{n,i}'}^4}\\
&\leq \frac{2^3}{\varepsilon^2} \sum_{i=2}^n \pa{a^4 \esp{Y_{n,i}^4}+b^4 \esp{{Y_{n,i}'}^4}}\\
&\leq \frac{2^3(a^4+b^4)}{\varepsilon^2} \sum_{i=2}^n \esp{Y_{n,i}^4}.
\end{align*}
Since $Y_{n,i}=n^{-1/2} \pa{\1{\rperm(i)<i}C_{i,\rperm(i)}+ \1{\Pi_n^{-1}(i)<i}C_{\Pi_n^{-1}(i),i}},$
\begin{eqnarray*}
\esp{Y_{n,i}^4}&\leq&\frac{2^3}{n^2} \esp{\1{\rperm(i)<i}C_{i,\rperm(i)}^4+ \1{\Pi_n^{-1}(i)<i}C_{\Pi_n^{-1}(i),i}^4 }\\
&\leq&\frac{2^3}{n^2} \sum_{j=1}^{i-1} \pa{\esp{\1{\rperm(i)=j}C_{i,j}^4}+\esp{\1{\Pi_n^{-1}(i)=j}C_{j,i}^4}}\\
&\leq&\frac{2^4}{n^2} \esp{C_{1,2}^4}.
\end{eqnarray*}
We thus obtain that
$$\sum_{i=2}^n\esp{\pa{aY_{n,i}+b Y_{n,i}'}^2\1{\left| aY_{n,i}+b Y_{n,i}'\right|>\varepsilon}}\leq \frac{2^7(a^4+b^4)}{\varepsilon^2 n} \esp{C_{1,2}^4},$$
where the right-hand side tends to $0$ as soon as $\esp{C_{1,2}^4}<+\infty$.\\
This last condition is ensured by $\pa{\mc{A}_{\varphi,Mmt}}$, which allows to confirm that assumption $(ii)$ is also checked, and that
$$\loi{a M_n^\rperm\!\pa{\X_n}+b M_n^{{\Pi_n'}} \pa{\X_n}}\cvf{n\to+\infty}  \mc{N}\pa{0, \pa{a^2+b^2}\sigma^{2}_{P^1\otimes P^2}}.$$ 
This ends the proof of Lemma \ref{lm1}.

\medskip

Recall that we aim at proving that 
 $$d_{BL}\pa{\loi{M_n^\rperm\!\pa{\X_n}\middle| \X_n},\mc{N}\pa{0,\sigma^{2}_{P^1 \otimes P^2}}} \cvproba{n\to+\infty} 0.$$
 From Lemma \ref{lm1},  we deduce that for every $t$ in $\R$,
\begin{equation*}
\left\{ \begin{array}{l}
\proba{M_n^\rperm\!\pa{\X_n}  \leq t} \cv{n\to+\infty} \Phi_{0,\sigma^{2}_{P^1\otimes P^2}}(t), \\
\proba{M_n^\rperm\!\pa{\X_n}  \leq t, M_n^{\Pi_n'}\pa{\X_n}  \leq t} \cv{n\to+\infty} \Phi_{0,\sigma^{2}_{P^1\otimes P^2}}^2(t).
\end{array}\right.
\end{equation*}
Setting $M_n=M_n^\rperm\!\pa{\X_n}$ for the sake of simplicity, this leads to
\begin{equation}\label{cvfdr}
\left\{ \begin{array}{l}
\esp{\esp{\1{M_n \leq t} \middle |\X_n}}\cv{n\to+\infty} \Phi_{0,\sigma^{2}_{P^1\otimes P^2}}(t), \\
\esp{\pa{\esp{\1{M_n \leq t}\middle| \X_n}}^2}\cv{n\to+\infty} \Phi_{0,\sigma^{2}_{P^1\otimes P^2}}^2(t).
\end{array}\right.
\end{equation}
In a separable metric space, convergence in probability is metrizable (see \cite[Th. 9.2.1]{Dudleysuppl} for instance),  therefore it is equivalent to almost sure convergence of a subsequence of any initial subsequence. So, let us fix an initial extraction $\phi_0:\N\to\N$, which defines  a subsequence $\pa{M_{\phi_0(n)}}_{n\in\N}$ of $\pa{M_n}_{n\in\N}$. Let us denote by $(q_m)_{m\in\N}$ a sequence such that $\{q_m,m\in \mathds{N}\}=\mathds{Q}$. For any $m$ in $\N$, from \eqref{cvfdr}, we derive that 
\begin{equation*}
\left\{ \begin{array}{l}
\esp{\esp{\1{M_{\phi_0(n)} \leq q_m} \middle |\X_{\phi_0(n)}}}\cv{n\to+\infty} \Phi_{0,\sigma^{2}_{P^1\otimes P^2}}(q_m), \\
\esp{\pa{\esp{\1{M_{\phi_0(n)} \leq q_m}\middle| \X_{\phi_0(n)}}^2}}\cv{n\to+\infty} \Phi_{0,\sigma^{2}_{P^1\otimes P^2}}^2(q_m),
\end{array}\right.
\end{equation*}
which leads (by Chebychev's inequality) to 
\begin{equation}\label{eq1cvprob}
\esp{\1{M_{\phi_0(n)} \leq q_m} \middle |\X_{\phi_0(n)}}\cvproba{n\to+\infty}  \Phi_{0,\sigma^{2}_{P^1\otimes P^2}}(q_m).
\end{equation}
Therefore, there exist an extraction $\phi_1$ and a subset $\Omega_1$ of $\Omega$ such that $\proba{\Omega_1}=1$, and for every $\omega$ in $\Omega_1$, 
$$\esp{\1{M_{\phi_1\circ\phi(n)} \leq q_1} \middle |\X_{\phi_1\circ\phi(n)}}(\omega)\cv{n\to+\infty}  \Phi_{0,\sigma^{2}_{P^1\otimes P^2}}(q_1).$$
Now, let $m\geq 1$ for which there exist an extraction $\phi_m$ and a subset $\Omega_m$ of $\Omega$ such that $\proba{\Omega_m}=1$, and for every $\omega\in\Omega_m$, 
$$\esp{\1{M_{\phi_m\circ \phi_{m-1}\circ \ldots\circ\phi_0(n)} \leq q_m} \middle |\X_{\phi_m\circ \phi_{m-1}\circ \ldots\circ\phi_0(n)}}(\omega)\cv{n\to+\infty}  \Phi_{0,\sigma^{2}_{P^1\otimes P^2}}(q_m).$$
Then, from \eqref{eq1cvprob}, there also exist an extraction $\phi_{m+1}$ and a subset $\Omega_{m+1}$ of $\Omega$ such that $\proba{\Omega_{m+1}}=1$, and for every $\omega$ in $\Omega_{m+1}$, 
\begin{multline*}
\esp{\1{M_{\phi_{m+1}\circ\phi_m\circ \phi_{m-1}\circ \ldots\circ\phi_0(n)} \leq q_{m+1}} \middle |\X_{\phi_{m+1}\circ\phi_m\circ  \ldots\circ\phi_0(n)}}(\omega)\\
\cv{n\to+\infty}  \Phi_{0,\sigma^{2}_{P^1\otimes P^2}}(q_{m+1}).\end{multline*}
Setting $\tilde{\Omega}= \bigcap_{m\in\N} \Omega_m$ and for every $n$ in $\N$, $\tilde\phi(n)=\phi_n\circ\ldots \circ \phi_2\circ\phi_1(n)$, then $\proba{\tilde\Omega}=1$. Moreover, for every $\omega$ in $\tilde\Omega$, $m$ in $\N$, 
$$\esp{\1{M_{\tilde{\phi}\circ\phi_0(n)} \leq q_{m}} \middle |\X_{ \tilde{\phi}\circ\phi_0(n) }}(\omega)\cv{n\to+\infty}  \Phi_{0,\sigma^{2}_{P^1\otimes P^2}}(q_{m}).$$
Since $\Phi_{0,\sigma^{2}_{P^1\otimes P^2}}$ is a continuous distribution function, it can be proved that this follows
$$d_{BL}\pa{\loi{M_{\tilde{\phi}\circ\phi_0(n)}  \middle| \X_{ \tilde{\phi}\circ\phi_0(n)}},\mc{N}\pa{0, \sigma^{2}_{P^1\otimes P^2} }} \cvps{n\to+\infty} 0.$$
To conclude, we actually proved that
$$d_{BL}\pa{\loi{M_n^\rperm\!\pa{\X_n}\middle| \X_n},\mc{N}\pa{0,\sigma^{2}_{P^1 \otimes P^2}}} \cvproba{n\to+\infty} 0,$$
which, combined with \eqref{step1Perm}, leads to
$$d_{BL}\pa{ \loi{ \sqrt{n}U_{n}\!\pa{\X_n^\rperm} \middle | \X_n},   \mc{N}\pa{0,\sigma^{2}_{P^1 \otimes P^2}} }\cvproba{n\to+\infty} 0.$$

\paragraph{Third step: convergence of conditional second order moments} Recall that from \eqref{simpleU}, $U_{n}\!\pa{\X_n^\rperm} =\frac{1}{n-1} U_{n}^\rperm,$ where
\begin{eqnarray*}
U_{n}^\rperm&=&\sum_{i=1}^n \varphi\!\pa{X_i^1,X_{\rperm(i)}^2}-\frac{1}{n}\sum_{i,j=1}^n  \varphi\!\pa{X_i^1,X_j^2}\\
&=& \sum_{i,j=1}^n \pa{\1{ \rperm(i)=j}-\frac{1}{n}}  \varphi\!\pa{X_i^1,X_j^2}.
\end{eqnarray*}
Therefore,
\begin{equation}\label{dereq}
\esp{\pa{\sqrt{n}U_{n}\!\pa{\X_n^\rperm} }^2\middle| \X_n}=\frac{n^2}{(n-1)^2}\pa{ \frac{1}{n} \esp{\pa{U_{n}^\rperm}^2 \middle |\X_n}},
\end{equation}
and if $C_{i,j,k,l}=\pa{\esp{\1{ \rperm(i)=j}\1{ \rperm(k)=l}}-\frac{1}{n^2}} \varphi\!\pa{X_i^1,X_j^2}\varphi\!\pa{X_k^1,X_l^2}$,
$$\frac{1}{n} \esp{\pa{U_{n}^\rperm}^2 \middle |\X_n}=\frac{1}{n}\sum_{i,j=1}^n\sum_{k,l=1}^n C_{i,j,k,l}.$$
Firstly,
 $$\frac{1}{n}\sum_{\substack{i,j,k,l\in\{1,\ldots,n\}\\
\#\{i,j,k,l\}=4}}C_{i,j,k,l}=\frac{(n-2)(n-3)}{n^2} U_{n,1},$$
where
$$U_{n,1}=\frac{(n-4)!}{n!} \sum_{\substack{i,j,k,l\in\{1,\ldots,n\}\\
\#\{i,j,k,l\}=4}}  \varphi\!\pa{X_i^1,X_j^2}\varphi\!\pa{X_k^1,X_l^2}$$
is clearly a $U$-statistic of order $4$.
From the strong law of large numbers of H{\oe}ffding \cite{Hoeffding61suppl}, we thus have that 
$$\frac{(n-2)(n-3)}{n^2} U_{n,1}\cvps{n\to+\infty} \pa{\esp{\varphi\!\pa{X_1^1,X_2^2}}}^2.$$
Secondly,
$$\frac{1}{n}\sum_{\substack{i,j,k,l\in\{1,\ldots,n\}\\
\#\{i,j,k,l\}=3\\
i=j,i=l, j=k, \textrm{ or } k=l}}C_{i,j,k,l}= \frac{2(n-2)}{n^2} U_{n,2},$$
where
$$U_{n,2}=\frac{(n-3)!}{n!} \sum_{\substack{i,k,l\in\{1,\ldots,n\}\\
\#\{i,k,l\}=3}} \pa{ \varphi\!\pa{X_i^1,X_i^2}\varphi\!\pa{X_k^1,X_l^2}+\varphi\!\pa{X_i^1,X_l^2}\varphi\!\pa{X_k^1,X_i^2}}$$
is a $U$-statistic of order $3$ which converges almost surely, so 
$$\frac{2(n-2)}{n^2} U_{n,2}\cvps{n\to+\infty} 0.$$
Thirdly,
$$\frac{1}{n}\sum_{\substack{i,j,k,l\in\{1,\ldots,n\}\\
\#\{i,j,k,l\}=3\\
i=k, \textrm{ or } j=l}}C_{i,j,k,l}= -\frac{n(n-1)(n-2)}{n^3} U_{n,3},$$
where
$$U_{n,3}=\frac{(n-3)!}{n!} \sum_{\substack{i,k,l\in\{1,\ldots,n\}\\
\#\{i,k,l\}=3}} \pa{ \varphi\!\pa{X_i^1,X_k^2}\varphi\!\pa{X_i^1,X_l^2}+\varphi\!\pa{X_i^1,X_l^2}\varphi\!\pa{X_k^1,X_l^2}}$$
is a $U$-statistic of order $3$.
So, 
\begin{multline*}
-\frac{n(n-1)(n-2)}{n^3} U_{n,3}\\
\cvps{n\to+\infty}  -\esp{\pa{\esp{\varphi(X_1^1,X_2^2)\middle | X_1}}^2}-\esp{\pa{\esp{\varphi(X_1^1,X_2^2)\middle | X_2}}^2}.\end{multline*}
Fourthly,
$$\frac{1}{n}\sum_{\substack{i,j,k,l\in\{1,\ldots,n\}\\
\#\{i,j,k,l\}=2\\
i=j=k, i=j=l,\\ i=k=l,\textrm{ or } j=k=l}}C_{i,j,k,l}= -\frac{2(n-1)}{n^2} U_{n,4},$$
where
$$U_{n,4}=\frac{1}{n(n-1)} \sum_{1\leq i\neq j\leq n} \pa{ \varphi\!\pa{X_i^1,X_i^2}\varphi\!\pa{X_i^1,X_j^2}+\varphi\!\pa{X_i^1,X_i^2}\varphi\!\pa{X_j^1,X_i^2}}$$
is a $U$-statistic of order $2$, so $$-\frac{2(n-1)}{n^2} U_{n,4}\cvps{n\to+\infty}  0.$$
Fifthly,
$$\frac{1}{n}\sum_{\substack{i,j,k,l\in\{1,\ldots,n\}\\
\#\{i,j,k,l\}=2\\
i=j\neq k=l , \textrm{ or } i=l\neq j=k}}C_{i,j,k,l}= \frac{1}{n^2} U_{n,5},$$
where
$$U_{n,5}=\frac{1}{n(n-1)} \!\sum_{1\leq i\neq j\leq n}\!\!\! \pa{ \varphi\!\pa{X_i^1,X_i^2}\varphi\!\pa{X_j^1,X_j^2}+\varphi\!\pa{X_i^1,X_j^2}\varphi\!\pa{X_j^1,X_i^2}}$$
is a $U$-statistic of order $2$, so $$\frac{1}{n^2}U_{n,5}\cvps{n\to+\infty}  0.$$
Sixthly,
$$\frac{1}{n}\sum_{\substack{i,j,k,l\in\{1,\ldots,n\}\\
\#\{i,j,k,l\}=2\\
i=k\neq j=l}}C_{i,j,k,l}= \frac{(n-1)^2}{n^2} U_{n,6},$$
where
$$U_{n,6}=\frac{1}{n(n-1)} \sum_{1\leq i\neq j\leq n} \varphi^2\!\pa{X_i^1,X_j^2}$$
is a $U$-statistic of order $2$, so $$\frac{(n-1)^2}{n^2}U_{n,6}\cvps{n\to+\infty}  \esp{\varphi^2\!\pa{X_1^1,X_2^2}}.$$
Seventhly,
$$\frac{1}{n}\sum_{\substack{i,j,k,l\in\{1,\ldots,n\}\\
\#\{i,j,k,l\}=1}}C_{i,j,k,l}= \frac{n-1}{n^3} \sum_{i=1}^n \varphi\!\pa{X_i^1,X_i^2},$$
which almost surely tends to $0$ thanks to the strong law of large numbers.\\
By combining all these results, and the fact that
\begin{multline*}
\sigma^{2}_{P^1\otimes P^2} = \esp{\varphi^2\!\pa{X_1^1,X_2^2}} +\pa{\esp{\varphi\!\pa{X_1^1,X_2^2}}}^2\\
-\esp{\pa{\esp{\varphi(X_1^1,X_2^2)\middle | X_1}}^2}-\esp{\pa{\esp{\varphi(X_1^1,X_2^2)\middle | X_2}}^2},
\end{multline*} 
we finally obtain that
$$\frac{1}{n} \esp{\pa{U_{n}^\rperm}^2 \middle |\X_n}\cvps{n\to+\infty} \sigma^{2}_{P^1\otimes P^2},$$
and from \eqref{dereq}, we deduce that
$$\esp{\pa{\sqrt{n}U_{n}\!\pa{\X_n^\rperm} }^2\middle| \X_n}\cvps{n\to+\infty} \sigma^{2}_{P^1\otimes P^2}.$$
Since $d_{BL}\pa{ \loi{ \sqrt{n}U_{n}\!\pa{\X_n^\rperm} \middle | \X_n},   \mc{N}\pa{0,\sigma^{2}_{P^1 \otimes P^2}} }\cvproba{n\to+\infty} 0,$ this allows to conclude that
$$d_{2}\pa{ \loi{ \sqrt{n}U_{n}\!\pa{\X_n^\rperm} \middle | \X_n},   \mc{N}\pa{0,\sigma^{2}_{P^1 \otimes P^2}} }\cvproba{n\to+\infty} 0.$$

% ------------------------------------------------------------ %
\subsection{Proof of Corollary \ref{coroPerm}}
% ------------------------------------------------------------ %

Here, unlike the bootstrap approach, we only have in Theorem \ref{consistencyperm} a consistency result in probability. Thus, as for Proposition \ref{bootMonteCarlo}, we use an argument of subsequences. So let  $\phi_0:\N\to\N$ be an extraction defining a subsequence. By Theorem \ref{consistencyperm}, there exists an extraction $\phi_1$ such that $P$-a.s. in $(X_i)_i$,
\begin{equation}
\label{cvBLpsPerm}
\mathcal{L}\!\pa{\!\! \sqrt{\phi_1\!\!\circ\!\phi_0(n)}U_{\phi_1\circ\phi_0(n)}\!\!\pa{\X_{\phi_1\circ\phi_0(n)}^{\Pi_{\phi_1\circ\phi_0(n)}}\!}\! \middle | \X_{\phi_1\circ\phi_0(n)}\!}\! \cvf{n\to+\infty} \mc{N}\!\pa{\!0,\sigma^2_{P^1 \otimes P^2}\!}.
\end{equation}
In particular, applying \cite[Lemma 21.2]{vandervaartsuppl} on the event where the convergence is true, we obtain that for $\eta$ in $(0,1)$, $q_{\eta,\phi_1\circ\phi_0(n)}^\star\pa{ \X_{\phi_1\circ\phi_0(n)}}\cvps{n\to+\infty} \Phi^{-1}_{0,\sigma^2_{P^1 \otimes P^2}}(\eta),$
which ends the proof by \cite[Theorem 9.2.1]{Dudleysuppl}.

% ------------------------------------------------------------ %
\subsection{Proof of Theorem \ref{NivConsistencyPerm}}
% ------------------------------------------------------------ %

The proof of Theorem \ref{NivConsistencyPerm} for the tests of $\Gamma(q^\star)$ is very similar to the one of Theorem \ref{thniveauconsistance}, just replacing the argument of \eqref{cvqstar} by
$q_{1-\alpha,n}^\star\pa{\X_n}\cvproba{n\to+\infty}\Phi^{-1}_{0,\sigma^{2}_{P^1 \otimes P^2}}(1-\alpha),$ which is derived from Corollary \ref{coroPerm}.

Now for the tests with a Monte Carlo approximation of the quantiles, we use arguments similar to those of Proposition \ref{bootMonteCarlo}, still focusing on the upper-tailed tests of $\Gamma(q^\star_{MC})$.  We therefore aim here at proving that
\begin{equation}
\label{cvValCritPermMC}
\sqrt{n}U^{\star (\lceil (1-\alpha)(B_n+1) \rceil)} \cvproba{n\to+\infty} \Phi^{-1}_{0,\sigma^{2}_{P^1 \otimes P^2}}(1-\alpha).
\end{equation} 
Then, one can conclude as in the proof of Theorem \ref{thniveauconsistance}.

Let $F^\star_{n,\X_n}$ be the c.d.f of 
$\loi{\sqrt{n}U_{n},P_n^{\star} | \X_n}$, and let us first prove that 
\begin{equation}
\label{cvProbaSupVraieGauss}
\sup_{z\in\R}\left|F^\star_{n,\X_n}(z) - \Phi_{0,\sigma^{2}_{P^1 \otimes P^2}}(z)\right| \cvproba{n\to+\infty} 0.
\end{equation}
As Theorem \ref{consistencyperm} provides only a convergence in probability, similar arguments of subsequences as in the proof of Corollary \ref{coroPerm}, have to be used. So, let $\phi_0$ be an initial extraction and $\phi_1$ be the extraction such that \eqref{cvBLpsPerm} is satisfied. 
As convergence in the $d_{BL}$ metric is equivalent to a weak convergence (see \cite[Proposition 11.3.3]{Dudleysuppl} for instance), and as the limit is continuous, by \cite[Lemma 2.11]{vandervaartsuppl} we obtain that 
$$\sup_{z\in\R} \left|F^\star_{\phi_1\circ\phi_0(n),\X_{\phi_1\circ\phi_0(n)}}(z) - \Phi_{0,\sigma^{2}_{P^1 \otimes P^2}}(z)\right|\cvps{n\to+\infty} 0. $$
This being true for any initial extraction $\phi_0$, we obtain \eqref{cvProbaSupVraieGauss}.\\
Let $F^{\star B_n}_{n,\X_n}$ denote the empirical c.d.f of $\loi{\sqrt{n}U_{n},P_n^{\star} | \X_n}$ associated with the sample $\pa{\Pi_n^1,\ldots,\Pi_n^{B_n}}$, that is 
$$\forall z\in\R,\quad 
F^{\star B_n}_{n,\X_n}(z) = \frac{1}{B_n}\sum_{b=1}^{B_n} \1{\sqrt{n}U_{n}\!\pa{\X_n^{\Pi_n^b}}\leq z}.$$
Then, by the DKW inequality, we obtain as in the proof of Proposition \ref{bootMonteCarlo},
\begin{equation}
\label{cvProbaSupEmpVraie}
\sup_{z\in\R}\left|F^{\star B_n}_{n,\X_n}(z) - F^\star_{n,\X_n}(z)\right| \cvproba{n\to+\infty} 0.
\end{equation}
Finally, let
$$G^{\star B_n}_{n,\X_n}(z) = \frac{1}{B_n+1} \sum_{b=1}^{B_n+1} \1{\sqrt{n} U^{\star b} \leq z}.$$
Since $G^{\star B_n}_{n,\X_n}(z) = \frac{1}{B_n+1}\pa{\1{\sqrt{n}U_{n}(\X_n)\leq z} + B_n F^{\star B_n}_{n,\X_n}(z)}$, 
\begin{equation}
\label{cvpsSupEmpValCrit}
\sup_{z\in\R}\left|G_{n,\X_n}^{\star B_n}(z) - F^{\star B_n}_{n,\X_n}(z)\right|\leq \frac{2}{B_n+1} \cv{n\to+\infty} 0.
\end{equation}
Combining \eqref{cvProbaSupVraieGauss}, \eqref{cvProbaSupEmpVraie} and \eqref{cvpsSupEmpValCrit} leads to:
\begin{equation*}
\sup_{z\in\R}\left|G_{n,\X_n}^{\star B_n}(z) - \Phi_{0,\sigma^{2}_{P^1 \otimes P^2}}(z)\right| \cvproba{n\to+\infty} 0.
\end{equation*}
Since
$$\sqrt{n}U^{\star (\lceil (1-\alpha)(B_n+1) \rceil)} = \pa{G_{n,\X_n}^{\star B_n}}^{-1}(1-\alpha),$$ 
we obtain \eqref{cvValCritPermMC}, and as explained above, we conclude as in the proof of Theorem \ref{thniveauconsistance}.

%\smallskip
%%\subsection{Proof of Proposition \ref{MonteCarloniveauexact}}
%Now, let us prove the non asymptotic part of the result.
%%, let us notice that  for every $\alpha$ in $(0,1)$, $U_{n}\!\pa{\X_n}>U^{\star (\lceil (1-\alpha)(B+1) \rceil)}$ is equivalent to $p_{B,\varphi}^{\star +}\pa{\X_n}\leq \alpha$ with $p_{B,\varphi}^{\star +}\pa{\X_n}=\frac{1}{B+1}\sum_{b=1}^{B+1} \1{U^{\star b}\geq U^{\star B+1}},$  which defines a $p$-value. 
%Under $(H_0)$, from Proposition \ref{exactdistribution}, we deduce that $\X_n^{\Pi_n^1},\ldots,\X_n^{\Pi_n^B} ,\X_n$, and hence also $U^{\star 1},\ldots,U^{\star B+1}$, are exchangeable real-valued random variables. Then,
%\begin{align*}
%\mathds{P}\Big(U^{\star B+1}>&U^{\star \pa{\lceil (1-\alpha)(B+1)\rceil }}\Big)\\
%&= \proba{\sum_{b=1}^{B+1}\1{U^{\star b}<U^{\star B+1}}\geq (B+1)-\lfloor \alpha(B+1)\rfloor }\\
%&= \proba{\sum_{b=1}^{B+1}\1{U^{\star b}\geq U^{\star B+1}}\leq \lfloor \alpha(B+1)\rfloor}\\
%&= \proba{\sum_{b=1}^{B+1}\1{U^{\star b}\geq U^{\star B+1}}\leq \alpha (B+1)}.
%\end{align*}
%By exchangeability of $U^{\star 1},\ldots,U^{\star B+1}$, applying Lemma 1 in \cite{RomanoWolf}, whose proof is given in \cite{Arlot}, we finally obtain that
%$$\proba{U^{\star B+1}>U^{\star \pa{\lceil (1-\alpha)(B+1)\rceil }}}\leq \alpha,$$
%which ends the proof.

% ---------------------------------------------------------------------------------------------------- %
\section{Additional Results}
% ---------------------------------------------------------------------------------------------------- %
 
% ---------------------------------------------------------------------------------------------------- %
\subsection{About the non-degeneracy of the $U$-statistic}
\label{Apdx:dege}
% ---------------------------------------------------------------------------------------------------- %

We focus on the \emph{Linear case} with $\varphi=\varphi^w$ given by \eqref{defw}.
Define
\begin{multline*}
Z(x)=\int w(u,v) dN_{x^1}(u)dN_{x^2}(v)+ \esp{\int w(u,v) dN_{X^1}(u)dN_{X^2}(v)}\\
-\esp{\int w(u,v) dN_{x^1}(u)dN_{X^2}(v)}-\esp{\int w(u,v) dN_{X^1}(u)dN_{x^2}(v)}.
\end{multline*}

Recall that in this case, degeneracy is equivalent to stating that for $X=(X^1,X^2)$ with distribution $P^1\otimes P^2$, 
$Z(X)$ is a random variable which is almost surely null under $(H_0)$. Since $\esp{Z(X)}=0$, $Z(X)=0$ a.s. is equivalent to $\var{Z(X)}=0.$ Here we provide a computation of $\var{Z(X)}$.

Let us introduce $dM_1^{[1]}(u)$ and $dM_2^{[1]}(v)$ the mean measures of respectively $X^1$ with distribution $P^1$ and $X^2$ of distribution $P^2$ \cite[Chapter 5]{DVsuppl}, then one can rewrite
\begin{multline*}
Z(X)=\int w(u,v) dN_{X^1}(u)dN_{X^2}(v)+ \int w(u,v) dM_{1}^{[1]}(u)dM_{2}^{[1]}(v)\\-\int w(u,v) dN_{X^1}(u)dM_{2}^{[1]}(v)-\int w(u,v) dM_{1}^{[1]}(u)dN_{X^2}(v).
\end{multline*}

\noindent Therefore, $\esp{Z(X)}=0$, and
\begin{eqnarray*}
\var{Z(X)}&=&\esp{Z(X)^2}\\
&=&\int_{[0,1]^4} w(u,v) w(s,t) \esp{dN_{X^1}(u)dN_{X^1}(s)}\esp{dN_{X^2}(v)dN_{X^2}(t)} \\
&&- \int_{[0,1]^4} w(u,v) w(s,t) \esp{dN_{X^1}(u)dN_{X^1}(s)}dM_{2}^{[1]}(v)dM_{2}^{[1]}(t)\\
&&-\int_{[0,1]^4} w(u,v) w(s,t)  dM_{1}^{[1]}(u)dM_{1}^{[1]}(s)\esp{dN_{X^2}(v)dN_{X^2}(t)} \\
&&+ \int_{[0,1]^4} w(u,v) w(s,t)  dM_{1}^{[1]}(u)dM_{1}^{[1]}(s)dM_{2}^{[1]}(v)dM_{2}^{[1]}(t).
\end{eqnarray*}

By assuming that $\# X^1$ (resp. $\#X^2$) has second order moment, (see also Section \ref{ma} for comment on this assumption), one can introduce the second factorial moment measure associated with $X^1$ (resp. $X^2$), and denoted by $dM_1^{[2]}(u,s)$  (resp. $dM_2^{[2]}(v,t)$). Then straightforward computations show that
\begin{align*}
\ds{V}\textrm{ar}&(Z(X))=\int_{[0,1]^2} w(u,v)^2 dM_{1}^{[1]}(u) dM_{2}^{[1]}(v) \\
&+ \int_{[0,1]^3} w(u,v) w(u,t) dM_{1}^{[1]}(u) \pa{dM_{2}^{[2]}(v,t)-dM_{2}^{[1]}(v)dM_{2}^{[1]}(t)}\\
&+\int_{[0,1]^3} w(u,v) w(s,v)  \pa{dM_{1}^{[2]}(u,s)-dM_{1}^{[1]}(u)dM_{1}^{[1]}(s)} dM_{2}^{[1]}(v)\\
&+\int_{[0,1]^4} w(u,v) w(s,t)  \pa{dM_{1}^{[2]}(u,s)-dM_{1}^{[1]}(u)dM_{1}^{[1]}(s)}\times\\
&\quad\pa{dM_{2}^{[2]}(v,t)-dM_{2}^{[1]}(v)dM_{2}^{[1]}(t)}.
\end{align*}

In particular, for Poisson processes, $dM^{[2]}(u,s)=dM^{[1]}(u)dM^{[1]}(s)$ and 
$$\var{Z(X)}=\int_{[0,1]^2} w(u,v)^2  dM_{1}^{[1]}(u) dM_{2}^{[1]}(v) >0,$$
as soon as the Poisson processes have non zero intensities since for $j=1,2$, $dM_j^{[1]}(u)=\lambda_j(u)du$, with $\lambda_j$ the intensity of $X^j$.

% ---------------------------------------------------------------------------------------------------- %
\subsection{About the empirical centering assumption} 
\label{Hyp_CentrEmp}
% ---------------------------------------------------------------------------------------------------- %
Recall that \smallskip\\
$\pa{\mc{A}_{Cent}^*}\quad \textrm{\begin{tabular}{|l} For $x_1=(x_1^1,x_1^2),\dots, x_n=(x_n^1,x_n^2)$ in $\calX^2$,\\
$\sum_{i_1,i_2,i'_1,i_2'=1}^n h\pa{\pa{x_{i_1}^1,x_{i_2}^2},\pa{x_{i'_1}^1,x_{i'_2}^2}}=0.$\end{tabular}}$

\paragraph{}
On the one hand, in the \emph{Linear case}, that is if $h=h_\varphi$, then for $n\geq 1$ and for $x_1=(x_1^1,x_1^2), \dots, x_n=(x_n^1,x_n^2)$ in $\calX^2$, 
\begin{align*}
\sum_{i,i',j,j'=1}^n \!\! & h\pa{\pa{x_{i}^1,x_{i'}^2},\pa{x_{j}^1,x_{j'}^2}} \\
&=\frac 1 2 \sum_{i,i',j,j'=1}^n \pa{\varphi\left(x_i^1,x_{i'}^2\right)+\varphi\left(x_j^1,x_{j'}^2\right)-\varphi\left(x_i^1,x_{j'}^2\right)-\varphi\left(x_j^1,x_{i'}^2\right)} \\
&= \frac{n^2}{2}\!\pa{\!\sum_{i,i'=1}^n\!\varphi\!\left(x_i^1,x_{i'}^2\right) +\!\!\sum_{j,j'=1}^n\!\varphi\!\left(x_j^1,x_{j'}^2\right) -
\!\!\sum_{i,j'=1}^n\!\varphi\!\left(x_i^1,x_{j'}^2\right) -\!\!\sum_{j,i'=1}^n\!\varphi\!\left(x_j^1,x_{i'}^2\right)\!}\\
&= 0.
\end{align*}
So $(\mc{A}_{Cent}^*)$ is immediately satisfied in the \emph{Linear case}. 

\paragraph{}
On the other hand, $(\mc{A}_{Cent}^*)$ does not imply that $h$ is of the form $h_\varphi$. \\
Indeed, consider 
$$h\pa{\pa{x^1,x^2},\pa{y^1,y^2}}=\# x^1\cdot \# x^2\cdot \# y^1\cdot \# y^2\cro{\pa{\# x^1-\# y^1}\pa{\# x^2-\# y^2}}.$$

\begin{itemize}
\item The kernel $h$ is obviously symmetric.
%\begin{align*}
%h\!\big(\pa{y^1,y^2}&,\pa{x^1,x^2}\big)\\
%&= \# y^1\cdot \# y^2\cdot \# x^1\cdot \# x^2\cro{\pa{\# y^1-\# x^1}\pa{\# y^2-\# x^2}} \\
%&= \# x^1\cdot \# x^2\cdot \# y^1\cdot \# y^2\cro{ (-1)^2\cdot \pa{\# x^1-\# y^1}\pa{\# x^2-\# y^2}} \\
%&= h\pa{\pa{x^1,x^2},\pa{y^1,y^2}}.
%\end{align*}

\item The kernel $h$ satisfies $(\mc{A}_{Cent}^*)$. Indeed, let
$$f\pa{x^1,y^1}=\# x^1\cdot \# y^1\pa{\# x^1-\# y^1}.$$
First, notice that $f\pa{x^1,x^1} = 0$ and $f\pa{x^1,y^1}=-f\pa{y^1,x^1}$.

Moreover, $h\pa{\pa{x^1,x^2},\pa{y^1,y^2}} = f\pa{x^1,y^1} f\pa{x^2,y^2}$. 
Thus  
\begin{align*}
\sum_{i,i',j,j'=1}^n\!\!& h\pa{\pa{x_{i}^1,x_{i'}^2},\pa{x_{j}^1,x_{j'}^2}} \\
&= \sum_{i,i',j,j'=1}^n f\pa{x_i^1,x_j^1}f\pa{x_{i'}^2,x_{j'}^2}\\
&= \pa{\sum_{i,j=1}^n f\pa{x_i^1,x_j^1}}\pa{\sum_{i',j'=1}^n f\pa{x_{i'}^2,x_{j'}^2}}\\
&= \pa{\sum_{i=1}^n \underbrace{f\!\pa{x_i^1,x_i^1}}_{0} + \!\!\!\!\sum_{1\leq i<j\leq n} \!\!\!\underbrace{f\!\pa{x_i^1,x_j^1}\!+\!f\!\pa{x_j^1,x_i^1}}_0}\!\pa{\sum_{i',j'=1}^n \!\!f\!\pa{x_{i'}^2,x_{j'}^2}}\\
&= 0,
\end{align*}
and thus $(\mc{A}_{Cent}^*)$ is satisfied by $h$. 

\item The kernel $h$ cannot be written as an $h_\varphi$. \\
On the one hand, first notice that for any $\varphi:\calX^2 \to \R$, the difference 
$$D_{h_\varphi}:= h_\varphi\!\pa{\pa{x^1,x^2},\pa{y^1,y^2}} - h_\varphi\!\pa{\pa{\tilde x^1,x^2},\pa{y^1,y^2}}$$ 
does not depend on $y^1$. 
Indeed, 
\begin{eqnarray*}
D_{h_\varphi} 
&=&\frac 1 2\Big( \varphi\left(x^1,x^2\right) + \varphi\left(y^1,y^2\right) - \varphi\left(x^1,y^2\right) - \varphi\left(y^1,x^2\right) \\ 
&& - \varphi\left(\tilde x^1,x^2\right) - \varphi\left(y^1,y^2\right) + \varphi\left(\tilde x^1,y^2\right) + \varphi\left(y^1,x^2\right)\Big) \\
&=&\frac 1 2\pa{ \varphi\left(x^1,x^2\right) - \varphi\left(\tilde x^1,x^2\right) + \varphi\left(\tilde x^1,y^2\right) -\varphi\left(x^1,y^2\right)}.
\end{eqnarray*}

On the other hand, for the kernel $h$ introduced above, the difference $D_h$ does depend on $y^1$. 
Indeed 
\begin{eqnarray*}
D_h
&=& h\pa{\pa{x^1,x^2},\pa{y^1,y^2}} - h\pa{\pa{\tilde x^1,x^2},\pa{y^1,y^2}} \\
&=& \# x^2\cdot \# y^1\cdot \# y^2 \pa{\# x^2-\# y^2}\times\\
&&\cro{\# x^1\cdot \pa{\# x^1-\# y^1}-\# \tilde x^1\cdot \pa{\# \tilde x^1-\# y^1}},
\end{eqnarray*}
and if for instance, $\# x^1=\# y^2 = 1$ and $\# \tilde x^1 = \# x^2=2$, then 
$$ D_h = 2\# y^1\cro{-3 + \# y^1},$$
which clearly depends on $y^1$.

So finally, there does not exist any $\varphi$ such that $h=h_\varphi$. 

\end{itemize}

\end{document}